\documentclass{lmcs}
\pdfoutput=1

\usepackage{lastpage}
\lmcsdoi{20}{2}{17}
\lmcsheading{}{\pageref{LastPage}}{}{}%
{Jun.~08,~2021}{Jun.~24,~2024}{}

\usepackage[utf8]{inputenc}

\keywords{Computable Analysis, Reliable Numerics, Imperative Programming, Turing-Completeness, Decidability, Hoare Logic, Kleene Logic}
\usepackage[unicode]{hyperref}
\usepackage{xspace}
\usepackage{amssymb}
\usepackage{amsthm}
\usepackage{stmaryrd} 
\usepackage{bm}
\usepackage{proof}
\usepackage{mathpartir}
\usepackage{mathtools}
\usepackage{listings}
\lstdefinelanguage{typetwo}{
    basicstyle=\footnotesize,
    breaklines=true,
    breakatwhitespace=false,
  keywords={input, if, then, else, end, while, repeat, until, break, with, iterate, return, goto, line, or, diverge},
  keywordstyle=\color{black}\bfseries,
  identifierstyle=\color{black},
  sensitive=false,
  comment=[l]{//},
  morecomment=[s]{/*}{*/},
  commentstyle=\color{purple}\ttfamily,
  stringstyle=\color{red}\ttfamily,
  morestring=[b]',
  morestring=[b]",
  tabsize=4
}

\usepackage{algorithm}
\usepackage{algpseudocode}
\newcommand{\tassert}[1]{\ensuremath{\pmb{\big\{}#1\pmb{\big\}}}}

\newcommand{\asser}[1]{ \;\pmb{\{\vphantom{\frac{1}{2}}} #1 \pmb{\vphantom{\frac{1}{2}}\}\;} }
\newcommand{\sem}[1]{\eval{#1}}
\newcommand{\trans}[1]{\evalP{#1}}

\newcommand{\naive}{na\"{\i}ve\xspace}
\newcommand{\naively}{na\"{\i}vely\xspace}

\newcommand{\IR}{\mathbb{R}}

\newcommand{\IK}{\mathbb{K}}
\newcommand{\Kleene}{\IK}

\newcommand{\IZ}{\mathbb{Z}}
\newcommand{\IN}{\mathbb{N}}

\newcommand{\calF}{\mathcal{F}}
\newcommand{\calG}{\mathcal{G}}

\newcommand{\calV}{\mathcal{V}}

\newcommand{\REAL}{\mathsf{R}} 
\newcommand{\INTEGER}{\mathsf{Z}}
\newcommand{\LOGIC}{\mathsf{K}}

\newcommand{\dreal}{\REAL}
\newcommand{\dint}{\INTEGER}
\newcommand{\dbool}{\LOGIC}

\newcommand{\mycond}[3]{#1\;?\;#2{\;:\;}#3}
\newcommand{\ctt}{\mathit{true}}
\newcommand{\cff}{\mathit{false}}
\newcommand{\cuu}{\mathit{unknown}}

\newcommand{\mychoose}{\mathsf{choose}}
\newcommand{\mult}{\times} 
\newcommand{\TRUE}{\ctt}
\newcommand{\FALSE}{\cff}
\newcommand{\UNDEF}{\cuu}
\newcommand{\true}{\TRUE}

\newcommand{\ttrue}{\TRUE}
\newcommand{\tfalse}{\FALSE}
\newcommand{\tundef}{\UNDEF}
\newcommand{\elt}{\lessdot}

\newcommand{\WHILE}{\textup{\textbf{while~}}}
\newcommand{\LET}{\textup{\textbf{let~}}}
\newcommand{\ENDWHILE}{\textup{\textbf{end}}}

\newcommand{\FOR}{\textup{\textbf{for}}}
\newcommand{\TO}{\textup{\textbf{to}}}
\newcommand{\DO}{\textup{\textbf{do}}}
\newcommand{\ENDFOR}{\textup{\textbf{end}}}
\newcommand{\IF}{\textup{\textbf{if}}}
\newcommand{\THEN}{\textup{\textbf{then}}}
\newcommand{\ELSE}{\textup{\textbf{else}}}
\newcommand{\ENDIF}{\textup{\textbf{end}}}

\newcommand{\prog}{P}

\newcommand{\eskip}{\textup{\textbf{skip}}}
\newcommand{\eif}{\textup{\textbf{if}}}
\newcommand{\ethen}{\textup{\textbf{then}}}
\newcommand{\eelse}{\textup{\textbf{else}}}
\newcommand{\elet}{\textup{\textbf{let}}}
\newcommand{\einput}{\textup{\textbf{input}}}
\newcommand{\ereturn}{\mathbf{return}}
\newcommand{\eendif}{\textup{\textbf{end}}}
\newcommand{\ewhile}{\textup{\textbf{while}}}
\newcommand{\edo}{\textup{\textbf{do}}}
\newcommand{\eendwhile}{\textup{\textbf{end}}}

\newcommand{\specification}{specification\xspace}

\newcommand{\ieq}{=}

\newcommand{\eval}[1]{\left\llbracket #1 \right\rrbracket }

\newcommand{\evalP}[1]{\llparenthesis #1 \rrparenthesis }

\newcommand{\dom}{\operatorname{dom}}

\newcommand{\sign}{\operatorname{sign}}

\newcommand{\id}{\operatorname{id}}

\newcommand{\abs}{\operatorname{abs}}

\newcommand{\ERC}{\textup{ERC}\xspace}
\newcommand{\Structure}{\mathcal{S}}
\newcommand{\Theory}{\mathcal{T}}
\newcommand{\Logic}{\mathcal{L}}

\newcommand{\toto}{\rightrightarrows}

\newcommand{\dfeq}{\coloneqq}

\newcommand{\pcolon}{\mathpunct{\,:\subseteq}}

\newcommand{\swap}{\textup{\textbf{swap}}}

\newcommand{\GL}{\operatorname{GL}}

\newcommand{\iExp}{\operatorname{\mathtt{iExp}}}
\newcommand{\pExp}{\operatorname{\mathtt{pExp}}}
\newcommand{\ppExp}{\operatorname{\mathtt{posExp}}}
\newcommand{\pPivot}{\operatorname{\mathtt{Pivot}}}
\newcommand{\pDet}{\operatorname{\mathtt{Det}}}
\newcommand{\nround}{\operatorname{\mathtt{Round}}}
\newcommand{\fround}{\operatorname{\mathtt{BinRound}}}
\newcommand{\HeronSqrt}{\operatorname{\mathtt{HeronSqrt}}}
\newcommand{\Trisection}{\operatorname{\mathtt{Trisection}}}
\newcommand{\uniq}{\mathrm{uniq}}
\newcommand{\cont}{\mathrm{cont}}

\newcommand{\Root}{\mathrm{Root}}
\newcommand{\Det}{\mathrm{Det}}

\newcommand{\Pivot}{\operatorname{\mathrm{Pivot}}}

\newcommand{\Round}{\operatorname{\mathrm{Round}}}  

\newcommand{\pdom}[1]{\mathbb{P}(#1)}

\makeatletter
\def\smallunderbrace#1{\mathop{\vtop{\m@th\ialign{##\crcr
   $\hfil\displaystyle{#1}\hfil$\crcr
   \noalign{\kern3\p@\nointerlineskip}%
   \tiny\upbracefill\crcr\noalign{\kern3\p@}}}}\limits}
\makeatother

\newcommand{\COMMENTED}[1]{}

\newcommand{\subst}[2]{#2 / #1}

\newcommand{\all}[1]{\forall #1.\ }
\newcommand{\some}[1]{\exists #1.\ }
\newcommand{\usome}[1]{\exists! #1.\ }

\newcommand{\bop}{\mathbin{\odot}}
\newcommand{\uop}{\mathop{\star}}
\newcommand{\bopty}{\mathbf{BinTy}}
\newcommand{\uopty}{\mathbf{UnTy}}

\newcommand{\kneg}{\mathop{\hat\neg}}
\newcommand{\klor}{\mathbin{\hat\lor}}
\newcommand{\kland}{\mathbin{\hat\land}}

\newcommand{\semcond}{\mathrm{Cond}}
\newcommand{\semkond}{\mathrm{Kond}}
\newcommand{\semproj}{\mathrm{proj}}

\newcommand{\semupdate}{\mathrm{update}}
\newcommand{\semchoose}{\mathrm{choose}}

\newcommand{\fcomp}{\mathbin{\circ}}

\newcommand{\pset}[1]{\mathcal{P}( #1 )}
\newcommand{\blim}[1]{\mathbf{input}\;}
\newcommand{\blimreturn}[2]{\ereturn\; #1\; \mathbf{as} \; #2 \to -\infty}

\algblockdefx{Prog}{EndProg}[1]{\textbf{input} #1}[1]{\textbf{return} #1}
\algblockdefx{Limit}{EndLimit}[2]{$\blim{#1} #2$}[2]{$\blimreturn{#1}{#2}$}
\algblockdefx{Mywhile}{EndMywhile}[1]{$\ewhile$ #1 $\edo$}{\eendwhile}
\algblockdefx{Myfor}{EndMyfor}[1]{$\FOR$ #1  $\DO$}{\ENDFOR}
\algblockdefx{Myscope}{EndMyscope}{$\{$}{$\}$}
\algnewcommand{\LeftComment}[1]{\State {\color{darkgray} // #1}}
\algnewcommand{\RightComment}[1]{\hfill {\color{darkgray} // #1}}

\newcommand{\ileq}{\leq}

\newcommand{\lfp}{\mathbf{lfp}}
\newcommand{\wop}{\mathbf{W}}

\newcommand{\mydom}[1]{\mathrm{dom}(#1)}

\newcommand{\otype}[1]{\langle #1\rangle}

\newcommand{\myabs}[1]{\lvert #1 \rvert}
\newcommand{\myappend}{\#}
\newcommand{\pto}{\rightharpoonup}

\newcommand{\cround}{\text{round}}
\newcommand{\creal}{\text{real}}

\newcommand{\abbeq}{:\equiv}

\newcommand{\inlinedef}[1]{\textit{\textbf{#1}}}
\newcommand{\myprec}[1]{{2}^{#1}}

\title[Exact Real Computation]{Semantics, Specification Logic, and Hoare Logic of Exact Real Computation}

\author[Park]{Sewon Park
\lmcsorcid{0000-0002-6443-2617}}[a]%
\address{School of Computing, KAIST, Daejeon, Republic of Korea} 
\email{swelite@kaist.ac.kr, ziegler@kaist.ac.kr} 

\author[Brau{\ss}e]{Franz Brau{\ss}e
\lmcsorcid{0000-0002-2386-7489}}[b]
\address{University of Manchester}
\email{franz.brausse@manchester.ac.uk}

\author[Collins]{Pieter Collins
\lmcsorcid{0000-0002-8896-9603}}[c]
\address{Maastricht University, Maastricht, The Netherlands}
\email{pieter.collins@maastrichtuniversity.nl}

\author[Kim]{SunYoung Kim
}[d]
\address{Yonsei University, Seoul, Republic of Korea}
\email{sy831@yonsei.ac.kr}

\author[Kone\v{c}n\'{y}]{Michal Kone\v{c}n\'{y}
\lmcsorcid{0000-0003-2374-9017}}[e]%
\address{Aston University, Birmingham, UK}
\email{m.konecny@aston.ac.uk}

\author[Lee]{Gyesik Lee
\lmcsorcid{0009-0004-6589-2684} 
}[f]
\address{Hankyong National University, Anseong, Republic of Korea}
\email{gslee@hknu.ac.kr}

\author[M\"{u}ller]{Norbert M\"{u}ller
\lmcsorcid{0000-0003-3684-3029}
}[g]
\address{Universit\"{a}t Trier, Trier, Germany}
\email{mueller@uni-trier.de}

\author[Neumann]{Eike Neumann
}[h]
\address{Swansea University, Swansea, UK}
\email{neumaef1@gmail.com}

\author[Preining]{Norbert Preining
}[i]
\address{arXiv, Cornell Tech, Komatsu, Japan}
\email{norbert@preining.info}

\author[Ziegler]{Martin Ziegler
\lmcsorcid{0000-0001-6734-7875}}[a]

\allowdisplaybreaks
\begin{document}

\begin{abstract}
We propose a simple imperative programming language, \ERC, that features arbitrary real numbers as primitive data type, exactly.
Equipped with a denotational semantics, \ERC provides a formal programming language-theoretic foundation to the algorithmic processing of real numbers.
In order to capture multi-valuedness, which is well-known to be essential to real number computation, we use a Plotkin powerdomain and make our programming language semantics computable and complete:
all and only real functions computable in computable analysis can be realized in \ERC.
The base programming language supports real arithmetic as well as implicit limits; expansions support additional primitive operations (such as a user-defined exponential function).
By restricting integers to Presburger arithmetic and real coercion to the `precision' embedding $\mathbb{Z}\ni p\mapsto 2^p\in\mathbb{R}$, 
we arrive at a first-order theory which we prove to be decidable and model-complete.
Based on said logic as specification language for preconditions and postconditions, 
we extend Hoare logic to a sound (w.r.t. the denotational semantics) and expressive system for deriving correct \emph{total} correctness specifications. 
Various examples demonstrate the practicality and convenience of our language and the extended Hoare logic.
\end{abstract} 

\maketitle



\begin{center}
\tableofcontents
\end{center}

\section{Introduction}
Real numbers in computable analysis are expressed and manipulated exactly without rounding errors by using infinite representations  \cite{Wei00,PR89}.
For example, a real number $x$ is represented by an infinite sequence of integers $\phi_n$ such that $|x - \phi_n \cdot 2^{-n}| \leq 2^{-n}$ holds for all $n$. 
When we have \emph{correctly} computed a real number $x$ under this framework, we automatically get access to the real number $x$ via any arbitrarily accurate precision. This property makes this approach promising for application domains where high-precision numerical results are required.
Several implementations substantiate this approach of exact real number computation:
AERN~\cite{Aern}, Ariadne~\cite{CollinsAriadne}, Core2~\cite{core2}, 
Cdar~\cite{BlanckCdar}, iRRAM~\cite{iRRAM}, and realLib~\cite{Lam07} to name a few. 
They provide their own data types of real numbers to programmers keeping complicated infinite representations in their back-ends. 
Abstracting away such details, the programmers can deal with real numbers as they were abstract algebraic entities closely resembling the classical mathematical structure of real numbers with the absence of rounding errors.
According to this observation, it is expected that reasoning about programs' behaviours in this approach becomes more intuitive;
cf. \cite{KettnerMPSY08,10.1007/978-3-540-71070-7_2,boldo2009,boldo2011flocq}.

Imperative programming is a ubiquitous style of expressing computation that is widely adopted in numerical computations admitting rich theories of precondition-postcondition-based program verification methodologies \cite{10.1145/363235.363259,apt2019fifty,reynolds2002separation,nanevski2006polymorphism} and their successful applications such as  \cite{10.1007/978-3-642-37036-6_8}.
Hence, to conveniently achieve provably correct arbitrary accurate numerical computations, it is desired to extend them with exact real number computation.
To make such an extension rigorous, we first acquire a formal programming language with exact real computation on which the extended verification methodologies are based.
However, extending an existing formal programming language with exact real number computation, which is a seemingly simple quest, fails when it is done \naively with the classical structure of real numbers as it is not identical, though resembles much thanks to exact operations, the computational structure of real numbers \cite[Section~3]{Her99a}. For example, the real number comparison $x < y$ as a binary mapping is only partially computable causing non-termination when $x = y$ \cite[Theorem~4.1.16]{Wei00}. 
This property makes multi-valued functions and nondeterminism essential notions in exact real computation \cite{Luckhardt}. 
We explain this in more detail later in Section~\ref{ss:kleenean}.

The aim of this paper is to devise an imperative programming language $\ERC$
for exact real number computation that correctly deals with the partial and nondeterministic nature in the following sense.
We formalize $\ERC$'s formal domain-theoretic denotational semantics in such a way that precisely the computable partial real functions and integer multivalued functions in computable analysis can be defined in $\ERC$.
Based on the formal semantics, we devise a sound Hoare-style logic of $\ERC$ specifications. 
The language $\ERC$ is expressive enough to model core fragments of existing software and simple enough to stand as a design guideline for future exact real number computation software developments which wish to adopt the theory developed in this paper.
For this purpose, we develop our language and theory relative to primitive operator extensions such that other frameworks providing richer sets of primitive operators can still benefit from $\ERC$ to analyze their real number computations.

\subsection{Related Works}
In terms of formalizing programming languages for exact real number computation, this paper is deeply influenced by the following related works.
Real PCF \cite{Esc96} is a pure functional language that is an extension of PCF \cite{plotkin1977lcf} with real numbers, computable real number operations, and a parallel conditional construct for nondeterminism. 
Its application in expressing integration demonstrated the usefulness of formalizing exact real number computation \cite{edalat2000integration}.
In \cite{marcial2007semantics}, a domain-theoretic denotational semantics of a variant language is defined but using the Hoare powerdomain which cannot express termination.
There, it is observed that nondeterminism in computable analysis becomes continuous only with the Hoare powerdomain.
In this paper, we avoid the continuity problem and use the Plotkin powerdomain instead such that we do not need to have a separate tool for reasoning programs' termination. 
Later in Remark~\ref{r:hoare}, we discuss more about using a different powerdomain.

When we focus on imperative programming languages for exact real number computation, 
the Blum-Cucker-Shub-Smale model \cite{BCSS} aka real RAM 
\cite[\S1.4]{PS85} captures
real numbers as algebraic structure:
It supposes tests as total, and renders transcendental functions uncomputable.
The feasible real RAM \cite{BH98} is a computable variant of real RAM where the infeasible total real comparison test \cite{Bra03f} is replaced with a nondeterministic approximation to it. Its operational semantics influenced the design of iRRAM.
Amongst further related works including \cite{era,10.1007/978-3-642-37075-5_22,10.1007/978-3-319-08918-8_15}, \textbf{WhileCC} \cite{TZ04}, a simple imperative language extended with the computational structure of real numbers and a (countable) nondeterminism construct, which is suggested with its algebraic operational semantics, is close to our design.
We proceed from there, suggest a different set of primitives, and further devise domain-theoretic computable denotational semantics, a specification language, and reasoning principles.

\subsection{Kleenean-valued Comparison and Nondeterminism}
\label{ss:kleenean}

An inevitable side-effect of using infinite representations for exact operations 
is that the ordinary order relation of real numbers as a Boolean-valued function becomes partial.
Comparing $x < y$ does not terminate when $x = y$ whichever specific representation is used \cite[Theorem~4.1.16]{Wei00}. 
Instead of making an entire procedure fail when the same numbers are compared, some frameworks, including AERN and iRRAM,  introduce the concept of Kleenean $\IK \dfeq \{\ctt, \cff, \cuu\}$, a lazy extension of Boolean. The third value $\cuu$ denotes an explicit state of non-termination that is delayed until tested.
Consider the following Kleenean-valued comparison
\begin{equation}\label{eqn:elt}
x \elt y := \begin{cases}
\ctt &\text{if } x< y,\\
\cff &\text{if } y< x,\\
\cuu &\text{otherwise,}
\end{cases}
\end{equation}
as a replacement for the ordinary Boolean-valued comparison test.
When $x = y$, the comparison $x \elt y$ evaluates to $\cuu$ instead of diverging immediately. 
Though $\cuu$ is an admissible value that can be assigned to Kleenean-typed variables when later it is required to test if the value is $\ctt$ or $\cff$, for example as the condition of a loop statement, 
it then causes non-termination.

A remark worth making here is that the third truth value $\cuu$ is only denotational in the sense that an expression having $\cuu$ value cannot be observed during a computation.
Though we can reason that $x \elt y$ must be $\cuu$ when we know that $x = y$, we cannot observe $x \elt y$ having its value $\cuu$ and conclude $x = y$. In other words, we cannot write a program that branches differently on the three values of $\IK$.

Due to the non-termination in comparisons, multi-valued functions and their nondeterministic computations become essential \cite{Luckhardt} to compose a total but nondeterministic procedure from possibly non-terminating sub-procedures.
Suppose we want to approximate $x < y$ with a tolerance factor $2^{-n}$ and consider
the two Kleenean-valued comparisons $x \elt y + 2^{-n}$ and $y \elt x + 2^{-n}$.
We can test their truths in parallel until one turns out to be $\ctt$.
Even in the case when one of them is $\cuu$, hence the test faces non-termination, it is promised that the other comparison evaluates to $\ctt$. 
Therefore, the parallel procedure can terminate safely and inform us if $x < y + 2^{-n}$ or $y < x + 2^{-n}$.
When the two $\ctt$ conditions overlap $y - 2^{-n} < x < y + 2^{-n}$ it becomes nondeterministic as it will depend on how the real numbers and the parallel procedure are implemented which we abstract away.
We model this nondeterministic computation by a multi-valued function whose value is $\{\ctt, \cff\}$ at such inputs denoting the set of possible values of the nondeterministic computation.

As nondeterminism becomes essential, most exact real computation frameworks provide their users with various primitive constructs for nondeterminism.
In iRRAM, an operator $\mychoose$ is provided that receives multiple Kleenean expressions and returns the index (counting from $0$) of \emph{some} expression which evaluates to $\ctt$.
See that the above example can be expressed by \[\mychoose(x \elt y + 2^{-n}, y \elt x + 2^{-n}) \ieq 0.\] 

\subsection{Main Contributions and Overviews}

In this paper, we define an imperative programming language \emph{Exact Real Computation} (\ERC) providing three primitive data types $\LOGIC$, $\INTEGER$, and $\REAL$ that denote the Kleene three-valued logic $\IK = \{\ctt, \cff, \cuu\}$, the set of integers $\IZ$, and the set of real numbers $\IR$, respectively.
It is a simple imperative language with a strict hierarchy posed between its term language and command language where only the latter contains unbounded loops.
Our language is carefully designed to ensure the computability (in the sense of computable analysis) of its domain-theoretic denotational semantics by adopting the Kleenean typed comparison $\elt$ for real numbers replacing the ordinary total comparison and offering the nondeterminism construction $\mychoose(b_0,b_1,\ldots)$. 
We define the denotations of terms and commands to be functions to Plotkin powerdomains \cite{DBLP:journals/siamcomp/Plotkin76} to model nondeterministic computations.
The definitions of computable partial (multi-valued) functions are reviewed in Section~\ref{s:preliminaries}.

Construction of limits, which is an essential feature that makes exact real computation more expressive than algebraic computations \cite{Bra03f}, is available in a restricted way in \ERC.
Instead of introducing an explicit operator for computing limits, which then requires the language to support constructions of nontrivial infinite sequences within its term language, \ERC supports limits to be computed at the level of ``programs'', another layer over commands representing (multi-valued) functions. 
A real-valued program is provided with access to an integer variable $p$ such that when the program computes (possibly multi-valued) $2^{p}$ approximations to a real number, the single-valued real number, which is the limit of the program as $p \to -\infty$, becomes the function value that of the program denotes.
(A similar approach of equipping limits can be found in \cite[Section~9]{TZ99} and \cite[Definition~4.5.1]{TZ04}.) 
Hence, in $\ERC$, commands compose either an integer program expressing a partial integer multi-valued function or a real program expressing a partial real function where limit operations are taken by default. 
Later in \autoref{r:realmult}, we discuss more about this design choice. 

$\ERC$ is defined relative to primitive operator extensions such that it can be extended easily to other languages that provide further primitive operators.
For each finite sets $\calF$ of computable partial real functions and $\calG$ of computable partial (finite) integer multi-valued functions, $\ERC(\calF, \calG)$ is defined. 
We simply write $\ERC$ if the underlying extension sets are obvious or irrelevant. 

The base language without any extension $\ERC_0 := \ERC(\emptyset, \emptyset)$ does not provide integer multiplications or the \naive coercion $\IZ \ni z \mapsto z \in \IR$ as its term constructs.
Instead, it provides the \emph{precision embedding} $\imath : \IZ \ni p \to 2^p \in \IR$ as a way to convert an integer to real.
Although this looks like a bit of restriction at first glance, since integer multiplication is expressible using a loop in our command language, we claim that our programming language itself is still expressive by showing its Turing-completeness in Theorem~\ref{t:Complete} and various examples in Section~\ref{s:Programming}.
This restriction and the strict hierarchy between \ERC's term language and command language are motivated by the 
specification language of \ERC which for $\ERC_0$'s terms is adequate and admits a decidable theory.

We introduce a first-order language over the three types (sorts) expanded with $(\calF, \calG)$ as the specification language of \ERC.
We show that the first-order language can define the semantics of our nondeterministic term language.
Moreover, for the default language $\ERC_0$,
whose design is influenced by the structure suggested in \cite{Dries1986}, 
the theory of the specification language
becomes model-complete and decidable with regard to the canonical interpretation.
Nonetheless, due to the fundamental trade-off, the specification language for $\ERC_0$ 
is not expressive enough to define the semantics of the command language.
\begin{rem}
\label{r:Definability}
Consider the following three desirable features:
\begin{enumerate}
\item[i)] a programming language being Turing-complete over integers, and thus able to realize an algorithm whose termination is co-r.e. hard
\item[ii)] a specification language sufficiently rich to express the termination of said algorithm
\item[iii)] a sound and complete r.e. deductive system of said specification language.
\end{enumerate}
Obviously, not all three are simultaneously feasible.
For instance integer WHILE programs / Peano arithmetic satisfy (i) / (ii), but not (iii) \cite[\S6]{Cook78}.
\end{rem}
We propose \emph{choosing} (i) and (iii) over (ii) as the base design aiming for further developments of automatic reasoning of exact real computations using the decidability property. 
However, since our language is defined relative to primitive operator extensions, those who are not happy with this trade-off can take the other direction by adding, for example, the integer multiplication or the \naive coercion $\IZ \ni z \mapsto z \in \IR$ in $(\calF, \calG)$. 
The other parts of this work remain valid.

Based on the specification language, we extend the classical Hoare logic for $\ERC$ commands' total 
correctness 
\[
\tassert{\phi}  \;S\; \tassert{\psi}
\]
that says for any initial state satisfying the precondition $\phi$,
all possible executions of the command $S$ regarding the nondeterminism in \ERC 
terminate and each yields a state satisfying the postcondition $\psi$.
The precondition and postconditions are formulae in the specification language.
We prove the soundness of the extended Hoare logic with regard to the denotational semantics.
We demonstrate the extended Hoare logic to prove the correctness of a root-finding program.

To summarize our main contributions:
\begin{itemize}
\item
A choice of primitive operations over real numbers and their computable (multi-valued) semantics in Section~\ref{s:preliminaries}.
\item
A small imperative programming language with three basic data types $\REAL$, $\INTEGER$, and $\LOGIC$, and the primitive 
operations for exact (multi-valued) computations over the algebraic structure of $\IR$, $\IZ$, and $\Kleene$ regarding the said semantics in Sections~\ref{s:ERC}~and~\ref{s:Semantic}.
\item
Examples demonstrating programming with these operations and semantics, such as:
multi-valued integer rounding, determinant via Gaussian elimination with full pivoting (subject to full-rank promise), and root-finding in Section~\ref{s:Programming}
\item
Proof that this programming language is Turing-complete:
those and only those real functions and integer multi-valued functions 
that can be realized in \ERC are computable in the sense of computable analysis
in Theorems~\ref{t:computable}~and~\ref{t:Complete}.
\item
A many-sorted structure for the formal specification of \ERC whose first-order theory is
decidable and model-complete in the case of $\ERC_0$ in Section~\ref{s:Logic} and Theorem~\ref{t:Eike}.
\item 
An extended Hoare logic to enable formal verification of $\ERC$ programs' behaviours including their terminations in 
Section~\ref{ss:Hoare} and Theorem~\ref{t:hoaresound}.
\item
An example correctness proof using the extended Hoare logic of the root-finding example program
in Section~\ref{ss:Trisection}.

\end{itemize}

This paper is structured as follows.
In Section~\ref{s:preliminaries}, we review the definitions of computable partial functions and computable partial multi-valued functions in computable analysis and the definition of Plotkin powerdomain over flat domains.
We specify several computable partial functions and multi-valued functions 
that constitute our primitive operators.
Section~\ref{s:ERC} defines the formal definition of our language syntax and
 typing rules.
In Section~\ref{s:Semantic}, we define domain-theoretic denotational semantics using the Plotkin powerdomain and 
prove that for any well-typed terms, well-formed commands, and well-typed programs, their denotations are computable partial multi-valued functions (or partial real functions in the case of real programs).
Section~\ref{s:Programming} defines some programming abbreviations and introduces various example programs which are: square root using Heron's method, exponential function via Taylor expansion and iterative, integer rounding, matrix determinants via Gaussian elimination, and root-finding.
Later in the section, based on the integer rounding example, we prove the Turing-completeness of $\ERC$ saying that any computable partial real functions and integer multi-valued functions are realizable in $\ERC$.

We introduce in Section~\ref{s:Logic} a three-sorted structure and its first-order language for specifying 
nondeterministic programs in $\ERC$. 
We prove that its first-order theory is decidable and model-complete in the case of $\ERC_0$.
Section~\ref{ss:Hoare} extends the classical Hoare logic
from $\INTEGER$ to ($\LOGIC$ and) $\REAL$ and prove its soundness.
Section~\ref{ss:Trisection} demonstrates our verification method by proving the correctness of the aforementioned \emph{tri}section program for root finding.
We conclude this work with Section~\ref{s:Conclusion} suggesting future research including extending \ERC from operating on real numbers to operating on functions to express operators and other higher-order data types.

\section{Preliminaries}
\label{s:preliminaries}

To devise sound mathematical semantics for unbounded loops, 
we take a traditional domain-theoretic approach \cite[Chapter~2.4]{reynolds2009theories} and define the denotations of loops as the least fixed points of some continuous operators.

To model nondeterministic computations, we use the Plotkin powerdomain (which we review later in this section) 
$\pdom{Y_\bot}$ on flat domains $Y_\bot$. 
An element of the domain $p \in \pdom{Y_\bot}$ is a subset of $Y\cup\{\bot\}$.
Hence, a mapping $f: X \to \pdom{Y_\bot}$, which is a partial multi-valued function, models a nondeterministic computation whose return values are in $Y$ where $f(x)\subseteq Y\cup\{\bot\}$ denotes the set of all possible return values of the computation at input $x$.
When $f(x)$ contains $\bot$, it denotes the case where there is a nondeterministic branch that fails.
In this sense, we say $f$ is partial and defined only at $x$ where $\bot \not \in f(x)$.

In this section, we review the definitions of computable partial functions $f: X \to Y_\bot$ and computable partial multi-valued functions $f: X \to \pdom{Y_\bot}$. 
We specify some computable partial functions and multi-valued functions that later constitute the set of primitive operators in 
$\ERC$. 

Instead of defining the computability of domain-theoretic functions by relating them with the ordinary computability definitions from computable analysis, we define the computability of functions to flat domains or to powerdomains 
directly such that we do not have to handle many different computability notions at the same time.
Remark~\ref{rem:domain computability} may help the readers who are already familiar with computable analysis to translate the settings back and forth.

\subsection{Representations and Computable Partial Functions}
For a set $X$, a representation $\delta_X$ of it is a partial surjective function from the set of infinite sequences of integers to $X$.
\footnote{The domains of representations can be any Cantor or Baire space where type-2 computability theory can be built on; cf. \cite[Exercise~3.2.17]{Wei00}.}
When $\delta_X(\varphi) = x \in X$, we say $\varphi \in \IZ^\IN$ is a $\delta_X$-\inlinedef{name} of $x$. 
We often omit the prefix $\delta_X$ if it is obvious from the context or 
if it is irrelevant.

We work with the following representations of $\IZ, \IR, \IK$:
\begin{align*}
\delta_\IZ(\varphi) &= k :\Leftrightarrow \varphi(n) = k \text{ for all } n\in\IN, \\
\delta_\IR(\varphi) &= x :\Leftrightarrow \myabs{x - \varphi(n)\cdot 2^{-n}} \leq 2^{-n}\text{ for all }n \in \IN,\\
\delta_\IK(0\myappend 0\cdots 0\myappend 1\myappend 1\cdots) &= \ctt, \\
\delta_\IK(0\myappend 0\cdots 0\myappend -\!1\myappend -\!1\cdots) &= \cff, \\
\delta_\IK(0\myappend 0\cdots 0\myappend 0 \cdots) &= \cuu. 
\end{align*}
Here, $a \myappend \phi$ when $a \in \IZ$ and $\phi \in \IZ^{\IN}$ denotes the infinite sequence where $a$ is prepended to $\phi$ as its new head. In a name of some $b \in \IK$, 
the digit $0$ acts as a token saying that we do not know yet if $b$ is $\ctt$ or $\cff$,
the digit $1$ acts as a token saying that we now know that $b$ is $\ctt$, and 
the digit $-1$ acts as a token saying that we now know that $b$ is $\cff$. 
Note that for the name $0\myappend0\myappend0\cdots$ of $b = \cuu$, 
it is not possible to decide in a finite time if $b$ is $\cuu$ or 
if it will happen to be $\cff$ or $\ctt$ in the future.
In other words, the procedure deciding the value of $b$, if $b = \ctt$ or $\cff$, diverges when $b = \cuu$.

Though we present the specific representations and work with them in this paper, 
the theory of admissibility \cite[Chapter~3.2]{Wei00} ensures their universality with regard to the topologies:
the discrete topology of $\IZ$, the standard topology of $\IR$, and the generalized Sierpinski topology of $\IK$:
\[ \big\{\emptyset,\{\ctt\},\{\cff\},\{\ctt,\cff\},\{\cuu,\ctt,\cff\}\big\}\;.\]
As we are only interested in computability in this paper, the underlying representations can be replaced with any other computably equivalent representations.

When we have a function to a flat domain $f: X \to Y_\bot$, we consider it as a partial function which is properly defined at $\mydom{f} \dfeq  \{x \in X\mid f(x) \neq \bot\}$.
When $X$ and $Y$ are represented by $\delta_X$ and  $\delta_Y$, respectively, $f$ is said $(\delta_X, \delta_Y)$-\inlinedef{computable} if there is a type-2 machine \cite[Chapter~2.1]{Wei00} that
for all $x \in \mydom{f}$ transforms each $\delta_{X}$-name of $x$ to a $\delta_{Y}$-name of $f(x)$.\footnote{It also can be defined equivalently using oracle machines \cite{kapron1996new}.}
In this case,
we say the type-2 machine 
\inlinedef{realizes} $f$.
We drop the prefix $(\delta_X, \delta_Y)$ if it is obvious or irrelevant which specific representations underlie. 
An important fact to note here is that the realizer of $f: X \to Y_{\bot}$
does not need to diverge in the case $f(x) = \bot$. 
The realizer on any name of such $x$ not in the domain can either diverge or even compute a meaningless result.

We often call a function $f: X \to Y_\bot$ \emph{partial} to emphasize that $\bot$ may be a function value.
Often we call $f$ which does not admit $\bot$ as its function value \emph{total} when we want to emphasize it.
To make our notation more informative, let us write $f: X \to Y$ in the case $f: X \to Y_\bot$ is total.

For a multivariate $f: X_1 \times \cdots \times X_d \to Y_\bot$, its computability is defined by the existence of a type-2 machine with multiple inputs realizing it in a similar manner.
Equivalently, it can be defined with regard to the product representation for its domain: 
\[\delta_{X_1\times\cdots\times X_d}(\varphi) = (x_1, \cdots, x_d)
\text{ if and only if } \delta_{X_j}(n \mapsto \varphi(d \times n + j-1))  = x_j \text{ for each } j.\]
Again, other computably equivalent representations of products can replace this specific product representation.

The following are examples of computable partial and total functions:
\begin{exas}[Primitive operations]
\hfill
\label{exa:primitive1}

\begin{enumerate}
    \item For any represented set $X$, its identity function $\id : X \to X$ is computable.

    \item The constants $\ctt, \cff, \cuu \in \IK$ and $k \in \IZ$ are computable as constant functions from any represented set.

    \item The logic operations for Kleene three-valued logic $\kneg : \IK \to \IK$, $\klor, \kland : \IK \times \IK \to \IK$ defined as the following tables are computable. Note that $\cuu$ stands for indeterminacy.

\vspace{1em}
{
\centering
    \begin{tabular}{c||c|c|c}
$\kland$ & $\ctt$ & $\cff$ & $\cuu$ \\\hline\hline
$\ctt$   & $\ctt$ & $\cff$ & $\cuu$ \\\hline
$\cff$   & $\cff$ & $\cff$ & $\cff$ \\\hline
$\cuu$   & $\cuu$ & $\cff$ & $\cuu$ 
    \end{tabular}
}

    \vspace{1em}

{
\centering
    \begin{tabular}{c||c|c|c}
$\klor$ & $\ctt$ & $\cff$ & $\cuu$ \\\hline\hline
$\ctt$   & $\ctt$ & $\ctt$ & $\ctt$ \\\hline
$\cff$   & $\ctt$ & $\cff$ & $\cuu$ \\\hline
$\cuu$   & $\ctt$ & $\cuu$ & $\cuu$ 
    \end{tabular}
}

    \vspace{1em}

{
\centering
    \begin{tabular}{c||c|c|c}
$\kneg$ & $\ctt$ & $\cff$ & $\cuu$ \\\hline\hline
  & $\cff$ & $\ctt$ & $\cuu$ 
    \end{tabular}
}
    \vspace{1em}
    
\noindent This justifies the name \emph{Kleenean} for $\IK$.

    \item The operators of the Presburger arithmetic $+ : \IZ \times \IZ \to \IZ$, $- : \IZ \to \IZ$, and the $\IK$-valued comparison tests $\ieq, \ileq : \IZ \times \IZ \to \IK$ (which does not admit $\cuu$) 
are computable.
    
    \item The field operators $+, \times : \IR \times \IR \to \IR$, $ - : \IR  \to \IR$, and ${}^{-1} : \IR  \to \IR_\bot$ where the multiplicative inversion is not defined at $0$, i.e., $0^{-1} = \bot$, are computable. (Any realizer for the multiplicative inversion function  will diverge on $0^{-1} = \bot$.)
    
    \item The exponentiation $\imath : \IZ \ni p \mapsto 2^p \in \IR$ called 
    precision embedding is computable.
    
    \item The $\IK$-valued real number comparison $\elt : \IR \times \IR \to \IK $ from Equation~(\ref{eqn:elt}) is computable. 
    It is realized by a procedure that given some names $\phi_{x}, \phi_{y}$ of $x, y \in \IR$
    computes the sequence of integers:
    \[n \mapsto 
    \begin{cases} 
    	1 &\text{if } \phi_{x}(n) + 2 < \phi_{y}(n),\\
		-1 &\text{if } \phi_{y}(n) + 2< \phi_{x}(n),\\
		0 &\text{otherwise}.
	\end{cases}    
    \]
    We can easily verify that the resulting sequence is a name of $x \elt y$.

    \item A partial projection map $\semproj : X^d \times \IZ \to X_\bot$ defined by
    \[
        \semproj((x_0, \cdots, x_{d-1}), n) \dfeq \begin{cases}
        x_n &\text{if } 0\leq n <d, \\
        \bot &\text{otherwise},
        \end{cases}
    \]
    and a partial updating map $\semupdate : X^d \times \IZ \times X \to (X^d)_\bot$ defined by 
    \[
        \semupdate((x_0, \cdots, x_{d-1}), n, x) \dfeq \begin{cases}
        (x_0, \cdots, x_{n-1}, x, x_{n+1}, \cdots, x_{d-1})  &\text{if } 0\leq n <d, \\
        \bot &\text{otherwise},
        \end{cases}
    \]
    are computable. 
    \item
    A total projection map $\semproj_i : X_1 \times \cdots \times X_d \to X_i$ defined by $\semproj_i(x_1, \cdots, x_d) \dfeq x_i$ and a total updating map $\semupdate_i : X_1 \times \cdots \times X_d \times X_i \to X_1 \times \cdots \times X_{d}$ defined by
    $\semupdate_i((x_1, \cdots, x_d), y) \dfeq (x_1, \cdots, x_{i-1}, y, x_{i+1}, \cdots, x_d)$ updating the $i$'th entry $x_i$ with $y$, 
    are computable.
\end{enumerate}
\end{exas}

\subsection{Computable Infinite Multi-valued functions}
Partial multi-valued functions are 
non-empty set-valued functions 
$f: X \to \pset{Y_\bot}$ that model nondeterministic computations in the following sense.
For an input $x$, 
$f(x)$ is the set of all possible nondeterministic values that the computation modelled by $f$ can yield, 
including $\bot$ if there is a nondeterministic branch for which the computation fails.
We say a function of type $f : X \to \pset{Y_\bot}$ a (\emph{partial}) \emph{multi-valued function} which is properly defined at $\mydom{f}\dfeq\{x \mid \bot \not \in f(x)\}$.
We also often call a partial multi-valued function \emph{total} when its function values never contain $\bot$; i.e., $\mydom{f} = X$.

For a partial multi-valued function $f : X \to \pset{Y_\bot}$, when the underlying sets $X, Y$ are represented by $\delta_X$ and $\delta_Y$, respectively, we say that $f$ is $(\delta_X, \delta_Y)$-\inlinedef{computable} if there is a
type-2 machine that for all $x \in \mydom{f}$
transforms each $\delta_X$-name of $x$ to a $\delta_Y$-name of some $y \in f(x)$.
Note that the \inlinedef{realizer}'s behaviour on non-valid inputs is not specified.
It can either diverge or compute a meaningless value. 
We omit the prefix $(\delta_X, \delta_Y)$ if it is not necessary to explicitly specify the underlying representations.
Note that nondeterminism happens due to the fact that even for a fixed $x$, for its different names, the realizer can compute different $y, y' \in f(x)$ though the realizer itself is deterministic.

\subsection{Computable (Finite) Multi-valued functions and Plotkin Powerdomain}
Using the entire non-empty power-set $\pset{Y_\bot}$ 
causes the same problem that motivated \cite[Section~2]{DBLP:journals/siamcomp/Plotkin76}
that it does not give us a nice order theoretic property in defining denotational semantics.
Reminding that $\ERC$ only provides finite nondeterminism, and that as long as only real numbers and integers are allowed as inputs there is no necessarily infinite integer multi-valued function (see Fact~\ref{fact:prelim 2}), we safely restrict to Plotkin powerdomain \cite{DBLP:journals/siamcomp/Plotkin76}.
For a set $X$, the powerdomain on its flat domain $X_\bot$ is defined by
\[
\pdom{X_\bot} \dfeq \{
p \subseteq X\cup\{\bot\} \mid p \neq\emptyset, \text{ and } p \text{ is finite or contains } \bot
\} \subseteq \pset{X_\bot}
\]
endowed with Egli-Milner ordering:
\[
p \sqsubseteq q \Leftrightarrow
(\forall x\in p\,\exists y\in q.\,x \leq y) \land
     (\forall y\in q\,\exists x\in p.\,x \leq y)\;.
\]
Here, $x \leq y$ is the order on the flat domain that holds if and only if $x = \bot$ or $x = y$.
The constructed powerdomain for any flat domain $X_\bot$ is known to be a $\omega$-cpo with the least element $\{\bot\}$ \cite[Chapter~7.2]{reynolds2009theories}.
We say a partial multi-valued function of type $f : X \to \pdom{Y_\bot} \subseteq \pset{Y_\bot}$ \emph{finite}.

For a function to a flat-domain $f : X \to Y_\bot$, let us write $f^\ddagger : X \to \pdom{Y_\bot}$ 
 for the embedding $f^\ddagger(x) \dfeq \{f(x)\}$.
For a multivariate partial multi-valued function $f : X_1 \times \cdots \times X_d \to \pdom{Y_\bot}$,
its extension $f^\dagger :  \pdom{(X_1)_\bot} \times \cdots \times \pdom{(X_d)_\bot} \to \pdom{Y_\bot}$ is defined by
    \[
        f^\dagger (p_1, \cdots, p_d) \dfeq \bigcup_{x_1 \in p_1, \cdots, x_d \in p_d}\begin{cases}
        \{f(x_1, \cdots, x_d)\} &\text{if } \forall i.\; x_i \neq \bot, \\
        \{\bot\} &\text{otherwise.}
        \end{cases}
    \]
Furthermore, 
the extension by one argument $f^{\dagger_i} : X_1 \times \cdots \times  \pdom{(X_i)_\bot} \times \cdots \times X_d \to \pdom{Y_\bot}$ is defined by
    \[
        f^{\dagger_i} (x_1, \cdots, p, \cdots, x_d) \dfeq \bigcup_{x_i \in p}\begin{cases}
        \{f(x_1, \cdots,x_i, \cdots x_d)\} &\text{if } x_i \neq \bot, \\
        \{\bot\} &\text{otherwise.}
        \end{cases}
    \]
The embedding and extensions can be more naturally understood in the setting of category theory \cite{brookes1993monads,Jung}.

The composition of two partial multi-valued functions $f: X \to \pdom{Y_\bot}$ and $g: Y \to \pdom{Z_\bot}$ is expressed by $g^\dagger \circ f: X \to \pdom{Y}$.
This composition is continuous in both arguments \cite[Section~6]{DBLP:journals/siamcomp/Plotkin76}.
In other words, given $f, f_i : X \to \pdom{Y_\bot}$ and $g, g_i : Y \to \pdom{Z_\bot}$ 
where $f_i$ and $g_i$ are increasing chains with regard to the point-wise orderings, 
it holds that
\[
\bigsqcup_{i\in \IN} (g^\dagger \circ f_i) =  g^\dagger \circ \bigsqcup_{i\in \IN} f_i \quad \text{and} \quad
\bigsqcup_{i\in \IN} (g_i^\dagger\circ f) =  \big(\bigsqcup_{i\in \IN} g_i\big)^\dagger \circ f .
\]
Partial composition is also continuous:
Consider a multivariate multi-valued function $f : X_1 \times \cdots \times X_d \to \pdom{Y_\bot}$, a list of functions $g_j : X \to X_j$ for $j \neq i$, and a chain of multi-valued functions $g_{i, k} : X \to \pdom{(X_i)_\bot}$ in $k$. Then,
\[
\bigsqcup_{k\in \IN}f^{\dagger_i}\circ \langle g_1, \cdots, g_{i, k}, \cdots g_d \rangle 
=
f^{\dagger_i}\circ \langle g_1, \cdots, \bigsqcup_{k\in \IN}g_{i, k}, \cdots g_d \rangle \;.
\]
Here, for a finite collection of functions $f_1, \cdots, f_d$ sharing their domains, $\langle f_1, \cdots, f_d\rangle$ denotes $x \mapsto (f_1(x), \cdots, f_d(x))$.

For each set $X$, define $\semcond : \IK \times \pdom{X_\bot}\times \pdom{X_\bot} \to \pdom{X_\bot}$ by 
    \[\semcond(b, p, q) \dfeq \begin{cases}
            p & \text{if} \; b = \ctt, \\
            q & \text{if} \; b = \cff, \\
            \{\bot\} & \text{otherwise},
        \end{cases}\]
        a conditional combinator of $\IK$.
For any $b : X \to \pdom{\IK_\bot}$ and $S : X \to \pdom{X_\bot}$, 
due to the continuity of compositions, the following operator can be seen as continuous:
\[
\wop(b, S)(f : X \to \pdom{X_\bot}) \dfeq \semcond^{\dagger_1}\circ \langle b,f^\dagger \circ S, \id^\ddagger\rangle\ .
\]
The continuity ensures the existence of the least fixed point $\lfp\big(\wop(b, S)\big): X \to \pdom{X_\bot}$ which is used to define the denotations of loops later in Section~\ref{s:Semantic}.

The following are examples of computable (finite) multi-valued functions:
\begin{exas}[Primitive multi-valued operations]
\label{exa:primitive2}
\hfill
\begin{enumerate}
    \item When $f : X \to Y_\bot$ is a computable partial/total function, its embedding $f^\ddagger : X \to \pdom{Y_\bot}$ is a computable partial/total multi-valued function.

    \item When $f : X \to \pdom{Y_\bot}$ and $g : Y \to \pdom{Z_\bot}$ are computable partial multi-valued functions, their composition $g^\dagger \circ f : X \to \pdom{Z_\bot}$ is a computable partial multi-valued function.

    \item When $f_i : X \to \pdom{(Y_i)_\bot}$ and $g : Y_1 \times \cdots \times Y_d \to Z$ are computable multi-valued functions, their composition $g^\dagger \circ \langle f_1, \cdots, f_d\rangle : X \to \pdom{Z_\bot}$ is a computable multi-valued function.

    \item For each $n \in \IN$, the nondeterministic selection function $\semchoose_n : \IK^n \to \pdom{\IZ_\bot}$ defined by
    \[
    \semchoose_n(b_0, \cdots, b_{n-1}) \dfeq 
    \begin{cases}
    \{i \mid b_i = \ctt\} &\text{if there is }i\text{ s.t. } b_i =\ctt,\\
    \{\bot\} &\text{otherwise,}
    \end{cases}
    \]
    is a computable partial multi-valued function.
    
	It is realized by the procedure which when it receives $\phi_{0}, \cdots, \phi_{n-1}$
	as some names of $b_{0}, \cdots, b_{n-1}$,
	for each $p = 0, 1,\cdots,$ checks $\phi_{i}(p)$ until it finds $\phi_{i}(p) = 1$ implying $b_{i} = \ctt$.
	If the realizer finds such $i$, it returns $i\myappend i \myappend i  \myappend \cdots$ as the name 
	of the index $i \in \semchoose_n(b_0, \cdots, b_{n-1})$.
	Note that depending on the given names, the realizer can pick a different index.
	If there is no choosable index, the realizer will diverge.

    \item 
    For any computable $f : X \to \pdom{\IK_\bot}$ and $g, h : X \to \pdom{Y_\bot}$, the composite multi-valued function
    \[\semcond^{\dagger_1}\fcomp \langle f, g, h\rangle : X \to \pdom{Y_\bot}\]
    is computable.     
    \item
    Define $\semkond : \IK \times \pdom{\IR_\bot} \times \pdom{\IR_\bot} \to \pdom{\IR_\bot}$, a continuous extension of $\semcond$ by 
\[ 
\semkond(b, p, q) \dfeq \begin{cases}
                p & \text{if} \; b = \ctt, \\
                q & \text{if} \; b = \cff, \\
                \{x\} & \text{if}\; b = \cuu \;\wedge\; p = q = \{x\} ,\\
                \{\bot\}&\text{otherwise},
            \end{cases}
    \]
    with regard to the topology of $\IK$; cmp. \cite[Theorem~2.3.8]{Wei00} and ``parallel if'' from PCF \cite{plotkin1977lcf}.
    For any computable $f : X \to \pdom{\IK_\bot}$ and $g, h : X \to \pdom{\IR_\bot}$, the composition
    $\semkond^{\dagger_1}\fcomp \langle f, g, h\rangle : X \to \pdom{\IR_\bot}$
    is a computable partial multi-valued function.
    
    Here, we make a further explanation of this operation. 
    Though this operation at first glance seems to violate 
    the principle that we cannot make a case distinction on $b \in \IK$,
    since this operation is continuous at $b = \cuu$, it is not the case. 
    The operation can be realized without testing $b = \cuu$ or $p = q$. 
 	We describe one possible realizer of $\semkond^{\dagger_1}\fcomp \langle f, g, h\rangle$ to make it clearer.
	
	For any $u \in X$, write $r \dfeq f(u) \in \pdom{\IK_\bot}$, $p \dfeq g(u), q\dfeq h(u) \in \pdom{\IR_{\bot}}$ and assume $\bot \not\in \semkond^{\dagger_1}(r, p, q)$. 
	(Recall that we do not need to analyze 
	how our realizer behaves when $\bot \in \semkond^{\dagger_1}(r, p, q)$.)
    Since $\bot$ is not in any of $p, q, r$, 
    we are given names $\phi_{b},\phi_{x}, \phi_{y}$ of $b, x, y$, respectively, where $b \in r$, $x \in p$, and $y \in q$.
	For some $z \in \semkond^{\dagger_1}(r, p, q)$ and its name $\phi$, the realizer needs to compute
	$n \mapsto \phi(n)$.

	The realizer first tests if $\phi_{x}(n+2)$ and $\phi_{y}(n+2)$ differ by more than $2$. 
	In this case (i), since we can conclude that $x \neq y$, we test if $b = \ctt$ or $b  = \cff$ by 
	inspecting $\phi_{b}(0), \phi_{b}(1), \cdots$. 
	If $b = \ctt$, since $x \in \semkond^{\dagger_1}(r, p, q)$, let $x$ be the $z$ and $\phi_{x}$ be the $\phi_{z}$
	and let the realizer return $\phi_{x}(n)$ for the $\phi_{z}(n)$.
	We can let the realizer work similarly in the other case $b = \cff$.
	If $b$ was $\cuu$, the realizer's inspecting procedure will diverge. 
	However, this is against our assumption that $b = \cuu$ and $p \neq q$ makes $\bot  \in \semkond^{\dagger_1}(r, p, q)$.
		
	For the other case (ii) where 
	$\phi_{x}(n+2)$ and $\phi_{y}(n+2)$ are different by at most $2$,
	 from the two integers, the realizer can compute $k \in \IZ$ 
	where $\myabs{x - k \cdot 2^{-n}} \leq 2^{-n}$ and $\myabs{y - k \cdot 2^{-n}} \leq 2^{-n}$ hold.
	Without being concerned with $b$, the realizer returns the computed $k$ for $\phi_{z}(n)$ as a valid 
	$n$'th entry of some names of both $x$ and $y$.
	There are some hypothetical cases to analyze.
	If $x = y$, then the realizer will execute this second case (ii) indefinitely for all $n$
	and yield a name of $x$ (which can be different from $\phi_{x}$ or $\phi_{y}$).
	Regardless of $b$, even when $b = \cuu$, since $x \in \semkond^{\dagger_1}(r, p, q)$ 
	as long as $\bot \not\in\semkond^{\dagger_1}(r, p, q)$, we can conclude that the realizer is computing a correct name. 
	Otherwise, if $x \neq y$, after returning some finite entries of a name 
	that works for both $x$ and $y$, the realizer will execute the first case (i) which is already analyzed to be correct.
	
	Note that in the case $b = \cuu$, $p, q$ being singleton is crucial. 
	When we ease the condition to $p = q$ without forcing them to be a singleton, 
	there is a case where the nondeterministically picked inputs $x \in p$ and $y \in q$ are different.
	In this case, at some point, the realizer executes the first case (i) and diverges by inspecting $b$.

\end{enumerate}
\end{exas}

Observe from the definitions that a realizer of a partial (multi-valued) function can be seen also as a realizer of some other partial (multi-valued) functions  because 
(1) the out-of-domain behaviour of realizers is not specified 
and (2) for a partial multi-valued function, its realizer only has to compute some of the possible function values.
Specifically:
\begin{fact}
\label{fact:prelim}
\label{fact:prelim 1} 
    If a partial function $f: X \to Y_\bot$ is computable, any domain restriction
    $g : X \to Y_\bot$ such that $\forall x.\; g(x) \leq f(x)$ is computable by the same realizers. 
    Similarly, if a partial multi-valued function $f : X \to \pset{Y_\bot}$ is computable, 
    any domain restriction and function-value enlargement 
	$g : X \to \pset{Y_\bot}$ such that $\forall x.\; g(x) \sqsubseteq f(x)$ or $f(x) \subseteq g(x)$ is computable by the same realizers.
\end{fact}

For any real number or integer, the set of its names is compact with regard to the standard Baire space topology on $\IZ^\IN$. 
Hence, for any type-2 computable realizer of an integer partial multi-valued function $f : X_1 \times \cdots \times X_d \to \pdom{\IZ_\bot}$ where $X_i = \IZ$ or $\IR$, for fixed inputs $(x_1, \cdots, x_d)$, the image of the 
realizer of the compact set of names is also compact as any type-2 computable realizer is continuous.
That means even when $f(x_1, \cdots, x_d)$ is infinite, there are only finitely many different integers that its realizer actually computes; see \cite[Theorem~4.2]{brattka1995computable}. 
\begin{fact}
\label{fact:prelim 2} 
    Any computable multi-valued function $f : X_1\times \cdots \times X_d \to \pset{\IZ_\bot}$ 
    where $X_i = \IZ $ or $\IR$, admits finite refinement in the sense that 
    there is a computable finite multi-valued function $g : X_1 \times \cdots \times X_d \to \pdom{\IZ_\bot}$ 
    such that $\forall (x_1, \cdots, x_d).\; g(x_1, \cdots, x_d) \subseteq f(x_1, \cdots, x_d)$ holds.
\end{fact}

In addition to the above justification of our use of Plotkin powerdomain, the following is a note on the continuity issue regarding nondeterminism:
\begin{rem}
\label{r:hoare}
Note that $\semchoose_n : \IK^n \to \pdom{\IZ_\bot}$ is not (domain-theoretic) continuous when $\IK$ is ordered by $b_1 \leq b_2$ iff $b_1 = b_2$ or $b_1 = \cuu$.
This is observed in \cite{marcial2007semantics} and there it is suggested to use Hoare powerdomain instead.
However, our approach distinguishes (domain-theoretic) continuity and computability of primitives:
we do not use (domain-theoretic) continuity as a criterion for computability.
We use continuity to define the semantics of loops whose computability is studied 
separately using the notion of type-2 computability of realizers. 
In other words, as long as we have a fact that $\semchoose_n: \IK^n \to \pdom{\IZ_\bot}$ is computable, there is no need to make it (domain-theoretic) continuous by giving an order on $\IK$.
This way we still can use Plotkin powerdomain in defining our semantics such that we do not need to have separate operational semantics to express termination as in \cite{marcial2007semantics}.
\end{rem}

We conclude this section with a remark on how the formulations of computable partial functions and multi-valued functions 
in this paper are related to the traditional formulations in computable analysis.
\begin{rem}
\label{rem:domain computability}
In the traditional setting of computable analysis:
\hfill
\begin{enumerate}
    \item 
         The computability of a partial function presented by $f \pcolon X \to Y$ is equivalent to the computability of $f_\bot: X \to Y_\bot$ defined by $x \mapsto f(x)$ if $x \in \dom(f)$ and $\bot$ if $x \not\in \dom(f)$. Similarly, the computability of $f: X \to Y_\bot$ is equivalent to the computability of $f\restriction_{\mydom{f}} \pcolon X \to Y$.
         Here, $\restriction$ is a notation for the domain restriction and
         recall that $\mydom{f} = \{x \mid f(x)\neq \bot\}$.

    \item
        The computability of a partial multi-valued function presented by $f \pcolon X \toto Y$, which is a partial non-empty set-valued function, is equivalent to the computability of $f_{\{\bot\}}: X \to \pset{Y_\bot}$.
        Similarly, the computability of $f : X \to \pset{Y_\bot}$ is equivalent to the computability of $f\restriction_{\mydom{f}} \pcolon X \toto Y$.
        Recall that $\mydom{f} = \{x \mid \bot \not \in f(x)\}$.

\end{enumerate}
\end{rem}

\section{Formal Syntax and Typing Rules of \emph{Exact Real Computation}} \label{s:ERC}
In this section, we provide the formal syntax and typing rules of $\ERC(\calF, \calG)$ for any $\calF$ and $\calG$ according to the following convention: 
\begin{conv}
\label{c:Expansion}
Fix a (possibly empty) finite set $\calF$ of computable partial functions to $\IR_\bot$ from products of $\IZ$, $\IR$, and $\IR^n$ for any $n \geq 1$, as well as a (possibly empty) finite set $\calG$ of computable finite partial multi-valued functions to $\pdom{\IZ}$ from products of $\IZ$, $\IR$, and $\IR^n$ for any $n \geq 1$.
\end{conv}
\noindent

\subsection{Syntax} \label{ss:Syntax}
$\ERC$ is an imperative programming language
comprising of the following axiomatized constituents:
data types (Section~\ref{sss:Syntax_Type}),
terms (Section~\ref{sss:Syntax_Term}),
commands (Section~\ref{sss:Syntax_Command}),
and programs (Section~\ref{sss:Syntax_Program}).

\subsubsection{Data Types}
\label{sss:Syntax_Type}

The data types that $\ERC$ provides are as follow:
\[
    \tau \Coloneqq \dbool \mid \dint \mid \dreal \mid \dreal[n]\text{ for each }n= 1, 2, \cdots
\]
See that $\ERC$ provides countably many data types: for each natural number $n \geq 1$, there is a data type $\dreal[n]$, which represents the set of arrays of real numbers of fixed length $n$.
Note $n$ here is a meta-level expression which is not a $\ERC$ ``term'' from Section~\ref{sss:Syntax_Term}.

\subsubsection{Terms}
\label{sss:Syntax_Term}
Although the type of a term is judged later by typing rules in Section~\ref{ss:Type},
here we follow the following conventions:
Write $m, \ell$ to denote terms which should be typed $\dint$,
$u, v$ to denote terms which should be typed $\dreal$, and
$b$ to denote terms which should be typed $\dbool$.
Arbitrary terms are denoted by $t$.
For some fixed countable set of variables $\calV$, 
the term language of $\ERC$ is defined as follows.
\begin{align*}
& t, m, \ell, u, v, b \Coloneqq  & & &\\
& \;\;\; 
k & & k\in \IZ; \;\dint\text{ literal}& \\
& \mid \ttrue \mid \tfalse \mid \tundef & & \dbool \text{ literal}& \\
& \mid [t_1, \cdots, t_n] & & \text{array literal}& \\
& \mid x                  & & x \in \calV;\; \text{variable}&\\
& \mid t[m] & & \text{array access} & \\
& \mid t_1 + t_2 \mid - t  & & \text{addition and additive inversion} &\\
& \mid m \leq \ell  \mid  m  \ieq  \ell &  & \text{integer comparison and equality test} &\\ 
& \mid u \mult  v \mid u^{-1}  & & \text{real multiplication and multiplicative inversion} &\\
& \mid u \elt v & & \text{real comparison} &\\
&  \mid \myprec{m} & & \text{precision embedding} & \\
& \mid \kneg b  \mid b_1 \kland b_2 \mid b_1 \klor b_2 & & \text{Kleene logic} &\\
& \mid f(t_1,\cdots, t_d) & & (f : X_1\times \cdots \times X_d \to \IR_\bot)\in\calF; \;\text{primitive function call} &\\
& \mid g(t_1,\cdots, t_d) & & (g : X_1\times \cdots \times X_d \to \pdom{\IZ_\bot})\in\calG; &\\
&  & & \quad\text{primitive multi-valued function call} &\\
& \mid \mychoose_n(b_0, b_1, \cdots, b_{n-1}) & & n\in\IN;\;\text{multi-valued choice} & \\
& \mid \mycond{b}{u}{v} & & \text{continuous conditional} &
\end{align*}

Note that we do not provide integer multiplication and any other coercion than $\myprec{m}$ unless they are provided through $\calF, \calG$.
For each natural number $n$, $\mychoose_n$ is a distinct term construct. However, since the $n$ is determined syntactically by the number of the arguments, we may omit it and instead, simply refer to $\mychoose()$.
We introduce $m  \ieq  \ell$ as a term construct instead of
as an abbreviation for $m\leq \ell \kland m\leq \ell$ because their
denotations which are defined later in Section~\ref{sss:Semantic_Term} do not agree.

\subsubsection{Commands}
\label{sss:Syntax_Command}
$\ERC$ provides a strict distinction between its term language and command language where loops belong only to the command language.
Commands in \ERC are inductively constructed as follows.
\begin{align*}
& S\;\;\;\;\Coloneqq \,\;\; \eskip & & \text{skip} &\\
& \qquad\quad\;\; \mid x \dfeq   t & & \text{variable assignment} &\\
& \qquad\quad\;\; \mid x[m] \dfeq   t & & \text{array assignment} &\\
& \qquad\quad\;\; \mid \elet\; x : \tau =  t & & \text{variable declaration} &\\
& \qquad\quad\;\; \mid S_1;\;S_2 & & \text{sequential composition} &\\
& \qquad\quad\;\; \mid \eif\; b \;\ethen\; S_1 \;\eelse\; S_2 \;\eendif
    & & \text{branching} &\\
& \qquad\quad\;\; \mid \ewhile\;b\;\edo\;S \;\eendwhile & & \text{loop} &\\
\end{align*}

The variable declaration $\elet$ statement introduces a variable with limited scope in a loop or a conditional branch; see the related typing rules in \autoref{fig:command_type}.
\subsubsection{Programs}
\label{sss:Syntax_Program}
Having data types, terms, and commands defined,
we can finally define what a program in $\ERC$ is.
A program in $\ERC$ is defined to be in either of the two forms:
\[
\begin{array}{l}
\textbf{input}\;x_1 : \tau_1, x_2 : \tau_2, \cdots, x_n : \tau_n \\
\quad S \\
\textbf{return}\;t
\end{array}
\quad\quad\quad
\begin{array}{l}
\blim{p} x_1 : \tau_1, x_2 :\tau_2, \cdots, x_n : \tau_n \\
\quad S \\
\blimreturn{t}{p}
\end{array}
\]
Here, $x_i$ and $p$ are variables varying over $\calV$ and $\tau_i$ varies over $\{\dint, \dreal, \dreal[1], \dreal[2], \cdots\}$. 
The variable $p$ that is separated from the other input variables is the precision parameter.
The intended behaviour of a program in the right kind is, taking $x_1, \cdots, x_n$ as inputs, 
to compute the limit of the command $S$ and the return term $t$ as $p$ goes to $-\infty$.
This is formally defined later in Section~\ref{s:Semantic}.

\subsection{Typing Rules}
\label{ss:Type}
Here we define well-typedness of terms (Section~\ref{sss:Type_Term}), well-formedness of commands (Section~\ref{sss:Type_Command}), and well-typedness of programs (Section~\ref{sss:Type_Program}) under a \inlinedef{context}
which is a mapping from a finite set of variables to their corresponding types. 
We use a list of assignments $\Gamma \dfeq  x_1 : \tau_1, x_2 : \tau_2, \cdots, x_n : \tau_n$ to represent the context mapping a variable $x_i$ to its data type $\tau_i$ for $i = 1, \ldots, n$.

\subsubsection{Well-Typed Terms}
\label{sss:Type_Term}

Well-typedness of a term $t$ in $\ERC$ to $\tau$ under a context $\Gamma$ is
written as $\Gamma \vdash t :\tau$.
\autoref{fig:term_type} shows \ERC's type inference rules.
Note that $+$ and $-$ are used both in integer and real number operations. 

\begin{figure*}[!ht]
\caption{Typing rules for terms 
\label{fig:term_type}}
\fbox{
\noindent\begin{minipage}{0.97\textwidth}
\begin{mathpar}
\infer{\Gamma \vdash k : \dint}{%
}
\and
\infer{\Gamma \vdash \ttrue : \dbool}{%
}
\and
\infer{\Gamma \vdash \tfalse : \dbool}{%
}
\and
\infer{\Gamma \vdash \tundef : \dbool}{%
}
\and
\infer{\Gamma \vdash [t_1, \cdots, t_n] : \dreal[n]}{%
\Gamma\vdash t_i : \dreal
\quad(\text{for } i =1, \cdots, n)
}
\and
\infer{\Gamma \vdash x : \tau}{%
\Gamma(x) = \tau
}
\and
\infer{\Gamma \vdash t[m] : \dreal}{%
\Gamma \vdash m : \dint &
\Gamma \vdash t : \dreal[n]
}
\and
\infer[(\bop, \tau, \tau') \in \bopty]
{\Gamma \vdash t_1 \bop t_2 : \tau'}{%
\Gamma \vdash t_1 : \tau
&
\Gamma \vdash t_2 : \tau 
}
\and
\infer[(\uop, \tau, \tau') \in \uopty ]
{\Gamma \vdash \uop t : \tau'}{%
\Gamma \vdash t : \tau 
}
\and
\infer{\Gamma \vdash f(t_1, \cdots, t_{d}) : \dreal}{%
(f : X_1 \times \cdots \times X_d \to \IR_\bot) \in \mathcal{F}&
\Gamma \vdash t_i : \langle X_i \rangle
\quad(\text{for } i =1, \cdots, d)
}
\and
\infer{\Gamma \vdash g(t_1, \cdots, t_{d}) : \dint}{%
(g : X_1 \times \cdots \times X_d \to \pdom{\IZ_\bot}) \in \mathcal{G}&
\Gamma \vdash t_i : \langle X_i \rangle
\quad(\text{for } i =1, \cdots, d)
}
\and
\infer{\Gamma \vdash \mychoose_n(b_0, b_1, \cdots, b_{n-1}) : \dint}{%
\Gamma \vdash b_i: \dbool\quad(\text{for } i =0, 1, \cdots, n-1)
}
\and
\infer{\Gamma \vdash (\mycond{b}{u}{v}) : \dreal}{%
\Gamma \vdash b : \dbool &
\Gamma \vdash u : \dreal &
\Gamma \vdash v : \dreal
}
\medskip
\end{mathpar}
\begin{align*}
\bopty &\dfeq \{ (+, \REAL, \REAL), (+, \dint, \dint),
(\leq, \dint, \dbool),
(\ieq, \dint, \dbool),
(\elt, \dreal, \dbool),(\times, \dreal, \dreal), \\
&\qquad
(\kland, \dbool, \dbool),
(\klor, \dbool, \dbool)
\}\\
\uopty &\dfeq \{
(-, \dreal, \dreal), 
(-, \dint, \dint), 
(2^{\square}, \dint, \dreal), 
(\kneg, \dbool, \dbool),
({\square}^{-1}, \dreal, \dreal)
\}\\
\langle \IZ \rangle &\dfeq \dint\\
\langle \IR \rangle &\dfeq \dreal\\
\langle \IR^n \rangle &\dfeq \dreal[n]\quad (\text{for }n = 1, 2, \cdots)
\end{align*}
\end{minipage}}
\end{figure*}

\subsubsection{Well-formed Commands}
\label{sss:Type_Command} 
Unlike terms, a command in \ERC may modify contexts. Let us denote
a command $S$ under a context $\Gamma$ being well-formed
and yielding a new context $\Gamma'$ as
$
\Gamma \vdash S \triangleright \Gamma'\, .
$
The Well-formedness of commands is defined by the inference rules in \autoref{fig:command_type}.

\begin{figure*}[!ht]
\caption{Typing rules for commands \label{fig:command_type}}
\fbox{\noindent\begin{minipage}{0.97\textwidth}
\begin{mathpar}
\infer{\Gamma \vdash \eskip  \triangleright \Gamma}{%
}
\and
\infer{\Gamma \vdash x\dfeq  t  \triangleright \Gamma}{%
\Gamma \vdash t : \tau &
\Gamma(x) = \tau &
}
\and
\infer{\Gamma \vdash x[m]\dfeq  t  \triangleright \Gamma}{%
\Gamma \vdash m : \dint &
\Gamma \vdash t : \dreal &
\Gamma(x) = \dreal[n] &
}
\and
\infer{\Gamma \vdash \elet\;x : \tau = t\; \triangleright \Gamma, x:\tau}{%
x\not\in \dom(\Gamma)&
\Gamma \vdash t : \tau
}
\and
\infer{\Gamma \vdash S_1;S_2 \triangleright \Gamma_2}{%
\Gamma \vdash S_1 \triangleright \Gamma_1 &
\Gamma_1 \vdash S_2 \triangleright \Gamma_2
}
\and
\infer{\Gamma\vdash\eif\;b\;\ethen\;S_1 \;\eelse\;S_2\;\eendif
\;\triangleright \Gamma}{%
\Gamma \vdash b : \dbool &
\Gamma \vdash S_1 \triangleright \Gamma_1 &
\Gamma \vdash S_2 \triangleright \Gamma_2
}
\and
\infer{\Gamma \vdash \ewhile\;b\;\edo\;S\;\eendwhile
\;\triangleright \Gamma}{%
\Gamma \vdash b : \dbool &
\Gamma \vdash S \triangleright \Gamma'
}\smallskip
\end{mathpar}
\end{minipage}}
\end{figure*}

The only construct that modifies a context
is $\elet\;x:\tau=t$ (variable declaration). When it is
executed under a context $\Gamma$, the command is well-formed
if $x$ is not already included in $\Gamma$ and the type of
the initializing term $t$ is of the declared type $\tau$.
After the execution, we get the new context
$\Gamma, x : \tau$.
However, when a context changes inside a branch or a loop, it gets restored once the block is finished.
In other words, the variables created inside of a branch or a loop only survive locally
(as common for example in C++).

Here we state a property of our typing rules that is needed later.
We omit detailed proof as it can be seen directly by structural induction on the type judgement.
\begin{lem}
\label{l:typing}
Whenever $\Gamma \vdash S \triangleright \Gamma'$, 
the new context $\Gamma'$ is an extension of $\Gamma$ in the sense that 
there is a context $\Delta$ such that $\Gamma' = \Gamma, \Delta$.
\end{lem}

\subsubsection{Well-Typed Programs}
\label{sss:Type_Program}
Consider the two different kinds of \ERC programs:
\[
\begin{array}{rl}
\prog_{\dint} \, \dfeq& \einput\;x_1 : \tau_1, x_2 : \tau_2, \cdots, x_{d} : \tau_d \\
& \quad S \\
& \ereturn\;t
\end{array}
\quad
\begin{array}{rl}
\prog_{\dreal} \, \dfeq  \,\,
& \blim{p} x_1 : \tau_1, x_2 : \tau_2, \cdots, x_d : \tau_d \\
& \quad S \\
& \blimreturn{t}{p}
\end{array}
\]
The first program $\prog_{\dint}$ has type $\tau_1 \times \cdots \times \tau_d \to \dint$ if there is a context $\Gamma'$ such that
\[
   x_1:\tau_1, x_2 : \tau_2, \cdots, x_d : \tau_d \vdash S \triangleright \Gamma'
   \quad\text{and}\quad \Gamma' \vdash t : \dint \,.
\]
Similarly, $\prog_{\dreal}$ has type $\tau_1 \times \cdots \times \tau_d \to \dreal$ if there is a context $\Gamma'$ such that
\[
   p:\dint, x_1 : \tau_1, \cdots, x_d : \tau_d \vdash S \triangleright \Gamma'
   \quad\text{and}\quad \Gamma' \vdash t : \dreal \,.
\]
Note that $p : \dint$ is not regarded as an input variable and 
recall that $\tau_i$ varies over $\{\dint, \dreal,
\dreal[1], \dreal[2], \cdots\}$ and cannot be $\dbool$.

We call an \ERC program an integer program if it is of the first kind ($\prog_{\dint}$) and a real program if it is of the second kind ($\prog_{\dreal}$).

\section[Denotational Semantics of ERC]{Denotational Semantics of \emph{Exact Real Computation}}
\label{s:Semantic}

We define multi-valued semantics for well-typed terms (Section~\ref{sss:Semantic_Term}), well-formed commands (Section~\ref{sss:Semantic_Command}), and well-typed programs (Section~\ref{sss:Semantic_Program}).
In this semantics, the objects of \ERC are assigned mathematical meanings that are arguably (i) closest possible to the intuition of real numbers as entities to be operated on exactly while simultaneously featuring (ii) Turing-completeness.

We start with denotations for data types and the definition of states.
\begin{defi}
\hfill
\begin{enumerate}
    \item Data types are interpreted as intended:
        \[
            \sem{\dbool} = \Kleene,\quad
            \sem{\dint} = \mathbb{Z},\quad
            \sem{\dreal} = \mathbb{R},\quad
            \sem{\dreal[n]} = \mathbb{R}^n.
        \]
    Recall that $n$ above is a natural number greater than $0$ which is not to be confused with an integer-typed term in $\ERC$.
     \item Contexts are interpreted as sets of assignments.
     Formally,
        given a context $\Gamma = x_1 : \tau_1, \ldots, x_n : \tau_d$,
        then
        \[\sem{\Gamma} \dfeq  
            \sem{\tau_1} \times \cdots \times \sem{\tau_d}
        .\]

    \item Given a context $\Gamma$, an element $\sigma\in\sem{\Gamma}$ is called
        \inlinedef{state} which is a specific assignment of variables
        contained in the domain of the context.
        \qed
\end{enumerate}
\end{defi}

\subsection{Denotations of Terms}
\label{sss:Semantic_Term}
A (well-typed) term's meaning under a state is, considering the nondeterminism in \ERC, the set of all possible values that the term can evaluate.
For example, a well-typed term $\Gamma \vdash t : \dreal$ can evaluate to multiple values in $\IR$ under a state $\sigma \in \sem{\Gamma}$.
Moreover, if $0$ is among these values, then the compound term $t^{-1}$ could be undefined in addition to its defined values derived from non-zero values of $t$. 
This is reflected by including $\bot$ in the denotation, which is thus an element of $\pdom{\sem{\dreal}_\bot}$.

\begin{defi}
\label{defi:STerms}
Recall the computable partial multi-valued functions specified in Examples~\ref{exa:primitive1} and Examples~\ref{exa:primitive2}.
Given a well-typed term $t$ such that $\Gamma\vdash t : \tau$,
we define its denotation as a partial multi-valued function
$\sem{\Gamma\vdash t : \tau} : \sem{\Gamma} \to \pdom{\sem{\tau}_\bot}$
inductively as follows:
\begin{align*}
\sem{\Gamma \vdash k : \dint} & \dfeq \; 
    \sigma \mapsto \{ k \}\\
\sem{\Gamma \vdash \ttrue : \dbool} & \dfeq \; 
    \sigma \mapsto\{ \ctt \}\\
\sem{\Gamma \vdash \tfalse : \dbool} & \dfeq \; 
    \sigma \mapsto\{ \cff \}\\
\sem{\Gamma \vdash \tundef : \dbool} & \dfeq \; 
    \sigma \mapsto\{ \cuu \}\\
\sem{\Gamma\vdash [t_1, \cdots, t_n] : \dreal[n]} &\dfeq \;
    \id^{\ddagger^\dagger}
    \fcomp \langle \sem{\Gamma \vdash t_1 : \dreal}, \cdots, \sem{\Gamma \vdash t_n : \dreal} \rangle\\
\sem{x_1:\tau_1, \cdots, x_i :\tau_i,\cdots,x_d : \tau_d \vdash x_i : \tau_i} & \dfeq \; 
    \semproj_i^{\ddagger} \\
\sem{\Gamma \vdash t[m] : \dreal} & \dfeq \;
    \semproj^{\ddagger^\dagger} \fcomp \langle \sem{\Gamma \vdash t : \dreal[n] }, \sem{\Gamma \vdash m : \dint}\rangle \\
\sem{\Gamma \vdash t_1 \bop t_2 : \tau'} & \dfeq \;
    \bop^{\ddagger^\dagger} \fcomp \langle \sem{\Gamma \vdash t_1 : \tau}, \sem{\Gamma \vdash t_2 : \tau} \rangle \\
& \text{for } (\bop, \tau, \tau') \in \bopty\\
\sem{\Gamma \vdash \uop t : \tau'} & \dfeq \;
    \uop{}^{\ddagger^\dagger} \fcomp \sem{\Gamma \vdash t : \tau}\\
& \text{for } (\uop, \tau, \tau') \in \uopty\\
\sem{\Gamma \vdash f(t_1, \ldots, t_d) : \dreal} & \dfeq \;
f^{\ddagger^\dagger}\fcomp \langle\sem{\Gamma \vdash t_1 : \otype{X_1}}, \cdots, \sem{\Gamma \vdash t_d : \otype{X_d}}\rangle \\
    &\text{where }(f : X_1\times\cdots\times X_d \to \IR_\bot)\in\calF\\
\sem{\Gamma \vdash g(t_1, \ldots, t_d) : \dint} & \dfeq \;
    g^\dagger\fcomp \langle\sem{\Gamma \vdash t_1 : \otype{X_1}}, \cdots, \sem{\Gamma \vdash t_d : \otype{X_d}}\rangle \\
    &\text{where }(g : X_1\times\cdots\times X_d \to \pdom{\IZ_\bot})\in\calG\\
\sem{\Gamma \vdash \mychoose_n(b_0, \cdots, b_{n-1}) : \dint} & \dfeq \;
    \semchoose_n^\dagger \fcomp \langle
\sem{\Gamma \vdash b_0:\LOGIC}, \cdots,\sem{\Gamma \vdash b_{n-1}:\LOGIC}\rangle \\
\sem{\Gamma \vdash (\mycond{b}{u}{v}) : \dreal}  & \dfeq \;
    \semkond^{\dagger_1} \fcomp \langle\sem{\Gamma \vdash b :\LOGIC}, \sem{\Gamma \vdash u : \REAL}, \sem{\Gamma \vdash v : \REAL}\rangle
\end{align*}
\qed
\end{defi}
\noindent We often simply write $\sem{t}$ instead of $\sem{\Gamma\vdash t:\tau}$ omitting $\Gamma$ and $\tau$ when they are obvious or irrelevant.
We also synonymously say
$t$ \emph{evaluates} to
$x$ under $\sigma$,
or $t$ \emph{has/contains} the element $x\in\eval{t}\sigma$.

Regarding the multi-valuedness, $\mychoose()$ and functions from $\calG$
are the only atomic constructs that yield multi-valuedness.
The denotations of all other constructs except the two and continuous conditionals
are defined by the embedding of single-valued mappings.
For the case of continuous conditionals, we can observe that $\semkond^{\dagger_1}(p, q, r)$ is a singleton when $p$, $q$, and $q$ are singletons.
Hence, it does not generate multi-valuedness.

Another remark on the multi-valuedness is that the denotation of integer equality $m_1  \ieq  m_2$ does not coincide with that of $m_1 \leq m_2 \kland m_2 \leq m_1$. 
As an example, consider a state $\sigma$ where $\sem{m_1}\sigma = \{0,2\}$ and $\sem{m_2}\sigma = \{1, 3\}$.
    Then, $\sem{m_1 \leq m_2 \kland m_2 \leq m_1}\sigma = \{\ctt, \cff\}$ 
    whereas $\sem{m_1 \ieq m_2}\sigma = \{\cff\}$. This justifies 
    introducing $m_1 \ieq m_2$ as a separate primitive.

Well-typedness in $\ERC$ does not prevent $\bot$: The sources of $\bot$ are (i) an array accessed by a wrong index ($\semproj$), (ii) the multiplicative inverse of $0$ ($0^{-1}$), (iii) function calls to partial (multi-valued) functions at points out of their domains, (iv) a choice operation without any choosable argument ($\semchoose$), and (v) a continuous conditional when the two arguments $u, v$ do not match in $\semkond(\cuu, u, v)$ as $u = v = \{x\}$ for some $x \in \IR$. 
We emphasize here again that the computational meaning of $\bot$ is not non-termination but is more general non-specified behaviour; recall Section~\ref{s:preliminaries}.

Though $\pdom{X_\bot}$ allows infinite sets (as long as they contain $\bot$), this case does not occur in our term language.
The denotation of any well-typed term at any state is finite unless the case is provided by functions in $\calG$.
However, we still use the powerdomain for our term language to use the same setting with the denotational semantics of our command language where indeed $\bot$-containing infinite sets arise.

The computability of the denotations follows directly from the computability of the primitive operations and of the compositions stated in Examples~\ref{exa:primitive1} and Examples~\ref{exa:primitive2}:
\begin{lem}
\label{l:Computable_Term}
The denotation of any well-typed term is a computable partial multi-valued function.
\end{lem}

\subsection{Denotations of Commands}
\label{sss:Semantic_Command}
As the denotation of a well-typed term under a state
is the set of all possible values that the term
may evaluate to, considering multi-valuedness in \ERC,
the denotation of a well-formed command under a state is
the set of all possible resulting states that the command may result in.
Hence, we let a well-typed command denote a function from the set of states to the restricted power-set of the resulting states:

\begin{defi}
\label{defi:SCommands}
Recall the (multi-valued) functions specified in Examples~\ref{exa:primitive1} and Examples~\ref{exa:primitive2}.
Given a well-formed command $S$ such that $\Gamma\vdash S \triangleright \Gamma'$,
we define its denotation to be a partial multi-valued function
$ \sem{\Gamma \vdash S \triangleright \Gamma'} :
        \sem{\Gamma} \to \pdom{\sem{\Gamma'}_\bot}$
inductively as follows:
\begin{align*}
\sem{\Gamma\vdash \eskip \triangleright \Gamma} & \dfeq \;
    \id^\ddagger \\
\sem{\Gamma  \vdash x_i\dfeq  t \triangleright \Gamma}  & \dfeq \;
    \semupdate_{i}^{\ddagger^\dagger} \fcomp \langle \id^\ddagger, \sem{t}\rangle \\
    &\text{ where } \Gamma = x_1:\tau_1, \cdots, x_i : \tau_i, \cdots, x_d : \tau_d \\
\sem{\Gamma \vdash x_i[m]\dfeq  t \triangleright \Gamma}  & \dfeq \;
    \semupdate_{i}^{\ddagger^\dagger} \fcomp \langle \id^\ddagger, \semupdate^{\ddagger^\dagger} \fcomp \langle 
    \semproj_i^\ddagger, \sem{m}, \sem{t}\rangle \rangle\\
    &\text{ where } \Gamma = x_1:\tau_1, \cdots, x_i : \REAL[n], \cdots, x_d : \tau_d\\
\sem{\Gamma \vdash \elet\;x : \tau = t\; \triangleright \Gamma'}  & \dfeq \;
    \id^{\ddagger^\dagger}\fcomp \langle \id^\ddagger, \sem{t}\rangle \\
\sem{\Gamma \vdash S_1; S_2 \triangleright \Gamma'} & \dfeq \;
    \sem{\Delta \vdash S_2 \triangleright \Gamma'}^\dagger \fcomp \sem{\Gamma \vdash S_1 \triangleright \Delta}
\\
   \sem{\Gamma \vdash\eif\;b\;\ethen\;S_1\; \eelse\;S_2 \;\eendif \triangleright \Gamma} & \dfeq \!
    \semcond^{\dagger_1}\fcomp\langle \sem{b}, \restriction_\Gamma^{\ddagger^\dagger} \fcomp \sem{\Gamma \vdash S_1\triangleright \Gamma'}, \restriction_\Gamma^{\ddagger^\dagger} \fcomp \sem{\Gamma \vdash S_2\triangleright \Gamma''}\rangle 
\\
\sem{\Gamma \vdash \ewhile \;b\;\edo \;S\;\eendwhile \triangleright \Gamma}
& \dfeq\;
\lfp\big(\wop(\sem{b}, \restriction_\Gamma^{\ddagger^\dagger} \fcomp\sem{ \Gamma \vdash S \triangleright \Gamma'})\big) 
\end{align*}
Here, $\restriction_{\Gamma} : \sem{\Gamma, \Delta} \to \sem{\Gamma}$ for any contexts $\Gamma, \Delta$ 
when $\Gamma = x_1 :\tau_1, \cdots, x_d : \tau_d$ is defined by
\[
(x_1, \cdots, x_d, y_1, \cdots, y_m) \mapsto (x_1, \cdots, x_d). 
\tag*{\qed}
\]
\end{defi}
\noindent We often write $\sem{S}$ to denote  $\sem{\Gamma\vdash S \triangleright \Gamma'}$ omitting $\Gamma$ and $\Gamma'$ when they are obvious or irrelevant.

When a term is assigned to a variable, an array variable, or a newly created variable as its initial value, 
the underlying state branches into multiple states for each single value in the denotation of the term. 
Instead of having a single state storing multi-values, we choose to have multiple states where each stores single-values.

The denotations of conditional statements and loops are well-defined due to Lemma~\ref{l:typing}
which justifies the post-composed restriction function $\restriction_{\Gamma}$ clearing the newly added local variables. 
Regarding loops and branches,
note that their denotations contain $\cuu$ when their conditions evaluate to $\cuu$.
We could also consider replacing $\semcond$ with $\semkond$ for the semantics of conditionals. However, first, their domains do not match: $\semkond$ requires the branches to be in $\pdom{\IR_\bot}$ whereas we want them to be in $\pdom{\sem{\Gamma}_\bot}$. Though one could imagine extending $\semkond$'s domain, recalling how $\semkond$ is realized in Examples~\ref{exa:primitive2}, this requires the resulting states of both branches to be inspected variable-wise.
For this work, we decide to keep the ordinary branching which makes a decision only on its condition and diverges in the case the condition diverges by being $\cuu$ instead of further inspecting the two resulting states.

Unlike for terms, the multi-valued computability of loops 
needs more attention:

\begin{prop}
\label{p:loop}
For a represented set $X$ and computable partial multi-valued functions 
$b : X \to \pdom{\IK_{\bot}}$ and $S : X \to \pdom{X_{\bot}}$,
the least fixed point
\[
\lfp\big(\wop(b, S)\big) : X \to \pdom{X_{\bot}}
\]
as partial multi-valued function is computable uniformly to $b$ and $c$
in the sense that the realizer for
$\lfp\big(\wop(b, S)\big) : X \to \pdom{X_{\bot}}$
can be defined with function calls to the realizers of $b$ and $S$
for any $b$ and $S$.
\end{prop} 
\begin{proof}
By the least fixed point theorem, the least fixed point is  
the limit of the chain
\[
\mathcal{W}^{(0)}_{b, S} \dfeq \sigma \mapsto\{\bot\}
\quad\text{and}\quad 
\mathcal{W}^{(n+1)}_{b, S} \dfeq 
\semcond^{\dagger_1}\fcomp\langle b, \mathcal{W}^{(n)}_{b, S} \fcomp S, \id^\ddagger\rangle
\]
which, intuitively, represents the process of unrolling the loop. 
For any two $\sigma, \delta \in X$,
$\delta\in\lfp\big(\wop(b, S)\big) \sigma$ if and only if there is a natural number $n$ where
$\delta \in \mathcal{W}^{(n)}_{b, S}\sigma$.
Furthermore, $\bot\in\lfp\big(\wop(b, S)\big) \sigma$ if and only if 
$\bot \in \mathcal{W}^{(n)}_{b, S}\sigma$ for all $n$.

Suppose $\tau_{b}$ realizes $b$ and $\tau_{S}$ realizes $S$.
We consider the type-2 machine $\tau_{\tau_{b}, \tau_{S}}(\phi_{\sigma})$ 
on its input $\phi_{\sigma}$ a name of some $\sigma \in X$ with the following description:
\begin{lstlisting}[language=typetwo,mathescape=true]
input $\phi_\sigma$
    repeat:
	   $\phi_b \dfeq \tau_b(\phi_{\sigma})$
	   iterate $i := 0, 1, 2, \cdots$ until $\phi_b(i) = -1$ or $\phi_b(i) = 1$
	   if $\phi_b(i) = 1$: 
		  $\phi_\sigma \dfeq \tau_S(\phi_{\sigma})$
	   if $\phi_b(i) = -1$: 
		  return $\phi_\sigma$
\end{lstlisting}
We need to prove that if $\bot\not\in\lfp\big(\wop(b, S)\big) \sigma$, 
$\tau_{\tau_{b}, \tau_{S}}(\phi_{\sigma})$ is defined and is a name 
of some $\delta \in \lfp\big(\wop(b, S)\big) \sigma$.
(Recall that we do not need to reason about $\tau_{\tau_{b}, \tau_{S}}(\phi_{\sigma})$
when $\bot\in\lfp\big(\wop(b, S)\big) \sigma$.)

Consider the indexed variants $\tau^{(n)}_{\tau_{b}, \tau_{S}}(\phi_{\sigma})$ described by:
\begin{lstlisting}[language=typetwo,mathescape=true]
input $\phi_\sigma$
    repeat:
        if n = 0 then diverge else n := n - 1
        $\phi_b \dfeq \tau_b(\phi_{\sigma})$
        iterate $i := 0, 1, 2, \cdots$ until $\phi_b(i) = -1$ or $\phi_b(i) = 1$
        if $\phi_b(i) = 1$: 
            $\phi_\sigma \dfeq \tau_S(\phi_{\sigma})$
        if $\phi_b(i) = -1$: 
            return $\phi_\sigma$
\end{lstlisting}
which is almost identical to $\tau_{\tau_{b}, \tau_{S}}(\phi_{\sigma})$
except that $\tau^{(n)}_{\tau_{b}, \tau_{S}}(\phi_{\sigma})$ 
repeats the main iteration only at most $n$ times and diverges 
unless it returns until then.
By induction on $n$, we can easily verify that for each natural number $n$, 
$\tau^{(n)}_{\tau_{b}, \tau_{S}}$ realizes 
$\mathcal{W}^{(n)}_{b, S}$.

Furthermore, it holds that as long as $\tau^{(n)}_{\tau_{b}, \tau_{S}}(\phi_{\sigma})$
is defined, it coincides with $\tau_{\tau_{b}, \tau_{S}}(\phi_{\sigma})$.
Since $\bot\not\in\lfp\big(\wop(b, S)\big) \sigma$, 
there exists $n$ where $\lfp\big(\wop(b, S)\big)\sigma = \mathcal{W}^{(n)}_{b, S}\sigma \not\ni \bot$
and $\tau^{(n)}_{\tau_{b}, \tau_{S}}(\phi_{\sigma})$ is a name of some $\delta \in \mathcal{W}^{(n)}_{b, S}\sigma$.
Therefore, by the above claim, $\tau_{\tau_{b}, \tau_{S}}(\phi_{\sigma})$
coincides with $\tau^{(n)}_{\tau_{b}, \tau_{S}}(\phi_{\sigma})$
which is a name of $\delta \in \lfp\big(\wop(b, S)\big)\sigma$.
\end{proof}

In the lemma, we say ``computable uniformly to $b$ and $c$'' to denote that the realizer can be defined 
by function calls to the realizers of $b$ and $c$ without depending on any further specifics of (the realizers of) $b$ and $c$. 
We put a justification for this notion of uniformity. 
Type-2 computability provides a utm theorem stating that 
there exists an indexing of type-2 machines by $\IZ^\IN$ and 
a universal type-2 machine for that indexing.
Hence, the above uniformity means that the realizer can be actually 
defined as a multivariate function that receives the indices for the realizers of $b$ and $c$ as its input where the function calls are replaced with compositions by the universal machine.
In other words, the least fixed point of 
$\wop(b, c)$ 
is computable uniformly in the following sense that the operation
\[
\text{indices of the realizers of }b, c \text{ and a name of }\sigma
\mapsto
\text{a name of }\delta \in\lfp\big(\wop(b, c)\big)\sigma \]
is computable for any $b, c, \sigma$.
We can also make this notion formal by representing computable partial multi-valued functions and making the above operation a computable mapping from them using the same type-2 machine. However,  since the cardinality of the set of computable partial multi-valued functions exceeds that of $\IZ^\IN$, it requires us to use a more general notion of representation, namely, multi-representation \cite{schroder2002effectivity} or assemblies over Kleene's second algebra.

Following up on Lemma~\ref{l:Computable_Term} and Proposition~\ref{p:loop}, we establish the computability result for our command language:
\begin{lem}
\label{l:Computable_Command}
The denotation of any well-formed command is a computable partial multi-valued function. 
\end{lem}
\begin{proof}
The restriction operations are computable.
Due to the computability of the primitives, compositions, and fixed points of loops,
the denotations of well-formed commands are computable partial multi-valued functions.
\end{proof}

%

\subsection{Denotations of Programs}
\label{sss:Semantic_Program}
Having defined the meanings of terms and commands, we are now ready to define denotations of well-typed \ERC programs.

\begin{defi} 
\label{defi:Sprograms}
The denotation of a well-typed integer program
\[
\begin{array}{rl}
\prog_{\dint} \, \dfeq& \einput\;x_1 : \tau_1, x_2 : \tau_2, \cdots, x_d : \tau_d \\
& \quad S \\
& \ereturn\;m
\end{array}
\]
is
a partial integer multi-valued function 
\[
    \sem{\prog_{\dint}} :\sem{\tau_1}\times\cdots\times\sem{\tau_d}
    \to \pdom{\IZ_\bot}
\]
defined by $\sem{\prog_{\dint}} \dfeq \sem{m}^\dagger \fcomp \sem{S}$.

The denotation of a well-typed real program 
\[\begin{array}{rl}
\prog_{\dreal} \, \dfeq  \,\,
& \blim{p} x_1 : \tau_1, x_2 : \tau_2, \cdots, x_d : \tau_d \\
& \quad S \\
& \blimreturn{u}{p}
\end{array}
\]
is a partial real function 
$\sem{\prog_{\dreal}} :\sem{\tau_1}\times\cdots\times\sem{\tau_d}
        \to \IR_\bot$
    that is defined by 
    \[\sem{\prog_{\dreal}}(x_1, \cdots, x_d) \dfeq 
    \begin{cases}
    x & \text{if } \forall p\in\IZ.\;\forall y \in \sem{t}^\dagger \fcomp \sem{S}(p, x_1, \cdots, x_d).\; y \neq \bot \land |y - x| \leq 2^{p}, \\
    \bot &\text{otherwise. }
    \end{cases}
    \]
\end{defi}

In words, a real program $\prog_\dreal$ denotes a single-valued function $f$ whose function value $f(x_1, \cdots, x_d)$ is the limit of the multi-valued sequence $\sem{t}^\dagger \fcomp \sem{S}(p, x_1, \cdots, x_d)$ as $p \to -\infty$ in the sense that for each $i\in\IN$, whichever $y_i$ is picked from $\sem{t}^\dagger \fcomp \sem{S}(-i, x_1, \cdots, x_d)$ nondeterministically, $y_i$ is a $2^{-i}$ approximation to $f(x_1, \cdots, x_d)$ and not $\bot$.
The obligated rate of convergence $2^p$ is chosen to make it coincide with the precision embedding $\imath : \IZ \ni p \mapsto 2^p \in \IR$.

\begin{thm}
\label{t:computable}
The denotation of a well-typed real \ERC program is a computable partial real function. 
The denotation of a well-typed integer \ERC program is a computable partial  integer multi-valued function. 
\end{thm}
\begin{proof}
The computability of terms, commands, compositions,
and limit operations \cite[Theorem~4.3.7]{Wei00} yields the computability of the denotations of programs.
\end{proof}

When a real program on its input does not generate a converging sequence as its precision parameter $p$ heads to $-\infty$, its function value is denoted by $\bot$. This is the main motivation for us to make $\bot$ to represent
not only divergence or non-termination but more general non-specified behaviour.
When a sequence is not converging, operationally, there is no possible way to decide that the sequence is not converging. 
Hence, after reading some finite portion of the sequence, believing that the rest of the sequence converges, an approximation to the limit has to be made. 
That means even when we realize later that the sequence does not follow our expectations and decides not to converge,
there can be an output printed already as the limit of the sequence which turns out to be completely meaningless as the sequence does not converge after all.
Without giving this freedom of a procedure returning completely meaningless results on its non-valid inputs, limit operations cannot be realized.

Recall from Fact~\ref{fact:prelim 1} that in computable analysis, a machine realizing a (multi-valued) function can be seen to realize some other (multi-valued) functions at the same time.
Reflecting this property, we add one more layer to the denotations of programs:
\begin{defi}\hfill
\label{d:Realizer}
\begin{enumerate}
    \item 
    A well-typed real program $\prog$ \emph{realizes} a partial function $f : X_1\times\cdots\times X_d \to \IR_\bot$ if $\forall (x_1, \cdots, x_d).\; f(x_1, \cdots, x_d) \leq \sem{\prog}(x_1, \cdots, x_d)$ holds.
    Recall that $x \leq y$ if and only if $x = \bot$ or $x = y$ for $x, y \in \IR_\bot$.
    \item
    A well-typed integer program $\prog$ \emph{realizes} a partial (possibly infinite) multi-valued function $f : X_1\times\cdots\times X_d \to \pset{\IZ_\bot}$
    if $\forall (x_1, \cdots, x_d).\; f(x_1, \cdots, x_d) \sqsubseteq \sem{\prog}(x_1, \cdots, x_d)$ or $\sem{\prog}(x_1, \cdots, x_d) \subseteq f(x_1, \cdots, x_d)$ holds.
\end{enumerate}
\end{defi}
Note that though the denotations of integer programs are finite multi-valued functions, they can still realize \emph{infinite} multi-valued functions.
Directly from Fact~\ref{fact:prelim 2}, we get the computability of realizable partial (multi-valued) functions:
\begin{cor}
Realizable partial (multi-valued) functions are computable partial (multi-valued) functions.
\end{cor}

\section{Programming in \emph{Exact Real Computation}} \label{s:Programming}

In this section, we collect several examples of (multi-valued) functions
that are realizable in $\ERC$. 
They illustrate the purpose of $\ERC$ which is to allow and justify \naively implementing numerical algorithms, with computable analysis as the rigorous but hidden theoretical foundation.
They demonstrate how naturally the classical algorithms can be modified to incorporate the Kleenean-valued comparisons and the nondeterministic $\mychoose$.

The examples are square root using Heron's method, exponential function via two different ways, integer rounding, matrix determinants via Gauss elimination, and root finding. 
Heron's method, as an intuitive and efficient way to compute square roots of real numbers, is often expressed in other frameworks \cite{iRRAM,Mueller18,DBLP:conf/fsttcs/Konecny0T20}.
Through the example, the readers can easily compare $\ERC$ with the others; similarly, for Gauss elimination \cite[Section~2.2]{TZ04}.
The example of realizing the exponential function shows that the base language $\ERC_0$ can realize transcendental functions; this leads to a more detailed discussion of our specification language later in Section~\ref{s:Logic}.
Integer rounding, on the other hand, shows that $\ERC_0$ can nondeterministically extract names of real numbers.
This later forms a core routine in the proof of Turing-completeness in Section~\ref{ss:Complete}.
Root finding, a computational version of the intermediate value theorem is often picked as an example in computable or constructive analysis. 
The example further demonstrates a trick of using the extension sets $\calF$ and $\calG$ to model function calls to arbitrary external functions.
This example is chosen to be formally verified later in Section~\ref{ss:Trisection}.

The examples have been implemented using a shallow embedding of $\ERC$ into Haskell which we develop on top of the AERN library for exact real number computation \cite{Aern}.\footnote{The implementation can be found in \url{https://github.com/michalkonecny/aern2/blob/master/aern2-erc/src/ERC/Examples.hs}\;.} 
There, though internally real numbers and Kleeneans are represented by infinite sequences, they are hidden from the users where the users see them as abstract entities. Hence, without a nondeterminism monad, it provides the nondeterministic $\mychoose$ operator which can choose different indices for the same inputs whose internal sequences are different. In other words, the $\mychoose$ operator is deterministic about the internal sequences of the arguments which are inspected with some fixed interleaving procedure, but is nondeterministic when the internal sequences are abstracted away.  

\subsection{Programming Abbreviations}
Several abbreviations are used in our 
example programs. While deferring their introductions 
to each example program, here we make some remarks for those used importantly throughout this paper.

\subsubsection{Extension by denotations}
Recall that the programming language $\ERC(\calF, \calG)$ is defined relative to the sets $\calF$ and $\calG$. 
Beyond modelling primitive operator extensions, the sets $\calF$ and $\calG$ provide a trick on reasoning user-defined function calls (but without recursive calls) from some richer exact real computation software. 
Instead of equipping $\ERC$ with a feature of calling programs within programs, which removes $\ERC$'s strict barrier between its term language and command language, in order to keep our term language as simple as possible, when the users want to model a program calling another program, they can use $\ERC(\calF,\calG)$ with $\calF$ or $\calG$ extended with the denotations of the programs.

When we want to extend the set $\calF$ with the denotation of a real program $\prog$ in $\ERC(\calF, \calG)$, we write $\ERC(\calF\cup\{\prog\}, \calG)$ to denote $\ERC(\calF \cup\{\sem{\prog}\}, \calG)$.
Similarly, we write $\ERC(\calF, \calG\cup\{\prog\})$
when $\prog$ is an integer program in $\ERC(\calF, \calG)$.
These notations are well-defined since the denotation of a real program is a computable partial real function and the denotation of an integer program is a computable partial integer multi-valued function.
When the sets are extended by the denotation of a program $\prog$, we write $\prog(t_1, \cdots, t_d)$ for $\sem{\prog}(t_1, \cdots, t_d)$ in the extended term language.

\subsubsection{Variable scope}
\ERC allows local variable creations inside conditional and loop statements.
This mechanism for local variables can be used to define arbitrary variable scopes.
For a well-formed command $\Gamma \vdash S \triangleright \Gamma'$, write
\[
\{S\}
\]
for 
\[
\Gamma \vdash \eif \; \ctt \; \ethen \;S \; \eelse \; \eskip \triangleright \Gamma\;
\]
We can easily check:
\[
\restriction_\Gamma^{\ddagger^\dagger} \fcomp\sem{\Gamma \vdash S \triangleright \Gamma'} = \sem{\Gamma \vdash \{S\}
\triangleright \Gamma}.
\]
\subsubsection{Term Abbreviations}
We assume obvious term abbreviations such as $t_{1}/t_{2}$ for $t_{1}\times t_{2}^{-1}$, 
$t_{1} - t_{2}$ for $t_{1} + (- t_{2})$, and so on.
For each integer constant $k$, we write $k$ to refer to the $\myabs{k}$ times repeated additions of either $\myprec{0}$ or $-\myprec{0}$
when we want it to represent the real number $k$.

\subsubsection{Command Abbreviations}
Though our default term language is restrictive without having integer multiplications or the ordinary integer-to-real coercion, 
the default command language can recover most of the restrictions.
For instance, given two integer terms $m_{1}$ and $m_{2}$, 
though we cannot have its multiplication as a term, 
we can have a command that 
computes the multiplication and assigns the result into a variable say 
$y : \dint$:
\def\baselinestretch{1.15}\selectfont%
{\fontsize{10}{15}
\begin{algorithmic}[0] 
\Myscope
\State $\LET x_1 : \dint = m_1;\; \LET x_2 : \dint = m_2;\; \LET k : \dint = 0;$
\State $y \dfeq 0;$
\State \IF\; $0 \leq x_2$
\State \THEN\; \WHILE $k + 1 \leq x_2$\; \DO\; $y \dfeq y + x_1;\; k\dfeq k + 1$\; \ENDWHILE 
\State \ELSE\;\;\, \WHILE $x_2 + 1 \leq k$\; \DO\; $y \dfeq y - x_1;\; k\dfeq k - 1$\; \ENDWHILE 
\State\ENDIF
\EndMyscope
\end{algorithmic}
}
The local variables $x_{1}$ and $x_{2}$ are used to choose values from possibly multi-valued $m_{1}$ and $m_{2}$ and to use them consistently throughout the computation.
For a context $\Gamma$ such that $\Gamma \vdash m_{1} : \dint$, $\Gamma \vdash m_{2} : \dint$, and 
$\Gamma(y) = \dint$, 
\[
\Gamma \vdash y \abbeq m_{1} \times m_{2} \triangleright \Gamma
\]
abbreviates the above well-formed command whose local variables are selected to not conflict with $\Gamma$.
The notation $\abbeq$ is chosen to make it be distinguished from the usual assignment construct $\dfeq$ to emphasize that $m_{1} \times m_{2}$ is not a term. 
The following semantic equation can be verified:
\[
\sem{y \abbeq t_{1} \times t_{2}}\sigma = \bigcup_{x_{1} \in \sem{t_{1}}\sigma \land x_{2} \in \sem{t_{2}}\sigma}
\begin{cases}
\{\bot\} &\text{if } x_{1} = \bot \lor x_{2} = \bot, \\
\{\sigma[y \mapsto x_{1} \times x_{2}]\} &\text{otherwise.}
\end{cases} 
\]

We can define various command abbreviations similarly for example
\[
\Gamma \vdash y \abbeq u^m \triangleright \Gamma
\] 
when $\Gamma \vdash u : \dreal$, $\Gamma \vdash m : \dint$, and $\Gamma(y) = \dreal$
for the repeated multiplications satisfying
\[
\sem{y :\equiv u^m} \sigma = \bigcup_{x_{1} \in \sem{u}\sigma \land x_{2} \in \sem{m}\sigma}
\begin{cases}
\{\bot\} &\text{if } x_{1} = \bot \lor x_{2} = \bot, \\
\{\sigma[y \mapsto x_{1}^{x_{2}}]\} &\text{otherwise.}
\end{cases} 
\]

\subsection{Example Programs} \label{ss:Examples}
We annotate assertions next to commands to convey the behaviours, including the domains (preconditions), of the introduced programs.
However, we remind the readers that the annotations are not part of $\ERC$.

\subsubsection{Square Root Function via Heron's Method}
\label{ss:HeronSqrt}

Heron's method approximates the square root of a given real number $0 \leq x$
up to any desired absolute error $2^p$, $p\in\IZ$, by calculating a contracting sequence of upper and lower approximations $z_{n} \dfeq
x/y_n \leq \sqrt{x}  \leq y_n$ iteratively taking the average $y_{n+1}\dfeq (y_n+z_n)/2$.
The $\ERC_0$ program $\HeronSqrt$ illustrated in \autoref{fig:HeronSqrt1} realizes the square root function via Heron's method.

\begin{algorithm}
\caption{$\HeronSqrt : \REAL \to \REAL$ \label{fig:HeronSqrt1}}

\def\baselinestretch{1.15}\selectfont%
{\fontsize{10}{15}
\begin{algorithmic}[0] 
\Limit{p}{x:\dreal}
 \RightComment {$0 \leq x$}
\State $\LET y:\REAL=1;\;\; \LET z:\REAL=x/y$;
\RightComment {$z \leq \sqrt{x} \leq y$}
\Mywhile{$\mychoose\big(y-z\elt \myprec{p} \;,\; \myprec{p-1} \elt y-z\big) \ieq 1$}
    \State $y \dfeq  (y+z)/2$;\;\; $z\dfeq  x/y$
    \RightComment {$z \leq \sqrt{x} \leq y$}
\EndMywhile
\EndLimit{y}{p}  \RightComment  {$|y-\sqrt{x}|<2^p$}
\end{algorithmic}
}
\end{algorithm}

Recall that the arguments to $\mychoose$ are indexed from $0$,
hence $\mychoose(b_0,b_1)=1$ means that the second argument $b_1$ must evaluate to $\ctt$.
The \ERC program employs the multi-valued $\mychoose()$ operation applied to two Kleenean-valued tests $y - z \elt \myprec{p}$ and $ \myprec{p-1} \elt y-z$ where at least one of the two must be $\true$ resulting in a total loop condition.
Moreover, when the loop terminates, the first test (corresponding to return value $0$) must (while the second test, corresponding to return value $1$, may or may not) be $\true$:
guaranteed to return an approximation $y$ to $\sqrt{x}$ up to absolute error $2^p$.

\subsubsection{Exponential Function via Taylor Expansion}
\label{ss:Exp}

The exponential function has a globally converging Taylor expansion
\begin{equation}
\label{eqn:Exp}
\exp(x) \;=\; \sum\nolimits_{j=0}^\infty x^j/j!
\end{equation}

A \ERC real program computing it must return, given a dedicated precision parameter $p\in\IZ$ and argument $x\in\IR$, an approximation to $\exp(x)$ up to error $2^p$.
For $|x|\leq 1$ and positive $n\in\IN$, the tail bound
\begin{equation}
\label{eqn:Exp2}
\big|\sum\nolimits_{j>n} x^j/j!\big| \;\leq\;
\sum\nolimits_{j>n} |x|^j/j! \;\leq\;
\sum\nolimits_{j>n} 2^{-j+1} \;=\; 2^{-n+1}
\end{equation}
justifies the straightforward algorithm as in the \ERC program $\pExp$ in \autoref{fig:Exp1}.
The program introduces a command abbreviation
$\elet\ x_1, \cdots, x_d : \tau = t$
for $\elet\ x_1 : \tau = t; \cdots; \elet\ x_d : \tau = t$.
In the program, the two variables $j:\INTEGER$ and $j_r:\REAL$  retain identical values but have distinguished types, written as $j \equiv j_r$ in the assertions.

\begin{algorithm}[t]
\caption{$\pExp : \REAL \to \REAL$ \label{fig:Exp1}}
\def\baselinestretch{1.15}\selectfont%
{\fontsize{10}{15}
\begin{algorithmic}[0] 
\Limit{p}{x:\dreal}  \RightComment{$-1 \leq x\leq 1$}
\State $\LET j:\INTEGER=1; \;\; \LET j_r,\ f :\REAL=1;$
\RightComment{$j\equiv j_r, \;\; f\equiv j!$}
\State $\LET y:\REAL=1; \;\;\LET z:\REAL=x;$
\RightComment{$z\equiv x^j$}
\Mywhile{$j\leq -p+1$}
    \State $y \dfeq  y+z/f$;
    \State $j \dfeq  j+1$;
    \State $j_r \dfeq  j_r+1$;
    \State $z\dfeq  z\mult x$;
    \State $f \dfeq  f\mult j_r$
\EndMywhile
\EndLimit{y}{p}  \RightComment{$|y-\exp(x)|\leq 2^p$}
\end{algorithmic}
}
\end{algorithm}

The following wrapper $\ppExp$ for $\pExp$ in \autoref{fig:Exp}, which is a program in $\ERC(\{\pExp\}, \emptyset)$ 
removes the restriction $|x|\leq1$ and instead works for all  $-1 \leq x$.
Furthermore, since for $x < 0$ it holds that $\exp(x) = 1/\exp(-x)$, and $\exp(y) = 1/\exp(y)$ at $y = 0$, we can construct a procedure that computes $\exp(x)$ for all $x \in \IR$ using the continuous conditional by 
\[\mycond{0 \elt x}{\ppExp(x)}{1/\ppExp(-x)}.\]

\begin{algorithm}[t]
\caption{$\ppExp : \REAL \to \REAL$ \label{fig:Exp}}
\def\baselinestretch{1.15}\selectfont%
{\fontsize{10}{15}
\begin{algorithmic}[0] 
\Limit{p}{x:\dreal}  \RightComment{$-1\leq x$}
\State $\LET z:\REAL=\pExp(1/2);$
\State $\LET y:\REAL=1;$
\Mywhile{$\mychoose(x\elt 1, 1/2  \elt x ) \ieq 1$}
    \State $y \dfeq  y\mult z$;
    \State $x \dfeq  x-1/2$
\EndMywhile
\EndLimit{y \mult\pExp(x)}{p}
\end{algorithmic}
}
\end{algorithm}

\subsubsection[Realizing the Exponential Function by Iterations]{Exponential Function via Iteration}
\label{ss:Exp2}

$\pExp$ employs unbounded sums, which leave the realm of the first-order language considered in Section~\ref{s:Logic} for formal specification and verification purposes. 
Here we consider an alternative, iterative approach to realize the exponential function in \ERC.
To this end recall for every $x\in[0;2]$ it holds \cite[\S3.6.3]{Mitrinovic}:
\[
    \exp(x) \;\overset{\scriptscriptstyle n\to\infty}{\longleftarrow}\;
    \left (1+\tfrac{x}{n}\right )^n \;\leq\; \exp(x) \;\leq\;
    \left (1+\tfrac{x}{n}\right )^{n+1}
    \;\overset{\scriptscriptstyle n\to\infty}{\longrightarrow} \exp(x)
\]
%

\begin{algorithm}[t]
\caption{$\iExp :  \REAL \to \REAL$ \label{fig:Exp2}}
\def\baselinestretch{1.15}\selectfont%
{\fontsize{10}{15}
\begin{algorithmic}[0] 
\Limit{p}{x:\dreal} \RightComment{$0\leq x\leq2$}
\State $\LET n:\INTEGER=1$;\;
\State $\LET n_r:\REAL=1$;
\State $\LET c:\REAL=1+x$;\;
        $\LET a:\REAL=c$;\;
        $\LET b:\REAL=a\mult c$;
\Mywhile{$\mychoose\big(b - a \elt \myprec{p} , \myprec{p-1} \elt b-a  \big) \ieq 1$}
    \State $n \dfeq  n+n$;
    \State $n_r \dfeq  2\mult n_r$;
    \State $c \dfeq  1+x/n_r$;\;
            $a \abbeq  c^n$;\;
            $b \dfeq  a\mult c$
\EndMywhile
\EndLimit{a}{p}
\end{algorithmic}
}
\end{algorithm}

This suggests the iterative \ERC program $\iExp$ in \autoref{fig:Exp2}.
It supposes $0\le x \le 2$ but can
be extended to the entire real line similar to Section~\ref{ss:Exp}.
Note that the loop is guaranteed to terminate
since $b - a$ converges to $0$ as $n$ grows.
The loop condition is total:
At least one of $b-a\elt \myprec{p}$ and $\myprec{p-1} \elt b-a $ must be $\true$.
Furthermore, when the  loop terminates,
it must hold $b - a < 2^{p}$,
hence $|a - \exp(x)| \le |b - a| \le 2^p$.

\subsubsection{Integer Rounding Multi-valued function} \label{ss:Round}

All real-to-integer rounding \emph{functions} (up, down, to nearest)
are discontinuous and therefore not computable. 
We thus relax the specification and instead consider
the following multi-valued function with \emph{overlap}:

\begin{equation} \label{eqn:Round}
\Round:\IR\ni x\mapsto \{k\in\IZ \mid x-1<k<x+1\}\subseteq\IZ
\end{equation}

\medskip
The $\ERC$ program $\nround$ in \autoref{fig:Round1} realizes $\Round$.
It returns, given a real $x$, multi-valued integer $k$ such that $|x - k| < 1$.

\begin{algorithm}
\caption{$\nround : \REAL \to \INTEGER$ \label{fig:Round1}}
\def\baselinestretch{1.15}\selectfont%
{\fontsize{10}{15}
\begin{algorithmic}[1]
\Prog{$x:\REAL$}
\State $\LET k:\INTEGER=0$;
\Mywhile{$\mychoose(x\elt 1,1/2\elt x) \ieq 1$}
  \State $k \dfeq  k+1;\; x \dfeq  x-1$
\EndMywhile;
\Mywhile{$\mychoose(-1 \elt x,x\elt-1/2) \ieq 1$}
  \State $k \dfeq  k-1;\; x \dfeq  x+1$
\EndMywhile
\EndProg{$k$}
\end{algorithmic}
}
\end{algorithm}

Intuitively, the number of steps made by $\nround$ is proportional
to the value of the argument $x$, that is,
exponential in its binary ($\approx$output) length because $\nround$ essentially counts up to $x$.

Recovering the rounded integer bit-wise via binary search seems
exponentially more efficient, but fails due to lack of continuity:
Extracting any digit of the binary expansion
(or one of the at most two possible ones) of a given real number is uncomputable \cite{Tur37}.
Instead, $\fround$ in \autoref{fig:Round2} realizes the idea that
\emph{some signed-}digit expansion  \cite[Definition~7.2.4]{Wei00}
of a given real number $x\in\IR$ can be determined computably.
There, the absolute value $\myabs{y}$ is an abbreviation for 
$\mycond{0 \elt y}{y}{-y}$.

\begin{algorithm}
\caption{$\fround : \REAL \to \INTEGER$ \label{fig:Round2}}
\def\baselinestretch{1.15}\selectfont%
{\fontsize{10}{15}
\begin{algorithmic}[1]
\Prog{$x:\REAL$}
\State 
    $\LET k,\ j,\ b:\INTEGER=0$;\;
    $\LET y:\REAL=x$;
\Mywhile{$\mychoose(\myabs{y} \elt 1, 1/2 \elt \myabs{y}) \ieq 1$}
    \label{line:round:while1S}
    \RightComment  {$x=y\cdot2^j$}
    \State $j \dfeq  j+1;\; y \dfeq  y/2$
\EndMywhile;
    \label{line:round:while1E}
    \RightComment  {$|y|\elt 1$}
\Mywhile{$1 \leq j$}
    \label{line:round:while2S}
    \RightComment {$x = (y+k)\cdot2^j\land|y|<1$}
    \State  $y \dfeq  y\mult 2$;
    \State  $b \dfeq  \mychoose (y\elt 0, -1\elt y\elt 1, 0 \elt y)-1$;
        \label{line:round:myselect}
    \Statex \RightComment{$b=$most-significant \emph{signed} binary digit of $y$}
    \State \IF\; $b \ieq -1$ \;\THEN\; $y\dfeq  y+1$\; \eendif;  \IF\; $b \ieq 1$ \;\THEN\; $y\dfeq  y-1\; \eendif;$
    \State $k \dfeq  k+k+b$;\; $j \dfeq  j-1$
        \label{line:round:split}
\EndMywhile
    \label{line:round:while2E}
    \RightComment {$|x-k|<1$}
\EndProg{$k$}
\end{algorithmic}
}
\end{algorithm}

Due to the multi-valuedness of the loop condition at Line~\ref{line:round:while1S}, when the loop (Lines~\ref{line:round:while1S}~to~\ref{line:round:while1E}) exits, the second argument $1/2 \elt y$
may still be $\true$, whereas the first
$|y| \elt 1$ must be $\true$.
In the integer loop (Lines~\ref{line:round:while2S}~to~\ref{line:round:while2E}), multi-valuedness
strikes only at Line~\ref{line:round:myselect}
which employs $\mychoose()$ with trinary argument.

\subsubsection{Determinant Function via Gaussian Elimination} \label{ss:Gauss}

The determinant of a $d\times d$ matrix $A=(a_{ij})_{_{i,j}}$ is given by Leibniz' formula
\begin{equation}
\label{eqn:Det1}
\det(A) \;=\; \sum\nolimits_{\pi} \sign(\pi)\cdot\prod\nolimits_{j=1}^d a_{j,\pi(j)}
\end{equation}
where the sum ranges over all $d!$ permutations $\pi:\{1,\ldots,d\}\to\{1,\ldots,d\}$.
Since \ERC conveniently relieves the programmer from numerical issues like cancellation,
this formula gives rise to a straightforward \ERC program. 
However one executing arithmetic operations exponential in $d$.

Common numerical approaches therefore transform $A$ to triangular form, whose determinant is simply the product of its diagonal elements \cite[\S2.3.3]{Recipes}.
Following Turing \cite{Tur48}, apply Gaussian Elimination to determine a \emph{LU factorization with full pivoting} $P\cdot A\cdot Q=L\cdot U$ of $A$,
where $P, Q$ denote permutation matrices and $L$ and $U$ are lower and upper triangular matrices, respectively.
In Gauss' original algorithm, such search either (i) returns the index of a non-zero entry in the given sub-matrix
or (ii) asserts that said sub-matrix is identically zero.
By iterating this process, 
Gaussian Elimination determines the rank $k\in\IN$ of the original matrix $A\in\IR^{d\times d}$ --- which depends discontinuously on $A$'s entries and hence cannot be computed.

We thus change the specification of the determinant
(and of the LUPQ factorization\footnote{LUPQ factorization $A\mapsto (L,U,P,Q)$ is not unique, hence a real matrix-tuple-valued \emph{multi}-function.} it builds on):
with the promise for the argument matrix $A$ to have full rank---otherwise, its determinant will vanish, anyway.
\begin{equation}
\label{eqn:Det2}
\Det_d : \GL(\IR^d) \ni A \mapsto \begin{cases}
\det(A) &\text{if } \det(A)\neq 0,\\
\bot &\text{otherwise}.
\end{cases}
\end{equation}

We also relax pivot search (within the lower-right submatrix) to become a partial multi-valued function:

\begin{equation}
\label{eqn:Pivot}
\Pivot_{d}(A, k)  \;=\; \big\{ (i,j) \::\: A\in\IR^{d\times d} \land  k\leq i,j<d \land A_{i,j}\neq0 \big\}_{\{\bot\}}
\;\in\; \pdom{\IZ^2_\bot}
\end{equation}
Here, $X_{\{\bot\}} \dfeq X$ if $X \neq \emptyset$ and $\{\bot\}$ if $X = \emptyset$.
So the argument is a real $d\times d$ matrix $A$ and an integer $k$, indicating that a pivot is to be sought for in the $(d-k)\times(d-k)$ non-zero sub-matrix $A[k\ldots d-1,k \ldots d-1]\dfeq \big(A[i,j]\big)_{k\leq i,j<d}$.
\begin{algorithm}
\caption{$\pPivot_d : \REAL[d\times d]\times\dint \to \INTEGER$ \label{fig:pivot}}
\def\baselinestretch{1.15}\selectfont%
{\fontsize{10}{15}
\begin{algorithmic}[1]
\Prog{$A:\REAL [d\times d],\; k:\INTEGER$}
\State $\LET i_0,\ j_0:\INTEGER=k; \; \LET x:\REAL=0$;
\State \FOR\;$i:\INTEGER=k \;\TO\; d-1$ \DO
\State \qquad   \FOR\; $j:\INTEGER=k \;\TO\; d-1$ \DO
\State \qquad \qquad   $x \dfeq  \max\big(x,\abs(A[i,j])\big)$ \ENDFOR\; \ENDFOR;
\State \FOR\; $i:\INTEGER=k \;\TO\; d-1$ \DO
\State \qquad   \FOR\; $j:\INTEGER=k \;\TO\; d-1$ \DO
\State  \qquad\qquad \IF\; $\mychoose(\abs(A[i,j])\elt x, x/2 \elt \abs(A[i,j])) \ieq 1$ 
\THEN \; $i_0\dfeq  i; j_0\dfeq  j$\;\ENDIF 
\State \ENDFOR\;  \ENDFOR
\EndProg{$(i_0,j_0)$}
\end{algorithmic}
}
\end{algorithm}

In the program $\pPivot_d$ in \autoref{fig:pivot}, 
we represent matrices with fixed dimensions using \ERC's one-dimensional arrays.
For each positive natural numbers $n$ and $m$, we write $A[n\times m]$, to 
express a $n\times m$-dimensional matrix, simulated by $A[k]$ where $k = n\times m$.
Accessing an entry $A[i, j]_{n}$ is then a shorthand for $A[i+j+\cdots+j]$
denoting $A[i + n\cdot j]$. The repeated addition is only possible as $n$ is a fixed natural number constant, not a $\ERC$ term while $i, j$ are. The subscript $n$ is omitted for better readability when it is obvious.
Similarly the pair $(i_0,j_0)$ of return values is thus understood as an abbreviation of the single integer $i_0+d \cdot j_0$ where $d$ is the dimension of the input matrix.

The program introduces more abbreviations.
For any integer terms $k, l$ and a command $S$, 
$\FOR\; n:\INTEGER=k \;\TO\; l \;\DO\; S\; \ENDFOR$
refers to 
\[
\LET n:\INTEGER=k; \;
\LET m :\INTEGER = l;\;
\WHILE (n\leq m) \;\DO\;  S ; \; n\dfeq  n+1 \; \ENDWHILE
\]
with $n, m$ not appearing in $S$. 
For a Kleenean typed term $b$ and a command $S$,
\[\IF\;b\; \THEN\;S\; \ENDIF\] 
is used as an abbreviation for
\[
\eif\;b\;\ethen\;S\;\eelse\;\eskip\;\eendif
\]
Moreover, $\max(u, v)$ is a term abbreviation for 
$\mycond{u \elt v}{v}{u}$.

Based on $\pPivot_d$, the $\ERC(\emptyset, \{\pPivot_d\})$ program $\pDet$ in \autoref{fig:Det} computes the non-zero
determinant $\Det_d$ from Equation~(\ref{eqn:Det2}) via LUP decomposition with full pivoting.
Note that the precision parameter $p\in\IZ$ is present but ignored since the result gets computed exactly.

\begin{algorithm}
\caption{$\pDet : \REAL[d\times d] \to \REAL$ \label{fig:Det}}
\algrenewcommand{\algorithmiccomment}[1]{\hfill// #1}
{\fontsize{10}{15}
\begin{algorithmic}[1]
\Limit{p}{A:\dreal[d\times d]}  \RightComment{ $A$ invertible, $p$ ignored}
\normalmarginpar
\State $\LET i:\INTEGER=0$;\;
    $\LET j:\INTEGER=0$;\;
    $\LET k:\INTEGER=0$;\;
\State $\LET p_i:\INTEGER=0$;\;
    $\LET p_j:\INTEGER=0$;\;
    $\LET det:\REAL=1$;   \RightComment{ret.val}
\Myfor{$k\dfeq  0\; \TO\; d-2$}
    \State    \RightComment{Convert $A[k..d-1,k..d-1]$ to reduced row echelon form:}
    \State \label{line:gauss:pivot}
        $(p_i,p_j) \dfeq  \pPivot(A,k)$; \RightComment{$pi,pj\geq k$ s.t. {$A[p_i,p_j] \neq 0$.}}
    \State $\det \dfeq  \det \mult  A[p_i,p_j]$;
    \State \FOR\; $j \dfeq  0$\; \TO\; $d-1$\; \DO\;
                $\swap(A[k,j],A[p_i,j])$\;
            \ENDFOR;
        \Statex \RightComment{Exchange rows \#$k$ and \#$ pi$}
    \State \IF \; $k\neq p_i$ \THEN\; $\det\dfeq  -\det$;  \RightComment{flip sign}
    \State \FOR\; $i \dfeq  0$\; \TO\; $d-1$\; \DO\;
	    	$\swap(A[i,k],A[i,p_j])$\;
            \ENDFOR;
        \Statex \RightComment{Exchange columns \#$k$ and \#$ p_j$}
    \State \IF\; $k\neq p_j$ \THEN\; $\det\dfeq  -\det$;  \RightComment{flip sign}
    \Myfor{$j \dfeq  k+1\; \TO\; d-1$}
        \RightComment{Scale row \#$k$ by $1/A[k,k]$ and}
		\State $A[k,j]\dfeq  A[k,j]/A[k,k]$;
            \RightComment{and subtract the $A[i,k]$-fold from}
      	\State \FOR\; $i \dfeq  k+1$\; \TO\; $d-1$\; \RightComment{from rows \#$i=k+1\ldots d-1$}
      	\State \qquad \DO\;
                    $A[i,j]\dfeq  A[i,j]-A[i,k]\mult  A[k,j]$\;
                \ENDFOR
    \EndMyfor; $A[k,k] \dfeq  1$;\;    \FOR\; $i \dfeq  k+1$\; \TO\; $d-1$\;\DO\;$A[i,k] \dfeq  0$\;\ENDFOR
\EndMyfor;
\State $\det \dfeq  \det \mult  A[d-1,d-1]$
\EndLimit{\det}{p}
\end{algorithmic}
}
\end{algorithm}

Note that pivot search in Line~\ref{line:gauss:pivot} of program $\pDet$ 
is guaranteed to succeed
in that the $(d-k)\times(d-k)$ submatrix $A[k\ldots d-1,k\ldots d-1]$
under consideration will indeed contain
at least one --- usually non-unique --- non-zero element
due to the promise/restriction that $A\in\GL(\IR^d)$.

\subsubsection{Simple Unique Root Finding Functional} \label{ss:Root}

The problem of finding (i.e. approximating) a root to a given real function $f$
occurs ubiquitously in numerics under various hypotheses.
We consider an algorithmic version of the Intermediate Value Theorem which is to find the unique root to a given continuous $f:[a;b]\to\IR$ satisfying $f(a)<0<f(b)$. 
This case is commonly treated using Bisection:
Determine the sign of $f(x)$ at the interval midpoint $x\dfeq (a+b)/2$
and recurse to either $[a;x]$ or to $[x;b]$ accordingly.
However, since equality is undecidable, the sign test fails in case $f(x)=0$.
\emph{Tri}section \cite[p.~336]{Her96b} instead
considers the signs of $f$ at both one third $x'\dfeq (2a+b)/3$
and at two third $x''\dfeq (a+2b)/3$ of the interval, in parallel;
and recurses to either $[a;x'']$ or to $[x';b]$ accordingly:
Now at most one of the two sign tests at $f(x')$ and $f(x'')$ can fail,
provided that $f$'s root is unique and $f(a)\cdot f(b) < 0$.
These hypotheses also avoid common counterexamples like \cite{Spe59} or \cite[Theorem~6.3.2]{Wei00}.
To summarize, we consider the (single-valued) root finding problem $\Root:f\mapsto x$
with $f(x)=0$ for $f$ satisfying the following first-order predicates:

\begin{eqnarray*}
\cont(f,a,b) \quad&:\equiv&\quad \\
\all{\epsilon>0} &\some{\delta>0} & \all{ x,x'}
\; a\leq x\leq x'\leq x+\delta\leq b \;\Rightarrow\; |f(x)-f(x')|\leq\epsilon \\
\uniq(f,a,b) \quad&:\equiv&\quad
\cont(f,a,b) \land
f(a)\cdot f(b) < 0 \land
\usome{x} a <x<b \land f(x) = 0
\end{eqnarray*}
The first condition  $\cont(f, a, b)$ says $f$ is continuous in the usual sense on the interval $[a; b]$
and the second condition  $\uniq(f, a, b)$ says $f$ is continuous in the interval $[a;b]$, admits a unique root in the interval $(a; b)$, and
the signs of $f(a)$ and $f(b)$ are different.


The $\ERC$ program $\Trisection$ in \autoref{fig:Trisection} is annotated with
precondition $\uniq(f,a,b)$
and postcondition
$\uniq(f,a,b)\land |b-a|\leq2^{p}$, as well as loop invariant
$\uniq(f,a,b)$.
Formally, as $\ERC$ does not accept function arguments, the program is defined in $\ERC(\{f\}, \emptyset)$.
Given $p$, the postcondition guarantees that $\Trisection$ indeed
returns an approximation to the root up to error $2^p$ when $\uniq(f,a,b)$ holds.
Therefore, as $p \to -\infty$ the program returns the root of $f$.
This program is chosen later in Section~\ref{ss:Trisection} to be formally proved for its correctness.

\begin{algorithm}
\caption{$\Trisection : \dreal \times \dreal \to \REAL$ in $\ERC(\{f\}, \emptyset)$ for any $f : \IR \to \IR_\bot$ \label{fig:Trisection}}
\def\baselinestretch{1.15}\selectfont%
{\fontsize{10}{15}
\begin{algorithmic}[1]
\Limit{p}{a : \dreal, b : \dreal}  \RightComment {$\uniq(f,a,b)$}
\Mywhile{$\mychoose\big(b-a\elt\myprec{p}  ,\,  \myprec{p-1} \elt b-a\big)  \ieq  1$}
    \RightComment{$\uniq(f,a,b) \land  {b-a>2^{p-1}}$}
	\State \IF \; $\mychoose(f(b/3 + 2\mult a/3)\mult f(b) \elt 0,\, f(a) \mult  f(2\mult b/3+a/3) \elt 0) \;\ieq\; 1$
	\State \THEN\;  
        $b\dfeq  2\mult b/3+a/3$
    \State \ELSE \;\;
         $a\dfeq  b/3+2\mult a/3$ \;
    \ENDIF
\EndMywhile
    \RightComment{${\uniq(f,a,b) \land |b-a|\leq2^{p} }$}
\EndLimit{a}{p}
\end{algorithmic}
}
\end{algorithm}

\subsection{Turing-Completeness} \label{ss:Complete}

Algorithms~\ref{fig:Round1}~and~\ref{fig:Round2} show the existence of a well-formed command which we abbreviate as
\[
\Gamma \vdash z :\equiv \cround(x) 
\triangleright \Gamma 
\]
when $\Gamma(x) = \dreal$ and $\Gamma(z) =  \dint$ 
that assigns the multi-valued rounding of $x$ into $z$ having its denotation
\[
\sem{\Gamma \vdash z :\equiv \cround(x) \triangleright \Gamma }\sigma = \{
\sigma[z\mapsto k] \mid \sigma(x) - 1 < k < \sigma(x) - 1\}.
\] 
Furthermore, though it is not available in the term language, the command language of \ERC provides the \naive coercion from integers to reals:
There exists a well-formed command which we write
\[
\Gamma \vdash x :\equiv \creal(z) \triangleright \Gamma 
\] 
when $\Gamma(x) = \dreal$ and $\Gamma(z) =  \dint$ 
that repeatedly add $1$ or $-1$ in a loop such that
\[
\sem{\Gamma \vdash x :\equiv \creal(z) \triangleright \Gamma }
\sigma = \{\sigma[x \mapsto \sigma(z)]\}.
\]

\ERC being an imperative language providing integers, arithmetical operations on them, compositions of commands, and general loops, it is type-1 Turing-complete over natural numbers:
\begin{prop}
\label{p:ordinarycomplete}
For any type-1 computable (or equivalently $\mu$-recursive)  
$f : \IN^d \pto \IN$, there exists a well-formed $\ERC$ command which we abbreviate as
\[
\Gamma \vdash 
y :\equiv f(x_1, \cdots, x_d) 
\triangleright \Gamma
\]
when $\Gamma(y) = \dint, \Gamma(x_i) = \dint$ for each $i$,
such that for all $\sigma\in\sem{\Gamma}$ where $(\sigma(x_1), \cdots, \sigma(x_d)) \in \dom(f)$,
\[
\sem{\Gamma \vdash 
y :\equiv f(x_1, \cdots, x_d) \triangleright \Gamma} \sigma = \{\sigma[y \mapsto f(\sigma(x_1), \cdots, \sigma(x_d))]\}.
\]
\end{prop}
\begin{proof}
We can construct the command $S$ inductively on the $\mu$-recursiveness of $f$. 
\end{proof}

From the fact that \ERC is type-1 Turing-complete, we can automatically ensure 
many interesting operations on integers where some of which are necessary to build our (type-2) Turing-completeness result.
Fix a standard numbering of finite sequences (including the empty sequence $\epsilon$) of integers $\zeta : \IN \to \IZ^{*}$ throughout this section. 
By the type-1 completeness, there are \ERC commands for (1) constructing 
the code of the empty string $\epsilon$, (2) appending a new tail to a sequence from its code, (3) accessing the tail of a sequence from its code assuming that the sequence is not empty, and (4) obtaining the length of a sequence from its code.
For any integer variables $y, x, x_1, x_2$ in $\Gamma$:
\begin{align*}
\sem{
\Gamma \vdash y :\equiv \epsilon\triangleright\Gamma} \sigma
&= \{\sigma[y\mapsto\zeta^{-1}(\epsilon)]\}
\\
\sem{\Gamma \vdash y :\equiv x_1 \myappend x_2 \triangleright\Gamma }\sigma
&= \{\sigma[y\mapsto\zeta^{-1}(\zeta(\sigma(x_1)) \myappend \sigma(x_2))]\}\\
\sem{
\Gamma \vdash y :\equiv \text{tail}(x)\triangleright\Gamma }\sigma
 = \{\sigma[y\mapsto m]\}
\text{ if } \zeta(\sigma(x)) &= a\myappend m \text{ for some }a \in \IZ^*\\
\sem{
\Gamma \vdash y :\equiv \text{length}(x)\triangleright\Gamma }
\sigma &= \{\sigma[y \mapsto m]\}
\text{ if } \zeta(\sigma(x)) \text{ has length }m
\end{align*}
Note that the command for accessing the tail of a sequence
$y :\equiv \text{tail}(x)$
may not well-behave when the sequence that $x$ represents is the empty sequence. 
Here, $x\myappend y$ when $x \in \IZ^*$ and $y\in \IZ$
denotes the finite sequence of integers where $y$ is appended to $x$.

For the last building block of our main theorem, we recall the type-1 characterization of type-2 machines:
\begin{prop}
\label{p:slice}
For any type-2 machine computing 
$\tau : \IZ^\IN\times \cdots \times \IZ^\IN \pto \IZ^\IN$, 
there exists a total type-1 computable 
(with regard to the standard numbering of finite sequences) 
$\bar\tau : \IZ^* \times \cdots \times \IZ^* \to \IZ^*$ such that for any $(\phi_{1}, \cdots, \phi_{d}) \in \dom(\tau)$, 
\[(\bar\tau(\bar\phi_1^{(m)}, \cdots, \bar\phi_d^{(m)}))_{m\in\IN}\]  (which 
can contain the empty sequence) is a chain with regard to the prefix ordering and it converges to $\tau(\phi_{1}, \cdots, \phi_{d})$.
Here, $\bar\phi^{(m)} \in \IZ^*$ denotes the length $m$ prefix of $\phi \in \IZ^\IN$.

\end{prop}
The type-1 $\bar\tau$ continuously approximates $\tau$ in the following sense. If $\phi \in \dom(\tau)$,
$\bar\tau$ tries to approximate $\tau(\phi)$ using finite approximations of $\phi$. 
The function $\bar\tau$ may fail initially due to for example lack of input precision and return the empty sequence. However, it eventually returns longer and longer prefixes of $\tau(\phi)$ converging to $\tau(\phi)$ consistently (i.e. monotonically with regard to the prefix ordering).
\begin{proof}
Let us first consider the case where $\tau$ is univariate.
As stated in \cite[Excercise~2.3.12]{Wei00}, 
due to the relation between domain-computability and type-2 computability,
for any partial $\tau$ computed by a type-2 machine, there exists a total function $T : \IZ^*\to \IZ^*$ which is monotone 
with regard to the prefix ordering, approximates $\tau$ in the sense that 
\[\lim_{n \to \infty} T(\bar\phi^{(n)}) = \tau(\phi)\]
for each $\phi \in \dom(\tau)$, and
the set
\[
\{
(u, v) \in \IZ^* \times \IZ^* \mid 
v \text{ is a prefix of }T(u)
\}
\]
is computably enumerable. 
Let $e : \IN \to \IZ^* \times \IZ^*$ be one computable 
enumeration. 
Then, we can define $\bar\tau$ as follows:
given $\bar\phi^{(n)}$, iterate $e(0), \cdots, e(n)$
to find $(u, v) \dfeq e(i)$ where $u$ is a prefix of $\bar\phi^{(n)}$
and the length of $v$ is the longest amongst such.
If such $(u, v)$ is located, return $v$. 
If there is no prefix of $\bar\phi^{(n)}$ found, return the empty sequence $\epsilon$.
This computational procedure describes the computable function $\bar\tau$. 

Now consider a computable $d$-variate $\tau$ and construct $\bar\tau$.
As $\tau$ is computable, there exists a computable $\tau' : \IZ^\IN \pto \IZ^\IN$ 
such that $\tau(\phi_1, \cdots, \phi_d) = \tau'\langle \phi_1, \cdots, \phi_d\rangle$ where $\phi \dfeq \langle \phi_1, \cdots, \phi_d\rangle$ is defined by interleaving. By the above proof, 
there exists $\bar\tau' : \IZ^* \to \IZ^*$ for $\tau'$ satisfying the conditions.
Given the prefixes $\bar\phi_1^{(n)}, \cdots, \bar\phi_d^{(n)}$, 
we construct $\bar{\phi}^{(n\cdot d)}$ by interleaving, feed it in $\tau'$, and let it be the return sequence of $\bar\tau(\bar\phi_1^{(n)}, \cdots, \bar\phi_d^{(n)})$. 
\end{proof}

\begin{thm}
\label{t:Complete}
Recall the definitions of 
computable partial (mutli-valued) functions in Section~\ref{s:preliminaries} and 
\ERC programs realizing them in Definition~\ref{d:Realizer}.
\begin{enumerate}
\item Every computable partial real function from integers and reals can be realized in $\ERC_0$.
\item Similarly, every computable partial integer multi-valued function from integers and reals can be realized in $\ERC_0$.
\end{enumerate}
\end{thm}
As stated in Proposition~\ref{p:ordinarycomplete}, $\ERC_0 $ can realize all $\mu$-recursive integer functions,
and consequently also every computable single real number $r\in\IR$,
since the latter is defined as limit $r=\lim_n \phi(n)/2^n$ 
for $|r-\phi(n)/2^n|\leq2^{-n}$.
Note that also every oracle Turing machine computable integer function
is oracle-recursive; see
\url{http://math.stackexchange.com/questions/2778974}.
This suggests realizing
a type-2 oracle machine computing a real function $f:\IR\to\IR_\bot$ 
according to \cite[\S2.3]{Ko91} in $\ERC_0$:
by replacing any oracle query ``$\phi(n)$'' 
for a $2^{-n}$--approximation $\phi(n)/2^n$ to the real argument $x\in\dom(f)$ 
by the $\ERC_0$ commands in Section~\ref{ss:Round} realizing $\Round(x\mult 2^n)$.
However said $\Round()$ is multi-valued, corresponding to a multi-valued oracle $\phi$;
and for these, oracle Turing computation does not seem to be known as equivalent
to oracle $\mu$-recursiveness.
Instead, we present in the sequel a different proof of Theorem~\ref{t:Complete}:

\begin{proof}
(1) Without loss of generality, suppose any computable partial $f : \IR^{d_1} \times \IZ^{d_2}  \to \IR_\bot$ where $d_1+ d_2 = d$. 
By definition,
there is a type-2 machine computing 
$\tau : \IZ^{\IN} \times \cdots \times \IZ^{\IN} \pto \IZ^{\IN}$
such that for each $(x_{1}, \cdots, x_d) \in \dom(f)$ and each $\phi_{i}$ a name of $x_i$, $\tau(\phi_1, \cdots, \phi_d)$ is defined and is a name of $f(x_1, \cdots x_d)$:
\[
\myabs{f(x_{1}, \cdots, x_{d}) - \tau(\phi_{1}, \cdots, \phi_{d})(n)\cdot 2^{-n}} \leq 2^{-n}\quad\text{holds for all }n.
\]
By Propositions~\ref{p:ordinarycomplete}~and~\ref{p:slice}, 
there exists a well-formed command in \ERC such that
when $x_i', y'$ are integer variables in $\Gamma$, for any $\sigma \in \sem{\Gamma}$
such that $\sigma(x_i')= m_i^{(n)} \in \IN$ a code for $\bar\phi_i^{(n)}$ (w.r.t. the fixed standard numbering $\xi$),
\[
\sem{
\Gamma \vdash y' :\equiv \bar\tau(x_1', \cdots, x_d')}\sigma 
 = \{\sigma[y' \mapsto y^{(n)}]\}
\]
holds for some $y^{(n)} \in \IN$ which is a code for a finite, possibly empty, prefix of $\tau(\phi_{1}, \cdots, \phi_{d})$.

\begin{algorithm}
\caption{A real program in $\ERC_0$ realizing 
$f$ \label{fig:completemachine} in the proof of Theorem~\ref{t:Complete}}
\def\baselinestretch{1.15}\selectfont%
{\fontsize{10}{15}
\begin{algorithmic}[0] 
\Limit{p}{x_1:\dreal, \cdots, x_{d_1} : \dreal, x_{d_1+1} : \dint, \cdots, x_d : \dint}
\State $\LET x_1',\ \cdots,\ x_d',\ y',\ y'' \abbeq \epsilon;$
\State $\LET \ell,\ t,\ n : \dint = 0;$
\Mywhile{$\ell  \leq -p\; \klor\; \ell \leq 0 $}
    \LeftComment{$x_i'$ is a code for the length $n$ prefix of some name $\phi_i$ of $x_i$}
    \LeftComment{$y'$ is a code for the length $\ell$ prefix of $\tau(\phi_{1}, \cdots, \phi_{d})$ for some name $\phi_i$ of $x_i$}
    \LeftComment{append $x_i'$:}
    \State $t \abbeq \text{round}(x_{1} \times \myprec{n});$
    \State $x_{1}' \abbeq x_1' \myappend t;$
	\State $\vdots$
    \State $t \abbeq \text{round}(x_{d_1} \times \myprec{n});$
    \State $x_{d_1}' \abbeq x_{d_1}' \myappend t;$
    \State $x_{{d_1}+1}' \abbeq x_{{d_1}+1}'\myappend x_{{d_1}+1};\; \cdots ; \; x_{d}' \abbeq x_{d}'\myappend x_{d};$
    \State $n \dfeq n + 1;$
    \LeftComment{$x_i'$ is a code for the length $n$ prefix of some name $\phi_i$ of $x_i$}
    \State $y' \abbeq \bar\tau(x_{1}', \cdots, x_{d}');$ 
	\State $\ell \abbeq \text{length}(y');$ 
	\LeftComment{$y'$ is a code for the length $\ell$ prefix of $\tau(\phi_{1}, \cdots, \phi_{d})$ for some names $\phi_i$ of $x_i$}
 \EndMywhile;
\LeftComment{$y'$ is a code for a length $\ell$ prefix of $\tau(\phi_{1}, \cdots, \phi_{d})$ for some names $\phi_i$ of $x_i$ and}
\LeftComment{$\ell > 0$ and $\ell -1 \geq -p$} 
\State $y : \dreal = 0;$
\State $y'' \abbeq \text{tail}(y');$
\RightComment{$\ell > 0$ ensures $\text{tail}(y')$ does not fail}
\State $y \abbeq \text{real}(y'');$
\RightComment{$\myabs{f(x_1, \cdots, x_d) -y\times2^{-(\ell-1)}} \leq 2^{-(\ell - 1)} \leq 2^p$ }
\EndLimit{y \times 2^{-(\ell-1)}}{p}
\end{algorithmic}
}
\end{algorithm}

We claim that the real program in $\ERC_0$ in \autoref{fig:completemachine} realizes $f$.
On its input variables $x_1, \cdots, x_d$, the program prepares integer variables $x_1', \cdots, x_d'$ to let 
each $x_i'$ stores a code of a finite prefix (of say length $n$) of a name of $x_i$. Starting from the empty sequences ($n=0$), the program keeps appending $x_i'$ and computes $\bar\tau(x_1', \cdots, x_d')$ in $y'$ which is a finite prefix (of say length $\ell$) of some name of $f(x_1, \cdots, x_d)$. 
The program iterates until the length $\ell$ of the obtained prefix $y'$ becomes enough to obtain a $2^{-p}$ approximation to $f(x_1, \cdots, x_d)$.

The variable $n$ increases as the number of iterations.
Each state that appears at the beginning of the $n$'th iteration 
satisfies that 
\[
x_i' \text{ is (a code for) the length $n$ prefix of some name $\phi_i$ of $x_i$}\quad\text{and}
\]
\[
\text{$y'$ is (a code for) the length $\ell$ prefix of $\tau(\phi_{1}, \cdots, \phi_{d})$ for some name $\phi_i$ of $x_i$}
\]
irrelevant to the multi-valuedness of the integer rounding.  
(Note that though due to multi-valuedness, $\phi_i$ in each annotated assertion does not need to be consistent.)

Suppose for now that the loop is total.
When the loop exits, due to $\ell -1 \geq -p$ and $\ell > 0$ the tail of $y'$, which is $y$, is the $\ell-1$'th (counting from $0$) entry of 
$\tau(\phi_{1}, \cdots, \phi_{d})$ for some names $\phi_i$ of $x_i$.
Hence, it holds that 
$\myabs{f(x_1, \cdots, x_d) -y\times2^{-(\ell-1)}} \leq 2^{-(\ell - 1)} \leq 2^p$ where $y\times2^{-(\ell-1)}$ is the return value. We conclude that the real program indeed realizes $f$.

To prove that the loop is total,
we show that for any natural number $k$, there is a natural number $n(k)$ such that for any names $\phi_i$ of $x_i$, 
$\bar\tau(\phi_1^{(n(k))}, \cdots ,\phi_d^{(n(k))})$ has length greater than or equal to $k$.
If this is true, the loop terminates after at most $n(-p)$ (or $n(1)$ if $p$ is non-negative) iterations.

Recall that  $A := \{(\phi_1, \cdots, \phi_d) \mid \phi_i\text{ is a name of }x_i\}$ is compact with regard to the standard topology of Baire space.
For each natural number $n$, define $A_n := \{(\phi_1, \cdots, \phi_d) \mid \phi_i\text{ is a name of }x_i \land
\bar\tau(\phi_1^{(n)}, \cdots ,\phi_d^{(n)})\text{ has length greater than or equal to $k$}
\}$ which is the set of names whose length $n$ prefix is enough to make $\bar\tau$ returns length $k$ prefix of $\tau$. Note that $A_n$ is open and by the convergence of $\bar\tau$, $(A_n)_{n\in\IN}$ forms a countable open cover of $A$. The compactness gives us a finite subcover of $A_n$ and we can choose the maximum index of the finite subcover to be $n(k)$.

\noindent (2) The integer multi-valued function case can be proved similarly except that now we only need to iterate until $\ell > 0$
and return $\text{tail}(y')$.
\end{proof}

\begin{rem}
\label{r:cc}
The completeness of \textbf{WhileCC*}'s algebraic semantics \cite{TZ04} is obtained in a similar way: 
using nondeterminism to obtain discrete approximation then proceed with ordinary computation.
However, their countable choice construct is not necessary for our case as we only consider real numbers;
using the total order of reals, we can rely on our finite choice construct to extract names.
\end{rem}

$\ERC$ defines integer multi-valued functions but not a real multi-valued function. 
Real programs are forced to take limit operations and compute single-valued real numbers. 
This is our design choice that can be chosen differently. 
Here, we state a possible alternative and explain our choice:
\begin{rem}
\label{r:realmult}
One straightforward way to allow real multi-valued functions is to ease our typing rule such that we let 
the program
\[
\begin{array}{rl}
P \dfeq & \einput\;x_1 : \tau_1, x_2 : \tau_2, \cdots, x_d : \tau_d \\
& \quad S \\
& \ereturn\;t
\end{array}
\]
have type $\tau_1 \times \cdots \tau_d \to \dreal$ if there is a context $\Gamma'$ such that
\[
   x_1:\tau_1, x_2 : \tau_2, \cdots, x_d : \tau_d \vdash S \triangleright \Gamma'
   \quad\text{and}\quad \Gamma' \vdash t : \dreal \,.
\]
We let its denotation be $\sem{P} \dfeq \sem{t}^\dagger \fcomp \sem{S} : \sem{\tau_1}\times\cdots\times\sem{\tau_d} \to \pdom{\IR_\bot}$ as it was an integer program.
From Lemmas~\ref{l:Computable_Term}~and~\ref{l:Computable_Command}, it can be seen immediately that the denotations of such real programs are computable multi-valued functions.
We can also extend the notion of realization in Definition~\ref{d:Realizer} accordingly.

However, this extension does not carry the completeness in Theorem~\ref{t:Complete} over.
Recall that the denotations of programs are finite multi-valued functions. 
Hence, a computable infinite multi-valued function is realizable in $\ERC$ only if it admits point-wise finite refinement that is computable.
In the integer case, it was feasible as any infinite integer multi-valued function admits point-wise finite refinement. 
(Recall Fact~\ref{fact:prelim 2}.)
However, this does not work for real multi-valued functions.
That means when we extend $\ERC$ with real multi-valued functions as above, we end up with a restricted completeness theorem. 
As easing the restriction as above is straightforward, we choose to keep the restriction and achieve a more natural completeness theorem.
\end{rem}

\section{Logic of \emph{Exact Real Computation}}
\label{s:Logic}
In this section, we devise a specification language and based on the language, devise an extended Hoare logic for formal verification of $\ERC$ programs.
Section~\ref{ss:Theory} proposes a three-sorted structure which we use to rigorously specify (multi-valued) functions with real arguments in Section~\ref{ss:Specification}.
It is carefully designed to be rich enough to allow arguing about computations in \ERC yet restricted such as to assert logical decidability as a guarantee to formal program verification; recall Remark~\ref{r:Definability}.
Section~\ref{ss:Hoare} extends the classical Hoare logic which is sound with regard to our denotational semantics.

\subsection{Three-Sorted Structure with Decidable Theory}
\label{ss:Theory}

\begin{defi}
\label{d:ERC_theory}
The \emph{Structure of} $\ERC_0$ is the three-sorted structure $\Structure_0$
combining the Kleene Algebra $(\Kleene,\cff,\ctt,\cuu,\kland,\klor,\kneg)$
with Presburger Arithmetic $(\IZ,0,1,+,-,\leq,2\IZ,3\IZ,4\IZ,\ldots)$
and ordered field $(\IR,0,1,+,-,\mult,<)$.
Available inter-sort functions are 
the \emph{binary precision} embedding $\imath:\IZ\ni p\mapsto 2^p\in\IR$
and its partial half-inverse $\lfloor\log_2\circ\abs\rfloor:\IR\setminus\{0\}\to\IZ$.
Here $k\IZ$ denotes the predicate on $\IZ$ which is $\ctt{}$ precisely for all integer multiples of $k\in\IN$.

The \emph{(specification) language of} $\ERC_0$ is the first-order language $\Logic_0$ over the signature of the structure $\Structure_0$;
and the \emph{theory of} $\ERC_0$ is the first-order theory $\Theory_0$ containing true statements of $\Logic_0$ in $\Structure_0$.

\end{defi}
We take the inequality on $\IZ$ as non-strict $\leq$, but that on $\IR$ strict $<$ as in \cite{Marker}.
This classical predicate $<\subseteq \IR\times\IR$ is the classical order relation of 
real numbers and is not to be confused with anything computational. 
Another important note is that the symbol $\bot$ in our denotational semantics does not appear in the structure of \ERC.

Theorem~\ref{t:Eike}(a) shows that the theory $\Theory_0$ is decidable.
This applies for example to preconditions and postconditions or loop invariants of $\ERC_0$ programs;
see Lemma~\ref{l:pred_trans2} and Remark~\ref{r:SideCondition}.
This is a significant advantage of \ERC, compared to traditional programming languages for discrete data:
Classical While programs over integers with multiplication, for instance, do suffer from G\"odel undecidability {\rm\cite[\S6]{Cook78}}.

\begin{thm} 
\label{t:Eike}
\hfill
\begin{enumerate}[label=(\alph*)]
\item
The Theory $\Theory_0$ of $\ERC_0$ is decidable.
\item
$\Theory_0$ is also `model complete' in that it admits elimination of quantifiers up to one (by choice either existential or universal) block ranging over integers.
\item
Each of the following expansions destroys the decidability of the first-order theory of the structure $\Structure_0$:
\begin{itemize}
\item expanding with integer multiplication 
\item expanding with the unary predicate $\IZ$ on $\IR$, or with the standard
coercion $\IZ\hookrightarrow \IR$
\item replacing the binary precision embedding $\imath$ with its unary counterpart $\jmath:\IN_+\ni n\mapsto 1/n\in\IR$
\item expanding simultaneously with the real exponential and sine function and with transcendental constants $\pi$ and $\ln2$
\end{itemize}
\end{enumerate}
\end{thm}
\begin{proof}
(c) Including integer multiplication recovers Peano arithmetic and G\"odel undecidability via Robinson's Theorem \cite{Robinson}.
A unary predicate $\IZ$ on $\IR$ allows to express integer multiplication via the reals;
similarly, for (any total extension of) the unary precision embedding $\jmath$.
Finally, the real transcendental functions and constants make the theory undecidable according to Richardson's Theorem \cite{Richardson}.

\medskip\noindent
(a)+(b) A celebrated result of
van den Dries \cite{Dries1986}
extends classical Tarski-Seidenberg quantifier elimination 
from the first-order theory of real-closed fields
to the expanded structure
\begin{equation}
\label{eqn:Expanded}
(\IR,0,1,+,-,\mult,<,2^{k\IZ}:k\in\IN,2^{\lfloor\log_2\circ\abs\rfloor})
\end{equation}
with axiomatized additional predicates $2^{k\IZ},\, k \in \IN$,
and truncation function to binary powers $2^{\lfloor\log_2\circ\abs\rfloor}$,
see also \cite{DBLP:journals/tcs/AvigadY07}.

Note that both the real-closed field $(\IR,0,1,+,-,\mult,<)$
and Presburger Arithmetic can be embedded into the expanded structure from Equation~(\ref{eqn:Expanded});
the latter interpreted as its multiplicative variant
$(2^{\IZ},1,2,\mult,<,2^{k\IZ}:k\in\IN)$
is called \emph{Skolem Arithmetic} \cite{AlexisBes}:
\begin{itemize}
    \item Replace quantifiers over Skolem integers 
    with real quantifiers subject to the predicate $2^{k\IZ}$ for $k\dfeq  1$;
    \item Consider $\imath:\IZ\to\IR$ as the restricted identity $\id_{2^\IZ}$ in $\IR$.
\end{itemize}
Then every formula $\varphi$ with or without parameters in our two-sorted structure translates
signature by signature to an equivalent one $\tilde\varphi$ over the expanded theory
where quantifiers can be eliminated, yielding equivalent decidable $\tilde\psi$
(which may involve binary truncation $2^{\lfloor\log_2\circ\abs\rfloor}$).

To translate this back to some equivalent $\psi$ over the two-sorted structure,
while re-introducing only one type of quantifier, observe that
for real $x$:
\begin{eqnarray*}
 x\in2^{k\IZ} &\Leftrightarrow&
\some{z\in\IZ} z\in k\IZ\;\wedge\;x=\imath(z); \\
x\not\in2^{k\IZ} &\Leftrightarrow&
\some{z\in\IZ} z\in k\IZ\;\wedge\;\imath(z)<x<\imath(z+k)\,.
\end{eqnarray*}
Similarly, replace real binary truncation $2^{\lfloor\log_2\circ\abs(x)\rfloor}$ with
``$\imath(z)$'' for some/every $z\in\IZ$ s.t. $\imath(z)\leq |x|<\imath(z)+1$
in case $x>0$, with 0 otherwise.

Since the Kleene Algebra $\Kleene$ as the third sort is finite, it does not affect decidability.
\end{proof}

\subsection{Specification Language}\label{ss:Specification}
We define a specification language for each $\ERC(\calF, \calG)$ by expanding 
the structure of $\ERC_0$ as follows:
\begin{defi}
\label{d:spec}
The \emph{structure of $\ERC(\calF,\calG)$} is the expansion $\Structure(\calF,\calG)$ of 
$\Structure_0$ with the graphs of $f\in\calF$ and $g\in\calG$ as new relations:
\begin{itemize}
    \item $\text{Graph}_{f} \dfeq \{(x_1, \cdots, x_d; y) \mid f(x_1, \cdots, x_d) = y \neq \bot\}$ and
    \item $\text{Graph}_{g} \dfeq \{(x_1, \cdots, x_d; y) \mid y \in g(x_1, \cdots, x_d) \land \bot \not \in g(x_1, \cdots, x_d)\}$.
\end{itemize}
The \emph{(specification) language of $\ERC(\calF,\calG)$}, which we write $\Logic(\calF, \calG)$ is the first-order language $\Logic_{0}$ extended with the new relation symbols.
The \emph{theory of $\ERC(\calF,\calG)$}, which we write $\Theory(\calF,\calG)$, is the sentences in the language that is true in  $\Structure(\calF,\calG)$.
\end{defi}

Let us write simply $\Structure$, $\Theory$, and $\Logic$ 
when the underlying $\calF$ and $\calG$ are obvious or not so relevant.
Though, formally, we only add graphs of the (multivalued) functions, we may still have their function applications as appropriate abbreviations.

The following lemma shows that our specification logic is adequate for the term language of \ERC:
\begin{lem} 
\label{l:pred_trans2}
For each well-typed term $\Gamma \vdash t : \tau$ in $\ERC$, 
its denotation 
$\sem{\Gamma \vdash t : \tau} : \sem{\Gamma} \to \pdom{\sem{\tau}_\bot}$ is definable in $\Logic$
in the sense that
there is a formula $\psi(x_1, \cdots, x_d; y)$ defining the graph of the denotation (as in Definition~\ref{d:spec}).
We represent a fixed length array, of say length $d$, with 
$d$ variables $x = (x_{(0)}, \cdots, x_{(d-1)})$ in an obvious way.
\end{lem}
\begin{proof}

Let $X, Y_i, Z \in \{\sem{\tau} \mid \tau \text{ is a data type}\}$.
Suppose $f_i : X \to \pdom{(Y_i)_\bot}$ and $g : Y_1 \times \cdots \times Y_d \to \pdom{Z_\bot}$ are definable partial multi-valued functions where
 $\psi_i(x; y_i)$ defines $f_i$ and $\psi(y_1, \cdots, y_d; z)$ defines $g$.
Then, its composition $g^\dagger \circ \langle f_1, \cdots, f_d\rangle : X \to \pdom{Z_\bot}$ is definable by
\[
\psi(x; z) :\equiv
\some{y_1,\cdots,y_d}
\psi_1(x; y_1) \land \cdots \land \psi_d(x; y_d) \land
\psi(y_1, \cdots, y_d; z)\;.
\]
Furthermore, if a partial function $f : X \to Z_\bot$ is definable by $\psi(x; z)$, its embedding $f^\ddagger$ is definable by the same $\psi$.

Hence, we only need to check the nontrivial cases, the definability of primitive functions that are absent in the structure $\ERC$:
(1) $\semchoose$, (2) $\semkond$, and (3) $\semproj$: 
\begin{enumerate}
    \item See that $\semchoose_n : \IK^n \to \IZ$ is definable by
\[
\psi(x_1, \cdots, x_n; y) :\equiv (y = 0 \land x_1 = \ctt) \lor \cdots \lor (y = n-1 \land x_{n-1} = \ctt)\;.
\]
    \item
    Suppose $\psi_b(x; y)$ defines $b : X \to \pdom{\IK_\bot}$,
    $\psi_f(x; y)$ defines
    $f: X \to \pdom{\IR_\bot}$,
    and $\psi_g(x; y)$ defines $g :X \to \pdom{\IR_\bot}$.
    Then, the composition $\semkond^{\dagger_1} \fcomp \langle b, f, g\rangle$ is definable by
\begin{align*}
\psi(x; z) :\equiv&
\big(\psi_b(x; \ctt) \Rightarrow \some{y}\psi_f(x; y) \big) \\
\land& \big(\psi_b(x; \cff) \Rightarrow \some{y}\psi_g(x; y)\big)\\
\land& \big(\psi_b(x; \cuu) \\
&\quad \Rightarrow \big(\some{y} \psi_f(x; y) \land \psi_g(x; y) \land \all{y_1, y_2} \psi_f(x; y_1) \land\psi_g(x; y_2) \Rightarrow y_1 = y_2 \big)\big)\\
 \land \;\;\;&  \!\!\!\!
\big((\psi_b(x; \ctt) \land \psi_f(x; z)) \lor 
(\psi_b(x; \cff) \land \psi_g(x; z)) \lor
(\psi_b(x; \cuu) \land \psi_f(x; z)) \big)\;.
\end{align*}

\item 
The partial projection map $\semproj : \IR^d \times \IZ \to \IR_\bot$ is defined by  
\[
\psi(x, k; y) :\equiv 0 \leq k<d \land (k=0 \Rightarrow y =x_{(0)}) \land \cdots \land  (k=d-1 \Rightarrow y = x_{(d-1)})\;.
\]
Note that when $k$ is an index out of range, 
there is no $y$ satisfying $\psi(x, k; y)$
since $\semproj(x, k, y) = \bot$.
\qedhere
\end{enumerate}
\end{proof}

\begin{defi}
\label{d:trans}
For a well-typed term $\Gamma \vdash t : \tau$ in $\ERC$, let us write 
$\trans{\Gamma \vdash t : \tau}$ for the formula defining $\sem{\Gamma \vdash t : \tau}$ according to Lemma~\ref{l:pred_trans2}.
\end{defi}
\noindent We write simply $\trans{t}$ for $\trans{\Gamma \vdash t : \tau}$ when 
the particular context $\Gamma$ is obvious or irrelevant.
We let the free variables representing the input values in $\trans{t}$ be synchronized with $\Gamma$ and omit them in the notation; i.e., 
$\trans{t}(y) \dfeq \trans{x_1 : \tau_1, \cdots, x_d : \tau_d \vdash t :\tau}(x_1, \cdots, x_d; y)$.

Note however that when we go beyond the term language of $\ERC_0$, 
there are real functions that can be defined in the programming language but not in our specification language:
\begin{exa}
\label{exa:Exp3}
The restricted exponential function $\exp:I=[0;1]\to\IR$
is uniquely characterized by the following formula:
\begin{gather}
\all{ x,y \in I}  \; x+y\in I\Rightarrow \exp(x+y)=\exp(x)\cdot \exp(y) \label{eqn:Exp4}   \\
\all{x,y\in I}  |\exp(x)-\exp(y)| \leq 3\cdot |x-y| \nonumber \\ 
\exp(1) \;=\; \lim_n (1+1/n)^n \;=\; \sum\nolimits_n 1/n! \label{eqn:Exp3}
\end{gather}
The first line is the well-known functional equation,
and the second one captures Lipschitz-continuity.

Although $\exp$ can be realized
as in Section~\ref{ss:Exp},
the defining Equation~(\ref{eqn:Exp3}) \emph{exceeds} the first-order language of \ERC by involving the transcendental constant $e$, which cannot be characterized algebraically.
\qed\end{exa}
In other words, our logic is not expressive enough for the command language of $\ERC$.
Propositional formulae in the specification language
can only define semi-algebraic subsets of Euclidean space;
and, according to Tarski-Seidenberg, real quantification does not
increase the expressive power. Integer quantification
can define countable unions of semi-algebraic subsets,
but no more according to Theorem~\ref{t:Eike}(b).             
According to Lindemann-Weierstrass, the graph of $\exp:[0;1]\to\IR$
is no countable union of semi-algebraic Euclidean sets,
hence impossible to \emph{define} in $\Structure_0$.

The expansion $\Structure\big(\{\exp\},\emptyset\big)$ on the other hand
makes $\exp$ trivial to define---but its first-order theory may violate decidability.
Such trade-offs are unavoidable according to Remark~\ref{r:Definability}.
Nonetheless,
specification and formal verification may suffice with less than
definability of the function under consideration: 
applications tend to be interested in solutions that satisfy 
given algebraic properties expressible in \ERC---such 
as the exponential functional Equation~(\ref{eqn:Exp4})---but 
do not necessarily make them unique, particularly 
transcendental or multi-valued ones.

\subsection{Hoare Logic}
\label{ss:Hoare}
Hoare logic is a well-known tool for formally proving the total correctness of a
program and agreement with the problem specification.
The following considerations are guided by \cite[\S3]{reynolds2009theories},
adapted and extended to \ERC with its three-sorted
structure and multi-valued semantics. Both complicate matters since, for instance,
a real guard variable in a while loop may strictly decrease during each iteration
yet remain bounded forever; furthermore merely evaluating the loop condition can
cause lack of termination when real equality occurs; see Remark~\ref{rem:Hoare}.
Since our language is simple imperative,
we adopt the following notion of \emph{total correctness specification}:

\begin{defi}
For a well-typed command $S$ in 
$\ERC(\calF, \calG)$ with $\Gamma \vdash S \triangleright \Gamma'$,
a (total correctness) \emph{specification} 
\[
\Gamma \vDash \tassert{\phi}  \;S\; \tassert{\psi} \triangleright\Gamma'
\]
is defined by 
\[
\all{\sigma \in \sem{\Gamma}} \bot \not\in \sem{S}\sigma \land \all{\delta\in\sem{S}\sigma} \psi(\delta).
\]
Here, $\phi$, $\psi$ are formulae in the specification language of $\ERC(\calF, \calG)$.
In the \emph{precondition} $\phi$, only the variables in $\Gamma$ appear free and
in the \emph{postcondition} $\psi$, only the variables in $\Gamma'$ appear free.
The notation says, for any $\sigma\in \sem{\Gamma}$ that satisfies $\phi$,
(i) $\bot \not \in \sem{S}\sigma$ and (ii) any $\delta \in\sem{S}\sigma$
satisfies $\psi$.
\end{defi}
\noindent Note that the definition of our specification does not publish any obligation to $S$ when the precondition $\phi$ is not met.

\begin{defi} 
\label{defi:Hoare_ERC}

A \emph{Hoare triple} for $\ERC(\calF, \calG)$ is of the form 
$\Gamma \vdash \tassert{\phi}\;S\;\tassert{\psi} \triangleright \Gamma'$
where $S$ is a well-typed command in $\ERC(\calF, \calG)$ such that $\Gamma \vdash S \triangleright \Gamma'$,
and $\phi, \psi$ are formulae in the specification 
language of $\ERC(\calF, \calG)$
such that
only the variables in $\Gamma$ are free in $\phi$ and only
the variables in $\Gamma'$ are free in $\psi$.
Hoare logic of $\ERC(\calF, \calG)$ is a formal system which
consists of the inference rules and axioms for constructing Hoare triples defined
in \autoref{fig:hoare}.
\qed
\end{defi}

\begin{figure*}[!ht]
\caption{The Hoare logic of $\ERC(\calF, \calG)$ \label{fig:hoare}}
\fbox{\noindent\begin{minipage}{0.97\textwidth}

\begin{mathpar}
\infer[\phi \Rightarrow \phi' \text{ and } \psi' \Rightarrow \psi]
    {
        \Gamma \vdash \tassert{\phi}\;S\;\tassert{\psi} \triangleright \Gamma'
    }
    {%
    \Gamma \vdash \tassert{\phi'}\;S\;\tassert{\psi'}\triangleright \Gamma'
    }
    \and
\infer
    {
        \Gamma \vdash \tassert{\psi} \;\eskip\;\tassert{\psi}  \triangleright \Gamma
    }
    {%
    }
\\
\infer
    {
        \Gamma \vdash \tassert{\some{w} \trans{t}(w)
                \land \all{w} \trans{t}(w)\Rightarrow \psi[w/x]\,}\;
         x\dfeq  t\;
        \asser{\psi} \triangleright \Gamma
    }
    {%
    }
\\
\infer
    {
        \Gamma \vdash
        \asser{\mathsf{ArrPre}(\psi,x,m,t)}
        \;x[m]  \dfeq  t\;
        \asser{\psi}
        \triangleright \Gamma
    }
    {%
    }
\\
\infer
    {
        \Gamma \vdash \tassert{\some{w} \trans{t}(w)
                \land \all{w} \trans{t}(w)\Rightarrow \psi[w/x]}\;
         \elet\;x : \tau = t\; \tassert{\psi} \triangleright \Gamma, x:\tau
    }
    {%
    }
\\
\infer
    {
        \Gamma \vdash\tassert{\phi} \; S_1;S_2 \;\tassert{\psi} \triangleright \Gamma_2
    }
    {
        \Gamma \vdash\asser{\phi}\; S_1 \;\asser{\theta} \triangleright \Gamma_1
        &
        \Gamma_1 \vdash \tassert{\theta} \;S_2 \; \asser{\psi} \triangleright \Gamma_2
    }
\\
\infer
    {
        \Gamma \vdash \asser{\phi \land (\trans{b}(\ctt) \lor \trans{b}(\cff)) \land \neg\trans{b}(\cuu)}
         \;
            \eif\;b\;\ethen\;S_1 \;\eelse\;S_2 \;\eendif \;
       \asser{\psi}            \triangleright \Gamma
}
    {
        \Gamma \vdash\asser{\phi\land\; \trans{b}(\ctt)}
        \; S_1 \;
        \asser{\psi}\triangleright \Gamma_1
        &
        \Gamma \vdash\asser{\phi\land \trans{b}(\cff)}
        \; S_2 \;
        \asser{\psi}\triangleright \Gamma_2
    }
\\
\infer
    {
        \Gamma \vdash \asser{I}\;
            \ewhile\;b\;\edo\;S \;\eendwhile \;
        \asser{I \land \trans{b}(\cff)} \triangleright \Gamma
    }
    {
                    \Gamma ,\xi, \xi' : \dreal \vdash \asser{\trans{b}(\ctt) \land I\land V = \xi
                    \land L = \xi'
                    }
                    \; S \;
                    \asser{I\land V \leq \xi - \xi' \land
                    L = \xi'
                    } \triangleright \Gamma''
    }
\\
\end{mathpar}
\noindent The rule for loops has the side conditions:
\begin{align*}
I \land \trans{b}(\ctt) &\Rightarrow L > 0\;,\\
I &\Rightarrow (\trans{b}(\ctt) \lor \trans{b}(\cff))\land \neg\trans{b}(\cuu)\;,
\\
I\land V \leq 0 &\Rightarrow\all{k}\trans{b}(k)\Rightarrow k=\cff \;,\text{ and }\\
\xi, \xi'& \text{ do not appear free in } I, V,L\;.
\end{align*}
\noindent where $L$ and $V$ are real-valued.

\vspace{1em}

\noindent In the case of array assignment
$\mathsf{ArrPre}(\psi,x,m,t)$ is defined as follows:
\begin{align*}
    \mathsf{ArrPre}(\psi,x,m,t) :\equiv \; & \big(\some{i, w} \trans{m}(i)\land\trans{t}(w)\big) \\
                & \land\, \all{i} (\trans{m}(i)\Rightarrow 0\leq i < d) \\
                & \land\, \all{i, w} (\trans{m}(i) \land \trans{t}(w) \land \all{y} \semupdate(x, i, w, y) \Rightarrow
                \psi[\subst{x}{y}])
\end{align*}
assuming $\Gamma \vdash x[m]  \dfeq  t\; \triangleright \Gamma$ and
$\Gamma \vdash x : \dreal[d]$ where
\begin{align*}
\semupdate(x, i, w, y) :\equiv \; & 0\leq i < d \land 
((y_{(0)} = x_{(0)} \land i \neq 0) \lor (y_{(0)} = w \land i = 0))\\
&\land \cdots \land ((y_{(d-1)} = x_{(d-1)} \land i \neq d-1) \lor (y_{(d-1)} = w \land i = d-1))
\end{align*}    
\end{minipage}}
\end{figure*}

The purpose of (Hoare) logic is to replace semantic arguments with formal proofs:
sequences of purely syntactic manipulations, starting with the axioms and following 
certain inference rules, that for example, a computer can verify.
Classical Hoare logic contains one exception:

\begin{rem}
\label{r:SideCondition}
The \emph{rule of consequence} for
precondition-strengthening and postcondition-weakening 
\[ \infer[\phi \Rightarrow \phi' \text{ and } \psi' \Rightarrow \psi]
    {
        \Gamma \vdash \tassert{\phi}\;S\;\tassert{\psi} \triangleright \Gamma'
    }
    {%
    \Gamma \vdash \tassert{\phi'}\;S\;\tassert{\psi'}\triangleright \Gamma'
    }
\] 
depends on the semantic side-conditions
$\phi \Rightarrow \phi'$  and  $\psi' \Rightarrow \psi$
which may or may not be feasible to verify algorithmically.
Over integers, algorithmic verifiability can fail according to G\"{o}del \cite[\S6]{Cook78}, 
but not in our specification language of $\ERC_0$ according to Theorem~\ref{t:Eike}(a).
\end{rem}


\begin{rem} 
\label{rem:Hoare}
\hfill
\begin{enumerate}

\item In the axiom for assignments,
recall from Definition~\ref{d:trans} that the precondition $\some{w}\trans{t}(w)$ ensures that
the denotation of $t$ is well-defined. And,
$\all{w}\trans{t}(w) \Rightarrow \psi[w/x]$
says that for each value in the denotation of $t$,
$\psi$ holds when we replace the variable $x$ with the value.

\item In the \WHILE~loop case,
    \begin{enumerate}
    \item the formula $I$ is the loop invariant and
    the term $V$ is the loop variant.
    The term $L$ is some invariant quantity that bounds by how much 
    $V$ decreases in each iteration.
    \item
    $\xi, \xi'$ are ghost variables that do not appear in
    $\Gamma$. They can be understood as meta-level universally quantified variables.
    \item
    The side condition says (i) each loop decreases
    $V$ by some positive invariant quantity $L$;
    (ii) as long as $I$ holds, the evaluation of
    $b$ is either $\ctt$ or $\cff$ (but not $\cuu$); and
    (iii) when $V$ is negative, it is guaranteed
    that the evaluation of $b$ is $\cff$.    
    \end{enumerate}
\item In the rule of $\IF$ conditionals,
    the precondition
    $(\trans{b}(\ctt)\lor \trans{b}(\cff)) \land \neg \trans{b}(\cuu)$
    says that the evaluation of
    $b$ is either $\ctt$ or $\cff$ (but not $\cuu$).
\item In the rule of array assignments,
$\semupdate(x, i, w, y)$ says $x = (x_{(0)}, \cdots, x_{(d-1)})$ is updated 
to $y = (y_{(0)}, \cdots, y_{(d-1)})$ by assigning $w$ at $0 \leq i < d$.
\qed
\end{enumerate}
\end{rem}

\begin{thm}
\label{t:hoaresound}
The Hoare logic of $\ERC(\calF, \calG)$ is sound; i.e., for any Hoare triple $\Gamma\vdash \tassert{\phi}\;S\;\tassert{\psi}\triangleright\Gamma'$,
it holds that $\Gamma\vDash \tassert{\phi}\;S\;\tassert{\psi}\triangleright\Gamma'$.

\end{thm}
\begin{proof}
See Appendix~\ref{append:hoare_sound}.
\end{proof}

\section{Example Formal Verification in \emph{Exact Real Computation}} \label{ss:Trisection}
The present section picks up from Section~\ref{ss:Root}
to illustrate formal verification in \ERC.
To emphasize, our purpose here is not to actually establish
correctness of the long-known Trisection method,
but to demonstrate the extended Hoare logic from Section~\ref{ss:Hoare} using a toy example.
Since Trisection relies on the Intermediate Value Theorem,
any correctness proof must make full use of real
(as opposed to floating point, rational, or algebraic) numbers.

Let us define some abbreviations such that the program in \autoref{fig:Trisection}
becomes of the form \[\blim{p}a : \dreal, b : \dreal\;S \;\blimreturn{a}{p}.\]
\begin{eqnarray*}
\tilde{t}_1 &:\equiv \;& b-a \elt \myprec{p}  \;,\;\;\tilde{t}_2 \dfeq \; \myprec{p-1} \elt b-a  \\
t_1 &:\equiv \;&  f(b/3+2\mult a/3)\mult f(b) \elt 0 \;,\;\; t_2 :\equiv f(a)\mult f(2\mult b/3 + a/3) \elt 0\\
b_1 &:\equiv \;& \mychoose(\tilde t_1, \tilde t_2)  \ieq  1,\;\;b_2 \dfeq \;
\mychoose(t_1, t_2)  \ieq  1 \\
S &:\equiv\;& \ewhile\; b_1
  		\; \edo \; S_1 \; \eendwhile \\[-0.5ex]
S_1 &:\equiv\;& \eif\; b_2  \;\ethen\; S_2\; \eelse\; S_3\; \eendif \\[-0.5ex]
S_2 &:\equiv\;& b \dfeq  2\mult b/3 + a/ 3 \\[-0.5ex]
S_3 &:\equiv\;& a \dfeq  b/3 + 2\mult a /3
\end{eqnarray*}

We work with $\Theory(\{f\}, \emptyset)$ 
which contains sentences saying that $f$ is a continuous function.
We want to verify that the program's denotation at its inputs $a, b$ which isolate a root of $f$ uniquely with a sign change:
\[
\uniq(f, a, b) :\equiv 
\cont(f, a, b) \land
f(a)\cdot f(b) < 0 \land
\usome{x} a <x<b \land f(x) = 0\;.
\]
Considering the limit taken at the end of every real program, we need to prove the
specification:
\[
\Gamma \vDash \tassert{p = p' \land\uniq(f, a, b)}\; S\; \tassert{\usome{z}
f(z) = 0\land a<z<b \land |a - z| \leq 2^{p'}} \triangleright \Gamma'
\]
where $\Gamma = p, p' : \dint, a, b : \dreal$ and $\Gamma' = p, p' : \dint, a, b : \dreal$.
The ghost variable $p'$ captures the initial value that the variable $p$ stores, considering that the value $p$ stores may throughout the computation (though it does not in this specific example.)

To simplify our presentation, let us write $t_1, t_2, \tilde{t}_1, \tilde{t}_2$ as valid formulae in our specification language where their occurrences of $\elt$ are implicitly replaced with $<$.
See that $\trans{b_1}(\ctt) \Leftrightarrow \tilde t_2$,
$\trans{b_1}(\cff) \Leftrightarrow \tilde t_1$,
$\neg\trans{b_1}(\cuu)$,
$\trans{b_2}(\ctt) \Leftrightarrow t_2$, $\trans{b_2}(\cff) \Leftrightarrow t_1$, and
$\neg\trans{b_2}(\cuu)$ hold.
Let us define $I\coloneqq p = p'\land \uniq(f,a,b)$ as a candidate for the loop invariant, $V\dfeq b-a-2^{p-1}$ as a candidate for the loop variant,
$L \dfeq 2^{p-2}$ be a candidate for a lower bound decrements,
$\tilde{P} \coloneqq \tilde t_2 \wedge I\wedge V=\xi \land L = \xi'$,
and $\tilde{Q} \coloneqq I \wedge V\leq \xi- \xi' \land L = \xi'$
in our \specification language with variables $\xi, \xi'$ of the sort $\IR$.
Let $\Delta \coloneqq p , p':\dint, a, b, \xi, \xi': \dreal$.

From the axiom for assignments, we have the triples:
\[
\Delta\vdash \tassert{\some{\omega} \trans{a/3 + 2\mult b/3}(\omega)
\land \all{\omega}\trans{a/3 + 2\mult b/3}(\omega)\Rightarrow
\tilde{Q}[\omega/b]}\;S_2\;\tassert{\tilde{Q}} \triangleright \Delta
\]
\[
\Delta\vdash \tassert{\some{\omega} \trans{2\mult a/3+ b/3}(\omega)
\land \all{\omega} \trans{2\mult a/3 +. b/3}(\omega)\Rightarrow
\tilde{Q}[\omega/a]}\;S_3\;\tassert{\tilde{Q}} \triangleright \Delta
\]
See that we can apply the rule of precondition weakening to get the following triples derived:
\[
\Delta\vdash \tassert{
\tilde{Q}[(a/3 + 2\mult b/3) / b]}\;C_2\;\tassert{\tilde{Q}} \triangleright \Delta,\quad
\Delta\vdash \tassert{
\tilde{Q}[(2\mult a/3+ b/3) / a]}\;C_3\;\tassert{\tilde{Q}} \triangleright \Delta.
\]

When we unwrap abbreviations, we can verify the following equivalences:
\begin{align*}
\tilde Q[(2\mult a/3+ b/3) / a] &\;\Leftrightarrow\; p = p' \land \uniq(f, (2\mult a/3+ b/3), b) \\
&\;\;\land\;\; b - (2\mult a/3+ b/3) - 2^{p-1} \leq \xi - \xi' \\
&\;\;\land\;\; 2^{p-2} = \xi' \\
\tilde P \land t_1 & \;\Leftrightarrow\; p = p' \land \uniq(f, a, b) \land  f(b / 3 + 2\times a / 3) \times f(b) < 0\\
&\;\;\land\;\;  2^{p-1} < b - a \land  b - a -2^{p-1} = \xi  \\
&\;\;\land\;\; 2^{p-2} =\xi'
\end{align*}

Due to the intermediate value theorem, 
if an interval $(a;b)$ contains
a root of $f$ uniquely,
and if $f(x)\mult f(y) < 0 $ for $a\leq x < y \leq b$,
then $(x; y)$ also contains the root of $f$ uniquely.
Hence, $\tilde P \land t_1 \Rightarrow \tilde Q[(2\mult a/3+ b/3) / a]$ holds.
And, similarly, 
$\tilde P \land t_2 \Rightarrow \tilde Q[(a/3+ 2\mult b/3) / b]$ holds.

Therefore, by the rule of precondition strengthening on the triples of $S_2, S_3$,
and the rule for conditionals, we get the triple:
\[
\Delta\vdash\tassert{\tilde{t}_2\wedge I\wedge \big(V=\xi\big) \wedge L = \xi'}\;S_1\;\tassert{ I \wedge(V\leq \xi-\xi') \wedge L = \xi'}\triangleright \Delta
\]
The side-conditions of the rule for while loops are quite trivial.
Hence, assuming that they are proven, we apply the rule of while loops,
and we get the following triple:
\[
\Gamma \vdash \tassert{I}\;S\;\tassert{I \land \tilde t_2} \triangleright
\Gamma'
\]
Using the rule of pre/postcondition strengthening/weakening, we can get the originally desired specification.

\section{Conclusion and Future work} \label{s:Conclusion}
We have formalized an imperative programming language for exact real number computation.
Its domain-theoretic denotational semantics, based on the Plotkin powerdomain modelling multi-valued computations, is carefully designed to be computable and complete, matching the intuition of operating on continuous data exactly without rounding errors and in agreement with (proofs in) calculus.
This enables a natural approach to formal program verification by adding real number axioms to Hoare logic.

The following considerations naturally suggest future further investigations:
\begin{itemize}
\item \textbf{Computational Cost:}
After the design of an algorithm comes its analysis,
in terms of computational cost as a quantitative indicator
of its practical performance. For realistic predictions,
real complexity theory \cite{Ko91,Wei03} employs bit-cost,
as opposed to unit cost common in algebraic complexity theory \cite{ACT}.
In \cite[Definition~2.4]{BH98} it is suggested that a logarithmic cost measure
where each operation is supposed to take time according to the
binary length of the integer (part of the real) to be processed.
More accurate predictions take into account the precision parameter $p$ for real programs
and for real number comparisons ``$x\elt y$''
the logarithm of the difference $\myabs{x - y}$.
\item \textbf{Full Mixed Data Types:}
\ERC as introduced here formalizes computing with data types $\INTEGER$, $\REAL$, and $\LOGIC$
as counterparts to mathematical $\IZ,\IR,\Kleene$ where $\LOGIC$ is permitted only for expressions and local variables.
A future extension will include also (multi-)functions with $\LOGIC$ type arguments and return values as well as arrays and a dedicated limit operator.
\item \textbf{Explicit Limit Operator:}
Constructing real numbers as the limits of converging sequences is done only at the level of programs in ERC. It was our design choice to separate the term language of ERC with loops and of course limits. A direction of extending our framework to allow constructing a limit within user
programs as a term has been also investigated \cite{clericalcca} where the authors suggest an expression-based programming language called Clerical that provides an explicit limit operator. Its prototype implementation can be found in \cite{BauerClerical}.
\item \textbf{Multi-valued Real Functions:}
The present version of \ERC formalizes computing mappings from reals to integers
in the multi-valued sense because any single-valued,
and necessarily continuous \cite[Theorems~4.3.1+3.2.11]{Wei00},
function with connected domain and discrete range must be constant.
On the other hand, computing real values is deliberately restricted to the single-valued case.
Defining approximate computation of real multi-valued functions
is delicate and still under exploration \cite{Mueller18,Konecny18}.

\item \textbf{Function\emph{al}s and Operators:}
The program in \autoref{fig:Trisection} mimics a function\emph{al} that receives a continuous function $f$ as an argument, accessible by point-wise black box evaluation.
Extending \ERC with function arguments enriched with quantitative continuity information, such as a modulus of continuity \cite{KC12}, is necessary to extend Theorem~\ref{t:Complete} (Turing-completeness of \ERC)
from functions to functionals.

\item \textbf{Automated Formal Verification:}
The decidability of the theory of the base language $\ERC_0$ according to Theorem~\ref{t:Eike}
includes its complete (and elegant) axiomatization:
Guaranteed to yield convenient \emph{automated} formal verification. 
This includes formalizing $\ERC$ itself in a formal system for example in the Coq proof assistant \cite{DBLP:conf/fsttcs/Konecny0T20,DBLP:journals/lmcs/SteinbergTT21} and further development of annotated ERC and verification condition extraction mechanisms. 

\item \textbf{Continuous Data Types beyond the Reals:}
Real computability theory has been extended to
topological T$_0$ spaces, real complexity theory
to co-Polish spaces \cite{Metric,Sch04}. Current
and future works similarly extend \ERC to continuous
data types beyond real numbers/functions,
such as the Grassmannian, tensors \cite{Tensors} and groups \cite{Haar}. 
As hinted in Remark~\ref{r:cc}, this may require extending our language with countable nondeterminism for the extended language to be complete as well.
\end{itemize}

\section*{Acknowledgement}
\noindent 
This work was supported by the National Research Foundation of Korea
  (grants NRF-2017R1E1A1A03071032 and NRF-2017R1D1A1B05031658),
  by the International Research \& Development Program of
  the Korean Ministry of Science and ICT\\ (grant NRF-2016K1A3A7A03950702),
  by the
  \protect\raisebox{-1pt}
  {\protect\includegraphics[height=8pt]{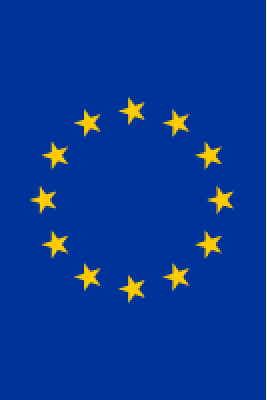}}~European Union's
  Horizon 2020 research and innovation
  program under the Marie Sk{\l}odowska-Curie grant agreement No 731143,
  and by the German Research Foundation (DFG) Grant MU 1801/5-1.

We thank Andrej Bauer, Cyril Cohen, Jeehoon Kang, Johannes Kanig, Alex Simpson, and Hongseok Yang for seminal discussions.
  Anonymous referees have provided valuable feedback to an earlier version.
We would like to thank Holger Thies for constructive comments on this work and proofreading.

  A preliminary version of this work was submitted as a part of the first author's PhD thesis \cite{sewonphd}.

\bibliographystyle{alphaurl}
\bibliography{lmcs_final}

\appendix
\section{Proof of the Soundness of the Hoare Logic of \ERC}
\label{append:hoare_sound}
We start the proof with the lemma:

\begin{lem}
    For a well-typed command $\Gamma\vdash\ewhile\;b\;\edo\;S\;\eendwhile \triangleright \Gamma$,
    define the sequences of set-valued functions on $\sem{\Gamma}$:
    \begin{itemize}
    \item
    $B_{b, S}^0(\sigma) :\equiv \{\sigma\}$
    \item
    $C_{b, S}^0(\sigma) :\equiv \emptyset$
    \item
    $B_{b, S}^{n+1}(\sigma):\equiv \bigcup_{\delta \in  B_{b, S}^n(\sigma)}
    \bigcup_{\substack{l \in \sem{b}\delta\\\delta'\in \restriction_\Gamma^{\dagger^\dagger}\fcomp\sem{S}\delta}}
    \begin{cases}
    \{\delta'\} &\text{if } l = \ctt \land \delta' \neq \bot,\\
    \emptyset &\text{otherwise}.
    \end{cases}$
    \item
    $C_{b, S}^{n+1}(\sigma) :\equiv
    C_{b, S}^n(\sigma)\cup\bigcup_{\delta \in B_{b, S}^n(\sigma)}
    \bigcup_{\substack{l \in \sem{b}\sigma\\ \delta'\in \restriction_\Gamma^{\dagger^\dagger}\fcomp\sem{S}\delta}}
    \begin{cases}
    \{\delta\} &\text{if } l = \cff, \\
    \emptyset &\text{if } l = \ctt \land \delta' \neq \bot,\\
    \{\bot\} &\text{otherwise}.
    \end{cases}$

    \end{itemize}
    Then, for all $n\in\IN$, it holds that
    $\mathcal{W}^{(n)}(\sigma) = C_{b, S}^n(\sigma) \cup \{\bot \mid \some{x} x\in B_{b, S}^n(\sigma) \}$.

    Intuitively, $B^n_{b, S}(\sigma)$ is the set of states that requires further execution after
    running the while loop on $\sigma$ for $n$ times, and $C^n_{b, S}(\sigma)$ is the set of states that have
    escaped from the loop (either because $\cff$ has been evaluated or $\bot$ has occurred)
    during running the loop for $n$ times.
\end{lem}
\begin{proof}
Let us drop the subscripts $b, S$ for the convenience in the presentation
and write $\mathbf{S}$ for $\restriction_\Gamma^{\ddagger^\dagger}\fcomp\sem{S}$.

We first prove the following alternative characterization of the sequence of
sets:
\[
B^{n+1}(\sigma) = \bigcup_{
\substack{\ell \in \sem{b}\sigma\\\delta \in \mathbf{S}\sigma}}
\begin{cases}
B^n(\delta) &\text{if } \ell = \ctt \land \delta \neq \bot, \\
\emptyset &\text{otherwise}.
\end{cases}
\]

It is trivial when $n = 0$. 
Now, suppose the equation holds for all $\sigma$ and for all $n$ up to $m$. 
Then the following derivation shows that the characterization is valid for $n = m + 1$ as well.
\begin{align*}
B^{m+2}(\sigma) &= \bigcup_{\gamma \in B^{m+1}(\sigma)}
\bigcup_{\substack{\ell\in\sem{b}\gamma \\ \delta \in \mathbf{S}\gamma}}
\begin{cases}
\{\delta\} &\text{if } \ell = \ctt \land \delta \neq \bot, \\
\emptyset & \text{otherwise}.
\end{cases}\\
&= \bigcup_{
  \gamma \in
\bigcup_{\substack{\ell'\in\sem{b}\sigma \\ \delta' \in \mathbf{S}\sigma}}
\begin{cases}
B^m(\delta') &\text{if } \ell' = \ctt \land \delta' \neq \bot, \\
\emptyset &\text{otherwise}.
\end{cases}
}
\hspace*{8pt}
\bigcup_{\substack{\ell\in\sem{b}\gamma \\ \delta \in \mathbf{S}\gamma}}
\begin{cases}
\{\delta\} &\text{if } \ell = \ctt \land \delta \neq \bot, \\
\emptyset & \text{otherwise}.
\end{cases}\\
&=
\bigcup_{\substack{\ell' \in \sem{b}\sigma \\ \delta' \in \mathbf{S}\sigma}}
\begin{cases}
\bigcup_{\gamma \in B^m(\delta')}\bigcup_{\substack{\ell \in \sem{b}\gamma \\ \delta \in \mathbf{S}\gamma}}\begin{cases}
\{\delta\} &\text{if } \ell = \ctt \land \delta \neq \bot, \\
\emptyset &\text{otherwise}.
\end{cases} &\text{if } \ell' = \ctt \land \delta' \neq \bot, \\
\emptyset &\text{otherwise}.
\end{cases} \\
&= \bigcup_{\substack{\ell'\in\sem{b}\sigma \\ \delta' \in \mathbf{S}\sigma}}
\begin{cases} B^{m+1}(\delta') &\text{if } \ell' = \ctt \land \delta' \neq \bot, \\
\emptyset &\text{otherwise}.
\end{cases}
\end{align*}

We now show the following characterization:
\[
C^{n+1}(\sigma) = \bigcup_{\substack{\ell\in \sem{b}\sigma \\ \delta\in\mathbf{S}\sigma}}
\begin{cases}
C^n(\delta) &\text{if } \ell = \ctt \land \delta \neq \bot,\\
\{\sigma \} &\text{if } \ell = \cff, \\
\{\bot\} &\text{otherwise}.
\end{cases}
\]
It is easy to show that the equation holds for $n = 0$.
Now, assume the equation holds for all $n$ up to $m$.
Then,

\begin{align*}
C^{m+2}(\sigma)
&= C^{m+1}(\sigma) \cup \bigcup_{\delta \in B^{m+1}(\sigma)} \bigcup_{\substack{\ell' \in\sem{b}\delta \\ \delta'\in\mathbf{S}\delta}}\begin{cases}
\{\delta\} &\text{if } \ell' = \cff, \\
\emptyset &\text{if } \ell' = \ctt \land \delta' \neq \bot, \\
\{\bot\} &\text{otherwise}.
\end{cases}\\
&= C^{m+1}(\sigma) \cup \\ 
&\;\;\bigcup_{\delta \in \bigcup_{\substack{\ell \in \sem{b}\sigma \\ \gamma\in\mathbf{S}\sigma}}
\begin{cases}
B^m(\gamma) &\text{if } \ell = \ctt \land \gamma \neq \bot, \\
\emptyset &\text{otherwise}.
\end{cases}}
\bigcup_{\substack{\ell' \in\sem{b}\delta \\ \delta'\in\mathbf{S}\delta}}\begin{cases}
\{\delta\} &\text{if } \ell' = \cff, \\
\emptyset &\text{if } \ell' = \ctt \land \delta' \neq \bot, \\
\{\bot\} &\text{otherwise}.
\end{cases}\\
&=C^{m+1}(\sigma) \cup \phantom{x}\\
&\mathrel{\phantom{=}}\bigcup_{\substack{\ell\in\sem{b}\sigma \\ \gamma \in \mathbf{S}\sigma}}
\begin{cases}
\bigcup_{\delta\in B^m(\gamma)}\bigcup_{\substack{\ell\in\sem{b}\delta \\ \delta' \in \mathbf{S}\delta}}
\begin{cases}\{\delta\} &\text{if } \ell' = \cff \\
\emptyset &\text{if } \ell' = \ctt \land \delta' \neq \bot, \\
\{\bot\} &\text{otherwise}.
\end{cases} &\text{if } \ell = \ctt \land \gamma \neq \bot, \\
\emptyset &\text{otherwise}.
\end{cases}\\
&=
\bigcup_{\substack{\ell\in\sem{b}\sigma \\ \gamma \in \mathbf{S}\sigma}}
\begin{cases}
C^m(\gamma)\cup & \text{if } \ell = \ctt \land \gamma \neq \bot, \\
\;\;\bigcup_{\delta\in B^m(\gamma)}\bigcup_{\substack{\ell\in\sem{b}\delta \\ \delta' \in \mathbf{S}\delta}}
\begin{cases}\{\delta\} &\text{if } \ell' = \cff \\
\emptyset &\text{if } \ell' = \ctt \land \delta' \neq \bot, \\
\{\bot\} &\text{otherwise}.
\end{cases} &\\
\{\sigma\} &\text{if } \ell = \cff, \\
\{\bot\} &\text{otherwise}.
\end{cases} \\
&=
\bigcup_{\substack{\ell\in\sem{b}\sigma \\ \gamma \in \mathbf{S}\sigma}}
\begin{cases}
C^{m+1}(\gamma) &\text{if } \ell = \ctt \land \gamma \neq \bot, \\
\{\sigma\} &\text{if } \ell = \cff, \\
\{\bot\} &\text{otherwise}.
\end{cases}
\end{align*}

Using the suggested characterization, we prove
    $\mathcal{W}^{(n)}(\sigma) = C_{b, S}^n(\sigma) \cup \{\bot \mid \some{x} x\in B_{b, S}^n(\sigma) \}$
    for all $n\in\IN$.
    When $n=0$, both are $\{\bot\}$.
    Suppose the equation holds for $n=m$.
    Then,
\begin{align*}
\mathcal{W}^{(n+1)}(\sigma) &=
\bigcup_{\substack{\ell \in \sem{b}\sigma \\ \delta \in \mathbf{S}\sigma}}
\begin{cases}
\mathcal{W}^{(n)}(\delta) &\text{if }\ell = \ctt \land \delta \neq \bot,\\
\{\sigma\} &\text{if } \ell = \cff,\\
\{\bot\} &\text{otherwise}.
\end{cases}\\
&=
\bigcup_{\substack{\ell \in \sem{b}\sigma \\ \delta \in \mathbf{S}\sigma}}
\begin{cases}
C^m(\delta) &\text{if }\ell = \ctt \land \delta \neq \bot,\\
\{\sigma\} &\text{if } \ell = \cff,\\
\{\bot\} &\text{otherwise}.
\end{cases}
\\
&\qquad\cup\bigcup_{\substack{\ell \in \sem{b}\sigma \\ \delta \in \mathbf{S}\sigma}}
\begin{cases}
\{\bot | \some{\gamma} \gamma \in B^m(\delta)\} &\text{if }\ell = \ctt \land \delta \neq \bot\\
\emptyset &\text{if } \ell = \cff,\\
\emptyset &\text{otherwise}.
\end{cases}\\
&= C^{m+1}(\sigma) \cup \left\{\bot \left| \some{\gamma}\gamma \in
\bigcup_{\substack{\ell \in \sem{b}\sigma \\ \delta \in \mathbf{S}\sigma}}
\begin{cases}
B^m(\delta) &\text{if }\ell  = \ctt \land \delta \neq \bot, \\
\emptyset &\text{otherwise}.
\end{cases}\right.\right\}\\
&= C^{m+1}(\sigma) \cup \{\bot \mid \some{\gamma} \gamma \in B^{m+1}(\sigma)\}\qedhere
\end{align*}
\end{proof}

We now proceed to the proof of Theorem~\ref{t:hoaresound}:
\begin{proof}
We can prove the statement by checking the soundness of each rule.

\begin{enumerate}
    \item (Assignment):

    Consider any state $\sigma$ which satisfies
    $\some{ w} \trans{t}(w)
                \land \all{w} \trans{t}(w)\Rightarrow \psi[\subst{x}{w}]$.
    Then, $\bot \not\in \sem{t}\sigma$ and
    for any $w \in \sem{t}\sigma$, $\psi[w/x]$ holds.

    Now, see that $\sem{x \dfeq  t}\sigma = \bigcup_{w \in \sem{t}\sigma}\{\sigma[x \mapsto w]\}$
    since $\bot \not \in \sem{t}\sigma$ and for all $w \in \sem{t}\sigma$, $\sigma[x \mapsto w]$ satisfies
    $\psi$.

    \item The rule variable declarations and the rule of array assignments can be verified in very similar manner as above
    and the rules of pre/postcondition strengthening/weakening, skip, and
    sequential compositions can be verified trivially.

    \item (Conditional):

    Consider any state $\sigma$ which satisfies
    $\phi \land (\trans{b}(\ctt) \lor \trans{b}(\cff)) \land \neg \trans{b}(\cuu)$.
    Then, $\sem{b}\sigma = \{\ctt, \cff\}, \{\ctt\}$, or $\{\cff\}$.
    Let us check the three cases:
    \begin{enumerate}
        \item when $\sem{b}\sigma = \{\ctt, \cff\}$:

        Then, $\sigma$ satisfies $\phi\land \trans{b}(\ctt)$ and $\phi\land\trans{b}(\cff)$.
        Therefore, (i) $\bot\not\in\sem{S_1}\sigma$, (ii) for all $\delta \in \sem{S_1}\sigma$ it holds that
        $\delta \vDash \psi$, (iii) $\bot\not\in\sem{S_2}\sigma$, and
        (iv) for all $\delta \in \sem{S_2}\sigma$ it holds that $\delta \vDash \psi$.

        Since $\bot \not\in\sem{b}\sigma$ and $\cuu \not\in \sem{b}\sigma$, the denotation becomes
        \[\sem{\eif\;b\;\ethen\;S_1\;\eelse\;S_2\;\eendif}\sigma =
        \restriction_\Gamma^{\ddagger^\dagger}(\sem{S_1}\sigma) \cup \restriction_\Gamma^{\ddagger^\dagger}(\sem{S_2}\sigma).
        \]
        Since the restriction operator does not create $\bot$, the denotation does not contain $\bot$.
        Also, since $\psi$ only consists of free variables that are in $\Gamma$, each resulting state after the restriction still satisfies $\psi$.

        \item when $\sem{b}\sigma = \{\ctt\}$:

        Then, $\sigma$ satisfies $\phi\land\trans{b}(\ctt)$.
        Hence, (i) $\bot\not\in\sem{S_1}\sigma$, (ii) each $\delta \in \sem{S_1}\sigma$ satisfies $\psi$.
        Since $\sem{b}\sigma = \{\ctt\}$, the denotation becomes
        \[\sem{\eif\;b\;\ethen\;S_1\;\eelse\;S_2\;\eendif}\sigma =
        \restriction_\Gamma^{\ddagger^\dagger}(\sem{S_1}\sigma).
        \]
        Again, since the restriction does not create $\bot$ and $\bot \not\in \sem{S_1}\sigma$,
        $\bot$ is not in the denotation.
        Since $\psi$ only consists of free variables that are in $\Gamma$, each resulting state after the restriction still satisfies $\psi$

        \item when $\sem{b}\sigma = \{\cff\}$,
        it can be done very similarly to the above item.
    \end{enumerate}

    \item (Loop):

    Consider any state $\sigma$ that satisfies
    $I$. Then, by the side-conditions,
    it also satisfies $(\trans{b}(\ctt) \lor \trans{b}(\cff)) \land \neg \trans{b}(\cuu)$.
    Hence, $\sem{b}\sigma = \{\ctt, \cff\}, \{\ctt\},$ or $\{\cff\}$ for any
    state $\sigma$ that satisfies $I$. Now, we fix a state $\sigma$ which
    satisfies $I$ hence satisfies the precondition.

    The core part of the proof is the statement:
    for any natural number $n$, it holds that
    (i) $\bot \not \in B_{b, S}^n(\sigma)$,
    (ii) $\bot \not \in C_{b, S}^n(\sigma)$,
    (iii) all $\delta$ in $B_{b, S}^n(\sigma)$
    or $C_{b, S}^n(\sigma)$ satisfies $I$, and
    (iv) all $\delta$ in $C_{b, S}^n(\sigma)$ satisfies $\trans{b}(\cff)$.

    At the moment, suppose that the above statement is true.
    Then, all we have to show is that $B^m_{b, S}(\sigma)$ becomes
    empty as $m \in\IN$ increases. Let us define
    $\ell_n \dfeq  \max\{V(\delta)\mid \delta \in B^n_{b, S}(\sigma)\}$
    and show that $\ell_n$ is strictly decreasing by some quantity
    that is bounded below, as $n$ increases. See that if it holds,
    there will be some $m$ that for all $\delta \in B^m_{b, S}(\sigma)$,
    $\sem{b}\delta = \{\cff\}$ and hence $B^{m+1}_{b, S}(\sigma) = \emptyset$.

    In order to prove it, we take the two steps:
    \begin{enumerate}
        \item
        If $B^1_{b, S}(\sigma) \neq \emptyset$, then for all $n\in\IN$ and
    for all $\delta \in B^n_{b, S}(\sigma)$, it holds that
    $L(\delta) = L(\sigma) >0$. In this case, let us write
    $\ell_0 = L(\sigma)$.

        \item
       If $B^{m+1}_{b, S}(\sigma) \neq \emptyset$,
    it holds that $\ell_{m+1} \leq \ell_m - \ell_0$.

    \end{enumerate}
    Now, we prove each statement:
    \begin{enumerate}
        \item
        $B^1_{b, S}(\sigma) \neq \emptyset$ only if
    $\ctt \in \sem{b}\sigma$ and there is some non-bottom $\delta \in \mathbf{S}$.
    Therefore, by the side-condition, $L(\sigma) > 0$.

    Suppose any $\delta \in B^{m+1}_{b, S}(\sigma)$ for any $m\in\IN$.
    See that it happens only if there is $\delta' \in B^m_{b,S}(\sigma)$
    such that $\ctt \in \sem{b}\delta'$ and $\delta \in \mathbf{S}\delta'$.
    Together with Item (iii), $\delta'$ satisfies $I$ and $\trans{b}(\ctt)$.
    Let us define $\hat\delta' \dfeq  \delta'\cup(\xi \mapsto V(\delta')\cup\xi' \mapsto L(\delta')$. Since $\hat\delta'$ satisfies the precondition in the premise,
    we have that for any $\hat \delta \in \mathbf{S}\hat\delta'$,
    $\hat\delta$ satisfies $I$ and $V\leq \xi-\xi'$ and $L = \xi'$.
    Hence, $L(\hat\delta) = L(\delta')$. Since $\xi', \xi$ are
    ghost variables, $L(\hat\delta) = L(\delta) = L(\delta')$.
    In conclusion, for any $\delta \in B^{m+1}_{b, S}(\sigma)$,
    the quantity $L(\delta)$ is identical to the quantity
    $L(\delta')$ for some $\delta' \in B^m_{b, S}(\sigma)$.
    Since, $B^0_{b, S}(\sigma) = \{\sigma\}$, we conclude that
    they are all identical to $L(\sigma)$.

    \item

    Suppose any $\delta \in B^{m+1}_{b, S}(\sigma)$ for any $m\in\IN$.
    See that it happens only if there is $\delta' \in B^m_{b,S}(\sigma)$
    such that $\ctt \in \sem{b}\delta'$ and $\delta \in \mathbf{S}\delta'$.
    Together with Item (iii), $\delta'$ satisfies $I$ and $\trans{b}(\ctt)$.
    Consider $\hat\delta' \dfeq  \delta'\cup(\xi \mapsto V(\delta')\cup\xi' \mapsto L(\delta')$
    which satisfies the precondition of the premise.
    Hence, $\delta\cup (\xi \mapsto V(\delta')\cup\xi' \mapsto L(\delta')$
    satisfies the postcondition. Hence,
    $V(\delta) \leq V(\delta') - L(\delta) = V(\delta') - \ell_0$.
    Hence, $\ell_{m+1} \leq \ell_m - \ell_0$.
    \end{enumerate}

    Now, we need to prove the aforementioned statement on $B^m_{b, S}$ and $C^m_{b, S}$:
    \begin{enumerate}
        \item (Base case): Recall that
        $B^0_{b, S}(\sigma) = \{\sigma\} \neq \{\bot\}$ and
        $C^0_{b, S}(\sigma) = \{\}$. Hence, the four conditions
        are all satisfied.

        \item (Induction step):
        Recall
        $B_{b, S}^{n+1}(\sigma):\equiv \bigcup_{\delta \in B_{b, S}^n(\sigma)}
    \bigcup_{\substack{l \in \sem{b}\delta\\\delta'\in\mathbf{S}\delta}}
    \begin{cases}
    \{\delta'\} &\text{if } l = \ctt \land \delta' \neq \bot\\
    \emptyset &\text{otherwise}.
    \end{cases}$.
    Since all $\delta\in B^n_{b, S}(\sigma)$ satisfies $I$,
    $\sem{b}\delta = \{\ctt\}, \{\cff\},$ or $\{\ctt,\cff\}$.
    In the case of $\ctt \in \sem{b}\delta$,
    $\delta$ satisfies the precondition of the premise.
    Hence, for all $\delta' \in \mathbf{S}\delta$, $\delta'$ is not $\bot$ and also satisfies $I$.
    The case of $\sem{b}\delta=\{\cff\}$ is not of interest.

    Recall
    $C_{b, S}^{n+1}(\sigma) :\equiv
    C_{b, S}^n(\sigma)\cup\bigcup_{\delta \in B_{b, S}^n(\sigma)}
    \bigcup_{\substack{l \in \sem{b}\sigma\\ \delta'\in \mathbf{S}\delta}}
    \begin{cases}
    \{\delta\} &\text{if } l = \cff \\
    \emptyset &\text{if } l = \ctt \land \delta' \neq \bot\\
    \{\bot\} &\text{otherwise}.
    \end{cases}$.

    Since all $\gamma \in C^n_{b, S}(\sigma)$ satisfies $I$ and $\trans{b}(\cff)$,
    we only need to care the rightmost part of the construction.
    Since all $\delta \in B^n_{b, S}(\sigma)$ satisfies $I$, by the
    side-condition, $\cuu$ and $\bot$ are not in  $\sem{b}\delta$.
    The $\delta$ is added to $C^{n+1}(\sigma)$ only if $\cff \in \sem{b}\delta$.
    Therefore, $\delta$ satisfies both $I$ and $\trans{b}(\cff)$.

    Also, in the case of $\ctt\in \sem{b}\delta$, since
    $\delta$ satisfies the precondition in the premise,
    $\bot \not\in \mathbf{S}\delta$. Therefore, $\bot\not\in C^{n+1}_{b, S}$.
    \qedhere
    \end{enumerate}
\end{enumerate}
\end{proof}

\end{document}